\documentclass{amsart}

\usepackage[all]{xy}
\usepackage{amsmath}
\usepackage{amsfonts,amssymb}

\newtheorem{satz}{Proposition}[section]
\newtheorem{Satz}[satz]{Theorem}
\newtheorem{definition}{Definition}[section]
\newtheorem{lemma}[satz]{Lemma}

\newtheorem{lemdef}[satz]{Lemma/Definition}
\newtheorem{kor}[satz]{Corollary}

\newtheorem{remark}[satz]{Remark}

\def\bem{{\bf Remark: }}

\newcommand{\tensor}{\otimes}

\newcommand{\map}[1]{\stackrel{#1}{\longrightarrow}}
\newcommand{\incl}[1]{\stackrel{#1}{\hookrightarrow}}
\newcommand{\eqweil}[1]{\stackrel{#1}{=}}
\newcommand{\congweil}[1]{\stackrel{{\textrm{\tiny #1}}}{\cong}}

\newcommand{\un}[1]{\ensuremath{\protect\underline{#1}}}

\DeclareMathOperator{\Coh}{Coh}

\def\1halb{\frac{1}{2}}

\def\tto{\twoheadrightarrow}
\def\pprime{{\prime\prime}}

\def\Qbar{{\overline{\bQ}}}

\def\Fbar{{\overline{\bF}}}

\def\GL{\textsf{GL}}
\def\SL{\textsf{SL}}
\def\Mat{\textsf{Mat}}

\def\Iw{\mathsf{Iw}}

\def\St{\mathrm{St}}

\def\Fib{\text{Fib}}
\def\Fibre{\text{Fibre}}

\DeclareMathOperator{\Hom}{Hom}
\DeclareMathOperator{\Spec}{Spec}
\DeclareMathOperator{\Ext}{Ext}
\def\uExt{\underline{\textrm{Ext}}}
\DeclareMathOperator{\Ker}{Ker}
\DeclareMathOperator{\Aut}{Aut}
\DeclareMathOperator{\Flag}{flag}

\DeclareMathOperator{\length}{length}
\DeclareMathOperator{\rank}{rank}
\DeclareMathOperator{\diag}{diag}
\DeclareMathOperator{\sign}{sign}
\DeclareMathOperator{\vol}{vol}
\DeclareMathOperator{\supp}{supp}
\DeclareMathOperator{\res}{Res}
\DeclareMathOperator{\Trace}{trace}
\DeclareMathOperator{\tr}{tr}
\DeclareMathOperator{\Frob}{Frob}
\def\proj{\mathit{proj}}
\DeclareMathOperator{\torsion}{torsion}

\def\gr{\mathrm{gr}}
\def\forget{\mathit{forget}}
\def\pr{\mathit{pr}}
\def\bl{\mathit{bl}}

\def\ext{\mathit{ext}}
\def\qext{\mathit{qext}}
\def\quot{\mathit{quot}}
\def\Bl{\mathrm{Bl}}

\def\Four{\mathcal{F}\mathrm{our}}

\DeclareMathOperator{\Hecke}{Hecke}
\DeclareMathOperator{\Bun}{Bun}
\def\univ{\textrm{\tiny univ}}

\def\OmegaExt{\Omega\text{--Ext}}
\def\OmegaP{\Omega\text{--Pl\"ucker}}

\def\sxymat{\xymatrix@C=0.8ex@R=0.8ex}
\def\grp{$\xymatrix{ R\times_{X}R  \ar[r]^-{\mu} & R \ar@<1ex>[r]^-{s}\ar@<-1ex>[r]_-{t} & X}$}
\def\dar{\ar@<-0.5ex>[r]\ar@<0.5ex>[r]}
\def\tar{\ar[r]\ar@<1ex>[r]\ar@<-1ex>[r]}
\newcommand{\dmap}[2]{\ar@<-0.5ex>[r]_-{#2}\ar@<0.5ex>[r]^-{#1}}
\newcommand{\dotarrow}[2]{\xymatrix{{#1}\ar@{..>}[r]&{#2}}} 
\def\cart{\ar@{}[dr]|{\square}}

\def\cC{\mathcal{C}}

\def\cE{\mathcal{E}}
\def\cF{\mathcal{F}}
\def\cG{\mathcal{G}}
\def\cH{\mathcal{H}}
\def\cI{\mathcal{I}}
\def\cJ{\mathcal{J}}

\def\cL{\ensuremath{\mathcal{L}}}
\def\cM{\mathcal{M}}

\def\cO{\mathcal{O}}

\def\cQ{\mathcal{Q}}

\def\cT{\mathcal{T}}

\def\cp{\mathfrak{p}}

\def\bA{{\mathbb A}}

\def\bD{{\mathbb D}}

\def\bG{{\mathbb G}}
\def\bP{{\mathbb P}}

\def\bV{{\mathbb V}}

\def\bN{{\mathbb N}}
\def\bZ{{\mathbb Z}}
\def\bQ{{\mathbb Q}}
\def\bR{{\mathbf R}}

\def\bF{{\mathbb F}}

\def\B{\textsf{B}}
\def\N{\textsf{N}}
\def\P{\textsf{P}}

\def\sfA{\mathsf{A}}

\def\sfE{\mathsf{E}}
\def\sfF{\mathsf{F}}
\def\sfK{\mathsf{K}}
\def\sfL{\mathsf{L}}
\def\sfN{\mathsf{N}}

\def\sfW{\mathsf{W}}

\begin{document}

\title[Some Hecke eigensheaves]{Coherent sheaves with parabolic structure and 
Construction of Hecke eigensheaves for some ramified local systems}     
\author{Jochen Heinloth}
\address{Mathematisches Institut der Universit\"at Bonn, Beringstr. 1, D--53111 Bonn, Germany}
\email{heinloth@math.uni-bonn.de}
\begin{abstract}
The aim of these notes is to generalize Laumon's construction \cite{Laumon_correspondance} of 
automorphic sheaves corresponding to local systems on a smooth, projective curve $C$ to  the case 
of local systems with indecomposable unipotent ramification at a finite set of points. 
To this end we need an extension of the notion of parabolic structure on vector bundles to coherent sheaves.
Once we have defined this, a lot of arguments from the article "On the geometric
Langlands conjecture" by Frenkel, Gaitsgory and Vilonen \cite{FGV3} carry over to our situation.
We show that our sheaves descend to the moduli space of parabolic bundles if the rank is $\leq 3$ 
and that the general case can be deduced form a generalization of the vanishing conjecture of 
\cite{FGV3}.
\end{abstract}
\maketitle


\UseComputerModernTips
\tableofcontents

\section*{Introduction}

Before explaining the main result (Theorem \ref{ergebnis}) of this article in more detail, I would like to recall the setting of the geometric Langlands correspondence as in \cite{Laumon_Bourbaki}.

Let $C$ be a smooth projective curve over a finite field $\bF_q$.
(As pointed out in \cite{FGV3} and \cite{Laumon_Bourbaki}, a lot of the arguments carry over to the case when $C$ is defined over the complex numbers.) 

In this situation the Langlands correspondence --- as proven by Lafforgue \cite{lafforgue} --- provides a bijection between irreducible $\ell$--adic local systems defined on some open subset $U\subset C$ and certain irreducible representations of $\GL_n(\bA)$ contained in the space $\cC^\infty\big(\GL_n(k(C))\backslash \GL_n(\bA)\big)$ called the space of automorphic functions. Here we denoted by $\bA:=\prod_{x\in C}^\prime K_x$ the ring of adeles of the function field $k(C)$ of $C$, and by $\cC^\infty\big(\GL_n(k(C))\backslash \GL_n(\bA)\big)$ the space of functions (with values in $\Qbar_\ell$) that are right invariant under some compact open subgroup of $\GL_n(\bA)$ (for notations see Section \ref{notations}). More precisely it is known (see e.g. \cite{Laumon_Cohomology_2}) that for any representation $\pi_\sfE$ corresponding to some local system $\sfE$ there is a compact open subgroup $\sfK$ such that $\pi_{\sfE}$ contains  a (up to scalar) unique  $\sfK$--invariant function $A_\sfE$. Further, this compact subgroup is determined by the ramification of $\sfE$. 
Finally, note that the group $\GL_n(\bA)$ does not act on the $\sfK$-invariant functions, but the algebra of $\sfK$-bi-invariant functions acts on these by convolution. This is the action of the $\sfK$-Hecke algebra. The function $A_\sfE$ is an eigenvector for this action, and it is determined by this condition. 

Drinfeld noted \cite{Drinfeld_fundamental_group} that this correspondence might have a geometric interpretation. First consider the case $\sfK=\GL_n(\cO)$. Weil explained that the double quotient $\GL_n(k(C))\backslash \GL_n(\bA) / \GL_n(\cO)$ can be identified with the set of isomorphism classes of vector bundles on $C$ (choose a trivialisation at all local rings of $C$ and at the generic point of $C$, the transition functions give an adele):
$$ \Bun_n(\bF_q) = \GL_n(k(C))\backslash \GL_n(\bA) / \GL_n(\cO).$$
Furthermore, Grothendieck explained that any complex $\sfA$ of $\ell$--adic sheaves on a scheme $X_{/\bF_q}$ gives rise to a function on the set of its points by
\begin{align*}
\Trace : \bD^b(X) &\to \prod_{n\in \bN} \text{Funct}(X(\bF_{q^n}))\\
  \sfA &\mapsto \text{tr}_\sfA(x) := \Trace(\Frob_{\bF_{q^n}},\sfA|_x)
\end{align*}
and that a perverse complex is determinded by this function (\cite{Laumon_transformation}).

Thus, Drinfeld expected that the above $A_\sfE$ should be of the form $\text{tr}_{\sfA_\sfE}$ for some irreducible perverse sheaf $\sfA_\sfE$ on the moduli space of vector bundles on $C$. He proved this for unramified local systems of rank $2$. Later Laumon (\cite{Laumon_correspondance}) gave a conjectural construction of $\sfA_\sfE$ for local systems of arbitrary rank, and recently Frenkel, Gaitsgory and Vilonen (\cite{FGV3},\cite{Gaitsgory}) proved that by Laumon's construction one indeed obtains a sheaf $\sfA_\sfE$.

Moreover, the action of the Hecke algebra also has a geometric interpretation in this case. Consider for example the characteristic function of the double coset $\GL_n(\cO_x)\left( {1 \atop 0}{0 \atop \pi_x} \right)\GL_n(\cO_x)$, where $\pi_x$ is a local parameter at some point $x\in C$. For a vector bundle $\cE$ the multiplication of the corresponding adele by an element of this set produces a subbundle $\cE^\prime \subset \cE$ such that the cokernel is $k(x)$. Further, every such subbundle can be obtained in this way. Drinfeld therefore considered the stack $\Hecke^1$ classifying pairs of bundles $\cE^\prime\subset \cE$ such that the cokernel has length $1$, i.e. $\deg(\cE^\prime)=\deg(\cE)-1 =: d-1$. This has forgetful maps
$$\xymatrix{ & \Hecke^1 \ar[dl]_{\pr_{\text{big}}} \ar[dr]^{\pr_{\text{small}}\times \quot} \\\Bun^d & & \Bun^{d-1} \times C}$$
With this definition the sheaf $\sfA_\sfE$ has the additional property that $$\bR (\pr_{\text{small}}\times \quot)_! \pr_{\text{big}}^* \sfA_\sfE \cong \sfA_\sfE \boxtimes \sfE,$$
and a similar definition works for more general Hecke stacks. One says that $\sfA_\sfE$ is a {\em Hecke eigensheaf}. 

Drinfeld also proved an analogous result for local systems of rank $2$ with unipotent ramification at a finite set of points $S\subset C(\bF_q)$ (see \cite{Drinfeld_galois_group}), this time producing a complex $\sfA_\sfE$ on the moduli space of vector bundles of rank $2$ with parabolic structure at $S$. The purpose of this article is to generalize this result. 

We will start with a local system $\sfE$ with unipotent ramification at a finite set of points $S\subset C(\bF_q)$, and we further have to assume that the ramification group at these points acts indecomposably, i.e. that the sheaf $j_*\sfE$ (where $j:C-S \to C$) has one-dimensional stalks at all points $p\in S$. 
This additional condition is the reason why for the moment we can only prove our main theorem for local systems of rank $\leq 3$.

In this case the corresponding automorphic function should be defined on the space $\GL_n(k(C)) \backslash
\GL_n(\bA) / \sfK_S$, where $\sfK_S= \prod_{x\in C-S} \GL_n(\cO_x) \times \prod_{x \in S} \Iw_x $ and $\Iw_x \subset \GL_n(\cO_x)$ 
is the subgroup of matrices which are upper triangular mod $x$.
As before we can interpret this set as vector bundles with the additional structure of a complete
flag of subspaces of the stalks at all points in $S$:
$$ \Bun_{n,S}(T) := \langle \big(\cE,(V_{i,p})_{i=1,\dots,n \atop p\in S}\big)\, |\, \cE\in \Bun_n; 0\subset V_{1,p} \subset \dots \subset V_{n,p} = \cE\tensor k(p) \rangle$$
This is usually called the stack of vector bundles with (quasi-)parabolic structure. Note that this can
also be described as:
$$ \Bun_{n,S}(T) := \langle \big(\cE,(\cE^{i,p})_{i=1,\dots,n \atop p\in S}\big)\, | \, \cE\in \Bun_n; \cE\subset \cE^{1,p} \subset \dots \subset \cE^{n,p} = \cE(p) \rangle$$
which has a simple generalization to coherent sheaves: one only has to replace ``$\subset$'' by 
arbitrary maps ``$\to$'' and to add the condition that the induced maps $\cE^{i,p}\to\cE^{i,p}(p)$ are the natural ones. This reformulation made our construction possible.

The first step of our construction is to recall that in principle a candidate for the automorphic function $A_\sfE$ is known, but we do not know of an explicit calculation of this function. Therefore, we have to prove an explicit formula (Proposition \ref{formel}). This motivates a generalization of Laumon's construction, and --- as a by--product of the notion of parabolic torsion sheaf --- we get a geometric interpretation of some Hecke operators for the group $\sfK$, i.e. of the Iwahori--Hecke algebra. Our main result is then the following:

\noindent {\bf Theorem \ref{ergebnis}. }{\em  For any irreducible local system $\sfE$ of rank $n \leq 3$ on $C-S$ with indecomposable unipotent ramification at $S$ there is an irreducible perverse sheaf $\sfA_\sfE$ on $\Bun_{n,S}$ which is an eigensheaf for the Iwahori--Hecke algebra.}

The strategy of the proof is the same as in \cite{FGV3}, using parabolic sheaves instead of coherent sheaves, but some additional problems arise from the ramification of $\sfE$. We reduce the theorem to an analogue (Prposition \ref{verschwinden}) of the vanishing conjecture of loc. cit. In particular, we show that the above theorem would follow for local systems of general rank if this analogue held in general.

The structure of the article is as follows. We start with the calculation of the Whittaker function for the Steinberg representation given in first section. This is an elementary calculation which served as motivation for our construction.

In the second section we introduce the notion of a coherent sheaf with parabolic structure and prove the results needed to give an analogue of Laumon's ``fundamental diagram'' and of Laumon's Whittaker sheaf $\cL_\sfE^d$. We then define two candidates for an automorphic sheaf. At the end of this section we define the geometric Hecke operators corresponding to operators of the Iwahori-Hecke algebra which are needed to give a precise formulation of our main Theorem \ref{ergebnis}.

After this short exposition of our results we try to clarify the notion of parabolic sheaves in Section 3. We explain the general structure of parabolic torsion sheaves. Further, we give an explicit description of the corresponding moduli stack, and finally we note some semicontinuity results. We then use these basic results to prove some properties of the Whittaker sheaf $\cL_\sfE^d$ (Section 4). Here we give a substitute for the Springer resolution in the case of parabolic sheaves which can be used to calculate this sheaf, and we prove a Hecke property of $\cL_\sfE^d$. The problem arising in the proof of these results is that in our situation the above resolution is not small and the ramification of $\sfE$ also generates additional cohomology. By simultaneously proving the Hecke property and the fact that $\cL_\sfE^d$ can be calculated via the resolution we see that the two effects cancel out.

In the fifth section we then compare the geometric construction of Section 2 with the calculation of the Whittaker function. The key idea here is to define an analogue of Drinfeld's compactification as given in \cite{FGV3}. However, we can not copy the proofs of loc. cit., which use results on the affine Grassmannian for which we do not know the corresponding statements for the affine flag manifold. Instead, we give an elementary proof of a much weaker result, sufficient for our purpose. 

With these results available we can follow the strategy of \cite{FGV3} again and apply Lafforgue's result to deduce the existence of a Hecke eigensheaf on the moduli space of parabolic vector bundles whenever we know that the two candidates constructed coincide. This is the content of Section 6.

In the last two sections we then prove a generalization of the vanishing theorem of 
\cite{FGV3} for local systems of rank $\leq 3$ and deduce the assumption needed to prove our theorem in Section 6. This is again very similar to the arguments in loc.cit., however we have to take care of the Iwahori-Hecke operators, for example we have to prove that some of them are central elements of the algebra (see Lemma \ref{kommutiert}).
\vspace{3ex}

\noindent {\bf Acknowledgments.} 
First of all I would like to thank my advisor G. Harder for teaching mathematics to me for 
such a long time. Furthermore I would like to thank G. Laumon for his help and encouragement during a wonderful visit to the Universit\'e Paris-Sud. Without his help I would never have been able to start this project. I would like to thank the mathematics department of 
Paris Sud for the kind hospitality, and the DAAD and the Graduiertenkolleg in Bonn for giving me the possibility to stay there. 



\setcounter{section}{-1}

\section{Notations and preliminary remarks}\label{notations}

We want to fix some notations used throughout this article.

\subsection{The curve and its rings}

We fix a smooth projective curve $C$ defined over a finite field $k=\bF_q$ and denote by: \begin{itemize}
\item $k(C)$ the field of rational functions on $C$.
\item $\cO_p$ (resp. $\widehat{\cO}_p$) the local ring (resp. the complete local ring) at a point $p\in C$.
\item $K_p:= \text{Quot}(\widehat{\cO}_p).$
\item $\bA:= \prod^\prime_{p\in C} K_p$ the ring of adeles of $k(C)$.
\item $\cO:= \prod_{p\in C} \widehat{\cO}_p$.
\item $\Omega:=\Omega_{C/k}$ the sheaf of differentials on $C$. 
\end{itemize} 

\subsection{Groups}\label{GLomega}

\begin{itemize}
\item We note by $\GL_n$ the algebraic group of invertible $n\times n$ matrices.
\item $\B_n\subset \GL_n$ the group of upper-triangular matrices.
\item $\N_n\subset \B_n$ the group of unipotent upper triangular matrices. 
\item $\P_1\subset \GL_n$ the subgroup fixing the subspace spanned by the first $n-1$ base vectors and acting trivially on the quotient by this subspace, i.e. $\P_1(R)= \{ \left( \begin{array}{cc} A & v \\ 0 & 1 \end{array} \right) \; |\; A\in GL_{n-1}(R), v \in R^{n-1} \}. $
\item $\Iw\subset \GL_n(\widehat{\cO}_p)$ the group of matrices which are upper triangular mod $p$.
\end{itemize} 

We will further fix a non-trivial additive character $\psi: \bF_q \to \Qbar_l^*$.

Choosing a meromorphic differential form $\omega$  this defines 
$$\Psi: \N_n(k(C)) \backslash \N_n(\bA) \to \Qbar_l^*$$ 
$$\Psi((U_p)_{p\in C}):= \prod_{p\in C} \psi \big(\Trace_{k(p)/\bF_q}(\res_p(\sum_{i=1}^{n-1} u_{p,i,i+1})\omega)\big),$$ 
where $u_{p,i,i+1}$ is the $i$-th entry of the first upper diagonal of the matrix $U_p$.

To avoid the choice of a meromorphic differential form we will (as in \cite{FGKV}) often replace the group $\GL_n \times C/C$ by the group $\GL_n^\Omega := \Aut(\oplus_{i=1}^n \Omega^{\tensor n-i})$.  More precisely, $\GL_n\times C =\Aut(\cO^{\oplus n})$ is the automorphism group of the trivial vector bundle over $C$, since for any ring $R$ the automorphisms of the trivial rank $n-$bundle over $\Spec(R)$ are the same as elements of $\GL_n(R)$. In the same way points of $\GL_n^\Omega$ are invertible matrices in which the $(i,j)$-th entry is a section of $\Hom(\Omega^{i-1},\Omega^{j-1})\cong\Omega^{j-i}$. In particular, the choice of a meromorphic differential $\omega$ induces a group isomorphism $\GL_n(\bA) \map{\sim} \GL_n^{\Omega}(\bA)$. Under this isomorphism $\Psi$ is induced from the canonical morphism $\GL_n^\Omega(\bA) \map{\sum \res} \bF_q \map{\psi} \Qbar_l^*$ defined by the sum of the residues of the upper diagonal entries.

\subsection{Fourier transform}

For the additive character $\psi:\bF_q \to \Qbar_l^*$ chosen above we denote by $\sfL_\psi$ the Artin-Schreier sheaf on $\bA^1$: 
Let  $AS:\bA^1 \map{x\mapsto x^q-x} \bA^1$  be the Artin-Schreier covering with structure group $\bF_q$, then $\Psi$ is the $\psi$-isotypic component of $AS_* \Qbar_l$. 

For a vector bundle $E\map{p} X$ of rank $n$ on a scheme (or algebraic stack) denote by $E^\vee \map{p^\vee} X$ the dual bundle and by $<,>: E\times_X E^\vee \to \bA^1$ the contraction. The Fourier transform defined in \cite{Laumon_transformation} is given by
\begin{eqnarray*}
 \Four : D^b(E) & \to & D^b(E^\vee) \\
 \sfK & \mapsto & \bR p^\vee_! (p^* \sfK \tensor <,>^*\sfL_\psi)[n]. 
\end{eqnarray*}

\subsection{The trace function of a complex}

For a complex $\sfK$ of $\Qbar_\ell$-adic sheaves on a scheme (or algebraic stack) $X$ we denote by $\tr_\sfK$ the function:
\begin{align*}
\tr_\sfK : \prod_{n>0} X(\bF_{q^n}) & \to \Qbar_\ell\\
 x & \mapsto tr_\sfK(x):=\Trace(\Frob_x,\sfK|_x).
\end{align*}

\subsection{Algebraic stacks}

For the general theory of algebraic stacks we refer to the book of Laumon and Moret-Bailly \cite{champs-algebriques}. In particular, an algebraic stack will be a stack that admits a smooth representable covering by a scheme.

We will view stacks as sheaves of categories for the fppf-topology.
Thus to define a stack $\cM$ we usually give the category
of $T$-valued points of $\cM$ and denote this as:
$$ \cM(T):= \langle \text{\em objects} \rangle, $$
where we use the brackets $\langle \qquad \rangle$ instead of $\{ \qquad \}$ to denote the category
of {\em objects} in which the only morphisms are isomorphisms of the {\em objects}.

Sometimes it is easier to give the $T$-valued points of a stack only for affine
schemes $T$ over the given base, which is equivalent to the data for all schemes by the descent condition for stacks. This point of view is used as definition  in loc. cit.

\subsection{Some remarks on generalized vector bundles}

Recall that for a flat algebraic group $G$ acting on a scheme $X$ there is a 
quotient stack $[X/G]$ classifying principal $G-$bundles together with a 
$G-$equivariant morphism to $X$. In this section we will be concerned with the particular case of a homomorphism of vector 
bundles $E_0\map{\phi} E_1$ and take $G:=E_0$ acting additively on $X:=E_1$: 

\begin{definition}(\cite{inc}) Let $E_0 \map{\phi} E_1$ be a homomorphism of vector bundles on a scheme (or an algebraic stack) $X$. Then the 
quotient stack $[E_1/E_0]$ is called a {\em generalized vector bundle} over $X$
\end{definition} 

\begin{lemma}\label{genvb}
Let $E_0 \map{\phi} E_1$ be a homomorphism of vector bundles on some scheme (or algebraic stack) $X$. The stack $[E_1/E_0]$ can be described as follows: 

For any affine scheme  $T=Spec(A)\map{f} X$ over $X$:                                      
$$ [E_1/E_0](T) = \left\langle \begin{array}{l}\text{objects}=\{s\in H^0(T,f^*E_1)\}  \textrm{ and for } s,t\in H^0(T,f^*E_1) \\ 
\Hom(s,t)=\{h \in H^0(T,f^*E_0) | s+\phi(h)=t \}\end{array} \right\rangle $$
Moreover, any quasi-isomorphism of such complexes gives rise to an equivalence of the corresponding stacks, thus the stack $[E_1/E_0]$ depends only on the class of the complex $E_0\to E_1$ in the derived category of coherent sheaves
on $X$.
\end{lemma}

\noindent{\bf Example:} Let $C\map{p}X$ be a smooth projective curve over some noetherian base scheme $X$, and let $\cF_1,\cF_2$ be coherent sheaves on $C$, flat
over $X$. By [EGAIII] the complex $\bR p_* (\cH om (\cF_1,\cF_2))$ can be represented by a homomorphism of vector bundles $\cE_0 \to \cE_1$ on $X$. 
By abuse of notation we denote by                                       
$\uExt(\cF_1,\cF_2)$ the corresponding generalized vector bundle on $X$. 

Note that this is well defined by the above lemma. The description of the categories of sections given in the lemma 
tells us that this stack classifies extensions $0\to \cF_2 \to  \cF \to \cF_1 \to 0$,
i.e. for any $T\map{f}X$:
$$ \uExt(\cF_1,\cF_2)(T) = \langle 0 \to f^*\cF_2 \to \cF \to f^*\cF_1 \to 0 \rangle. $$
\noindent{\bf Proof } (of Lemma \ref{genvb}): First note that the claimed description of $[E_1/E_0]$ defines a stack: 

1. We can glue morphisms, because sections of $E_0$ form a sheaf.

2. Any descent datum of objects is effective (i.e. we can glue objects): 
Let $U_i$ be an affine covering of the affine scheme $T$. A descent datum for this covering is a collection of objects $s_i\in \Gamma(U_i,E_1)$ together with morphisms $h_{ij} \in \Gamma(U_{ij},E_0)$
such that $s_i|_{U_{ij}} + \phi(h_{ij}) = s_j|_{U_{ij}}$ and  $h_{ik}|_{U_{ijk}}=h_{jk}|_{U_{ijk}}+h_{ij}|_{U_{ijk}}$. 

This implies that $h_{ij}$ is a 1-cocycle, and since $T$ is affine it must
be a coboundary, i.e. we can find $h_i\in H^0(U_i,E_0)$  with $h_i-h_j = h_{ij}$ on $U_{ij}$. Therefore we may define $s_{i}^\prime:= s_i - h_i$,
and this collection of sections glues to give $s\in H^0(T,E_1)$ with $s|_{U_i}= s_i^\prime$.

Thus we may define a morphism of stacks 
$$\left\langle \begin{array}{l}\text{objects}=\{s\in H^0(T,E_1|_T)\}  \textrm{ and for } s,t\in H^0(T,E_1) \\ 
\Hom(s,t)=\{h \in H^0(T,E_0) | s+\phi(h)=t \}\end{array} \right\rangle \to [E_1/E_0](T)$$ 
mapping a section $T\to E_1$ to the composition $T\to E_1 \to [E_1/E_0]$.
 
Since $H^1(T,E_0)=0$ for an affine $T$ any $s\in [E_1/E_0](T)$ is isomorphic to some $s^\prime \in H^0(T,E_1)$ and by definition any 
morphism between two elements $s,t$ in the image of this functor is given by a section of $H^0(T,E_0)$. Thus the morphism is an equivalence of stacks. 

The above description of the stack $[E_1/E_0]$ also shows that a quasi-isomorphism of complexes induces an equivalence of the categories of points of the corresponding stacks. 
\hfill $\square_{\textrm{\tiny Lemma}}$

\begin{lemma}\label{exakte_sequenzen}
Let $\xymatrix@R=1.4ex{ 0 \ar[r] & E_0^\prime \ar[d]^{\phi^\prime}\ar@{^(->}[r]^{i_0} & E_0 \ar@{->>}[r]^{p_0} \ar[d]^{\phi}
 & E_0^\pprime \ar[r] \ar[d]^{\phi^\pprime} & 0 \\
0 \ar[r] & E_1^\prime \ar@{^(->}[r]^{i_1} & E_1 \ar@{->>}[r]^{p_1}  & E_1^\pprime \ar[r] & 0 }$ 
be an exact sequence of (2 term-)complexes of vector bundles on some (quasi-separated) scheme $X$. Denote by $[E^\prime_1/E^\prime_0]\map{i}[E_1/E_0]\map{p} [E^\pprime_1/E^\pprime_0]$ the induced morphisms of the generalized bundles,  and let $s^\pprime: X \to [E_1^\pprime/E_0^\pprime]$
be a section. 

Then locally over $X$ the stack $p^{-1}(s^\pprime)= [E_1/E_0]\times_{[E^\pprime_1/E^\pprime_0]}X$  
is isomorphic to $[E^\prime_1/ E^\prime_0]$. More precisely such an isomorphism exists over any $U\to X$ such that 
there is a lift $s_1\in \Gamma(U,E_1)$ with $p(s_1)\cong s$.
\end{lemma} 

\noindent{\bf Remark:} We might state the above as ``$p^{-1}(s^\pprime)$ is a principal 
homogeneous space for $[E_1^\prime/E_0^\prime]$''. More generally, we will call a morphism of
stacks a {\em generalized affine space bundle} if it can be factored into a sequence of
maps each of them locally (over the target space) isomorphic to a generalized vector bundle. 

\noindent{\bf Proof:} We may assume that $X=U$, such that there exists $s_1\in H^0(U,E_1)$ with $p(s_1)=s^\pprime$ (e.g. we can take $U$ affine). 

Using the previous lemma, we find that $p^{-1}(s^\pprime)(T)$ is the category with:\begin{eqnarray*} \text{objects}= \{(s,h^\pprime)\in \Gamma(T,E_1)\times \Gamma(T,E_0^\pprime) \; | \;  p_1(s)+\phi^\pprime(h)=s_1\}&&\\
\Hom((s,h^\pprime),(t,g^\pprime))= \{ h \in \Gamma (T,E_0) \; | \; s+\phi(h)=t \textrm{ and } p_0(h_0)=h^\pprime-g^\pprime  \}. 
\end{eqnarray*}
Thus we define:\begin{eqnarray*}
[E_1^\prime/ E_0^\prime] & \to & p_1^{-1}(s_1) \\
         H^{0}(T,E_1^\prime) \ni  s^\prime & \mapsto & (i_1(s^\prime) + {s}_1,0)\\
         H^{0}(T,E_0^\prime) \ni h^\prime & \mapsto & i_0(h^\prime). 
            \end{eqnarray*}         
This is essentially surjective, since for affine $T$ and any $h^\pprime\in H^0(T,E_0^\pprime)$ there is an $h\in H^0(T,E_0)$ with $p_0(h)=h^\pprime$, and therefore any $(s,h^\pprime)\cong(s-\phi(h),0)$.  
Morphisms of two objects in the image of the above map are given by
$H^0(T,E_0^\prime)=\Ker(H^0(T,E_0)\to H^0(T,E^\pprime_0))$,
therefore this is an equivalence of categories.
                       \hfill$\square_{\textrm{\tiny Lemma}}$

{\bf Application:} We will apply this lemma in the following situation: Consider the morphism of stacks classifying
diagrams (with exact lines and columns) of torsion sheaves on a curve $C$:
$$ 
\left\langle\vcenter{\xymatrix@C=2ex@R=2ex{ {\cT}_1^\prime \ar@{^(->}[d]  & & \\ {\cT}_1 \ar@{->>}[d]\ar@{^(->}[r] & {\cT}_2 \ar@{->>}[r]\ar@{->>}[d] & {\cT}_3^\pprime \ar@{=}[d]\\ 
{\cT}^{\pprime}_1 \ar@{^(->}[r] & {\cT}^{\pprime}_2 \ar@{->>}[r] & {\cT}^{\pprime}_3}}\right\rangle
 \map{\forget_{\cT_2}}
\left\langle \vcenter{\xymatrix@C=2ex@R=2ex{ {\cT}_1^\prime \ar@{^(->}[d]  & & \\ {\cT}_1 \ar@{->>}[d] & & \\ 
{\cT}^{\pprime}_1 \ar@{^(->}[r] & {\cT}^{\pprime}_2 \ar@{->>}[r] & {\cT}^{\pprime}_3}}\right\rangle,
$$ where the degree of each torsion sheaf is fixed.

On the right hand stack the exact triangle of complexes
$$\bR \Hom(\cT_3^\pprime, \cT_1^\prime) \to \bR \Hom(\cT_3^\pprime, \cT_1) \map{p} \bR \Hom(\cT_3^\pprime, \cT_1^\pprime)$$
can be represented by an exact sequence of 2-term complexes of vector bundles.
There is a canonical $s^\pprime$ of $\bR \Hom(\cT_3^\pprime, \cT_1^\pprime)$ 
given by the extension in the lower line, and the projection map from $p^{-1}(s^\pprime)$ to the base stack is the map $\forget_{\cT_2}$.

Thus, by the above lemma, we see that the fibres of this morphism are isomorphic
to the stack $\uExt(\cT_3^\pprime,\cT_1^\prime)$. These stacks are generalized
affine spaces, in particular the \'etale cohomology of the fibres is one-dimensional.  

\subsection{A lemma used more than once...}

The following general lemma is stated in \cite{FGV3}, a similar calculation is done in \cite{brylinski}.
I would like to thank Sergey Lysenko for explanations about this:

\begin{lemma}\label{brylinski-lemma}                         
Let $\cE\map{p} X$ be a (generalized) vector bundle, and denote by $s_0:X\to \cE$ the zero-section of $\cE$.
Let further $\sfK\in D^b_{\text{\'et}}(\cE)$ be a complex of \'etale sheaves on $\cE$ such that the restriction of $\sfK$ to the complement of the zero-section descends to the projective bundle $\bP(\cE)$ (e.g. a $\bG_m$-invariant complex of sheaves on $\cE$).
Then $$\bR p_* \sfK = s_0^* \sfK .$$
\end{lemma}
\noindent{\bf Proof:} We may assume that $\cE$ is a vector bundle, since for a generalized vector bundle $[\cE_1/\cE_0]$ the functor $\bR p_*$ is defined via an acyclic representable covering of the bundle, i.e. by definition we may replace $[\cE_1/\cE_0]$ by $\cE_1$. 

Let $j: \cE^\circ:=\cE- s_0(X) \hookrightarrow \cE$ be the inclusion.  Then we have an exact triangle
$$  \to j_!j^* \sfK \to \sfK \to s_{0,*}s_0^* \sfK \map{[1]}.  $$

For the first term $j_!j^* \sfK$ we have to prove that $\bR p_* j_!j^* \sfK=0$.
If we can show this we are done, since the lemma is true for the last term, and the right hand map then gives the claimed isomorphism.

Write $\sfK^\circ:= j^* \sfK$. Then by assumption $\sfK^\circ \cong \proj^*(\overline{\sfK})$, where $\proj: \cE^\circ \to \bP(\cE)$ is the projection 
to the projectivized bundle and $\overline{\sfK}$ is a sheaf on $\bP(\cE)$.
To get a relation between $\cE$ and $\bP(\cE)$ , blow up the zero-section of $\cE$, and denote the blow up by $\Bl_{s_0}(\cE)$:
$$\xymatrix{
{\Bl}_{s_0}({\cE}) \ar[d]_{bl}\ar[drr]^{\pr_{\bP({\cE})}} & & \\
{\cE}\ar[d]_p & {\cE}^\circ \ar@{_(->}[l]^{j}\ar@{_(->}[ul]^{\tilde{j}} \ar[r]_-{\proj} & {\bP}({\cE})\ar[dll]^-{\overline{p}}\\
X & & 
}$$

Note that $\Bl_{s_0}(\cE)\map{\pr_{\bP(\cE)}} \bP(\cE)$ is the line bundle $\cO(-1)$ over $\bP(\cE)$. Let $s_{\bP(\cE)}:\bP(\cE) \to Bl_{s_0}(\cE)$
be the zero-section (i.e. the inclusion of the special fibre of the blow-up).

Since 
\begin{eqnarray*}
 j_!\, \proj^* (\overline{\sfK}) & = & \bR \bl_! \,\tilde{j}_! \,\proj^* (\overline{\sfK})\\
                   & \eqweil{bl \text{ projective}} & \bR \bl_* \, \tilde{j}_! \, \proj^* (\overline{\sfK}),
\end{eqnarray*}
we need to show that $\bR (p\circ \bl)_* ( \tilde{j}_!\, \proj^* (\overline{\sfK}))=0$. But this is easy, since --- as before --- there is an exact triangle on $Bl(\cE)$
$$  \to \tilde{j}_! \proj^* \overline{\sfK} \to \pr_{\bP(\cE)}^* \overline{\sfK} \to s_{\bP(\cE),*}\overline{\sfK} \to, $$ 
and the natural map induces
$$ \bR \pr_{\bP(\cE),*} \pr_{\bP(\cE)}^* \overline{\sfK} \congweil{proj. formula} \overline{\sfK} \cong \bR \pr_{\bP(\cE),*} s_{\bP(\cE),*} \overline{\sfK}.$$

Thus $\bR \pr_{\bP(\cE),*} \tilde{j}_! \proj^* \overline{\sfK} =0$, and therefore $\bR (p\circ \bl)_*  \tilde{j}_! \proj^* \overline{\sfK} =0$, since $p\circ bl$ factors through $\bP(\cE)$. 
\hfill $\square$


\section{The Whittaker function for the Steinberg representation}

As indicated in the introduction, for any local system $\sfE$ of rank $n$ on $C-S$ with indecomposable unipotent ramification at points in $S$ there is a particular function $f_\sfE$ on $\GL_n(\bA)$ which one expects to span the automorphic representation corresponding to $\sfE$.

In this section we will give a formula for this function, more precisely we will give an explicit formula for a function $W_\sfE$ from which $f_\sfE$ may be obtained by some explicit transformation. 
This formula served as motivation for our construction, whereas it is not needed to define the geometric construction. The reader might want to skip the simple, but lengthy calculation.

\subsection{The Whittaker space}

We will denote by $\cC^\infty(\GL_n(\bA))$ the space of functions $f$ on $\GL_n(\bA)$ with values in 
$\Qbar_l$ such that there exists a compact open subgroup $\sfK\subset \GL_n(\bA)$ (depending on $f$) such 
that $f(xk)=f(x)$ for all $x\in \GL(\bA), k\in \sfK$. The same notation will be used for other locally compact groups.

The space of functions 
$$\cC^\infty(\GL_n(\bA))^{\N_n(\bA),\Psi}:=\left\{ f\in \cC^\infty(\GL_n(\bA)) \left| {f(ug)=\Psi(u)f(g) \atop  \forall u\in N_n(\bA), g\in \GL_n(\bA)} \right. \right\}$$
is called {\em Whittaker representation} of $\GL_n(\bA)$. A subrepresentation $\pi$ of this representation of 
$\GL_n(\bA)$ is called a {\em Whittaker model} for the isomorphism class 
of $\pi$. 

Similarly let $\cC^{\infty}_{\text{cusp}}(\GL_n(\bA))^{\P_1(k(C))}$ be the space of functions which are $\P_1(k(C))$-invariant and cuspidal\footnote{A reader unfamiliar with this notion may ignore it for the moment, it will be explained again.} (see \cite{FGKV}). Recall the theorem of Shalika (\cite{Shalika}, 5.9) as stated in loc.cit.:

\begin{Satz}(Shalika) There is an isomorphism of representations of $\GL_n(\bA)$
\begin{eqnarray*}
\Phi: \cC^{\infty}(\GL_n(\bA))^{\N_n(\bA),\Psi} & \to & \cC^{\infty}_{\text{cusp}}(\GL_n(\bA))^{\P_1(k(C))}\\
\text{given by  } f & \mapsto & \Phi(f)(g):= \sum_{y\in \N_{n-1}(k(C))\backslash \GL_{n-1}(k(C))} f(\left(\begin{array}{ll} y & 0 \\ 0 & 1 \end{array}\right) g).
\end{eqnarray*} 
\end{Satz}

Since the character $\Psi$ is a product of characters of the groups $\N_n(K_p)$, we may construct functions in the Whittaker representation as products of functions on $\GL_n(K_p)$ which satisfy the analogous transformation condition for elements of $\N_n(K_p)$. Thus, using Shalika's theorem, the strategy to construct automorphic functions has been to construct functions in the Whittaker model, then to apply $\Phi$ and try to prove that the resulting function is not 
only invariant under the action of $\P_1(k(C))$ but really invariant under the
action of $\GL_n(k(C))$. 

In this chapter we will only be concerned with the local question, i.e.
with representations of $\GL_n(K_p)$ for one fixed prime $p$. The global Whittaker function 
$$W_\sfE (g) := \prod_{p\in C} W_{\sfE,p} (g_p)$$ 
corresponding to our local system will be given as the product of the local functions $W_{\sfE,p}$. These are given by the formula of Shintani and Casselman, Shalika (see \cite{FGKV}) for all $p\in C-S$.  Whereas for $p\in S$ the local factor is the Whittaker function of the Steinberg representation 
(twisted by the eigenvalue $\lambda_p$ of $\Frob_p$ on the one-dimensional stalk $(j_*\sfE)_p$) which is calculated below. 

\subsection{The Steinberg representation}

Fix a point $p\in S\subset C$ and choose a local parameter $\pi$ at $p$.
Let \begin{eqnarray*} 
\delta_\lambda: (K^*_p/\cO_p^*)^n & \to     &  \Qbar_l^*\\
             (\pi^{d_i})     & \mapsto &  \lambda^{\sum_i d_i} \prod_{i<j}  q^{-(d_i-d_j)}. 
\end{eqnarray*}
This  may be viewed as a character of $\B_n(K_p)$, by applying $\delta_\lambda$ to the diagonal entries of an element of $\B_n(K_p)$. In this interpretation $\delta_\lambda$ is the modulus character multiplied by $\lambda^{\text{valuation}(\det)}$.

The (twisted) Steinberg representation $\St_\lambda$ of $\GL_n(K_p)$ is the 
unique irreducible subrepresentation of the induced representation 
$$ \text{Ind}_{\B_n(K_p)}^{\GL_n(K_p)} \delta_\lambda := \left\{ f\in\cC^\infty(GL_n(K_p)) \left| { f(bg)=\delta_\lambda(b) f(g) \atop \forall b\in \B_n(K_p), g\in \GL_n(K_p)}
\right. \right\}.$$
Here again $\cC^\infty(\GL_n(K_p))$ denotes the $\Qbar_l-$valued functions which are invariant under some compact open subgroup. 
For this representation there is a unique (up to scalar) nontrivial $\Iw$-invariant vector, which is an eigenvector of the Iwahori-Hecke algebra
 \cite{Borel}. We denote this vector by $f_\Iw$.  
Furthermore we know that this representation has a Whittaker model, and we denote the $\Iw$-invariant vector in the Whittaker model by $W_\lambda$ and normalize it by the condition that $W_\lambda(1)=1$.                                             

\subsection{The Whittaker function -- statement of the formula}

For any $\un{d}=(d_1,\dots,d_n)\in\bZ^n$ denote by 
$\diag(\un{d}):=${\scriptsize$
\left(\begin{array}{ccc}
\pi^{d_1} &             &          \\ 
          &  \ddots  & \\ 
          &         & \pi^{d_n}
\end{array}\right)$} 
the diagonal matrix and by $\sigma$ the permutation matrix corresponding to the permutation $\sigma(e_i)=e_{\sigma(i)}$. 
 
\begin{satz}\label{formel} The unique $\Iw-$invariant function $W_\lambda$ in the Whittaker model of $\St_\lambda$, normalized by $W_\lambda(1)=1$, is given by: 
$$ W_\lambda(\diag(\un{d})\cdot \sigma)  = 
\left\{ \begin{array}{ll} 
      \frac{\sign(\sigma)\lambda^{\sum d_i}}{q^{\sum_{i<j}d_i-d_j}\vol(\Iw\sigma \Iw)}  &  d_i \geq d_{i+1}-\delta_{\sigma^{-1}(i)>\sigma^{-1}(i+1)} \\ 
     0   & otherwise. 
\end{array}\right.$$
Here $\delta_{\sigma^{-1}(i)>\sigma^{-1}(i+1)}=
\left\{ 
 \begin{array}{ll} 
        1 & \textrm{if } \sigma^{-1}(i)>\sigma^{-1}(i+1) {\text{ i.e. the entry in line i of $\sigma$}\atop \text{ is right of the entry below}} \\
        0 & else.
 \end{array}
\right.,$   

\noindent and the volume is normalized by $\vol(\Iw)=1$.
\end{satz}
\noindent {\bf Remark:} Since 
\begin{eqnarray*} 
\GL_n(K_p) & =  & \B_n(K_p)\GL_n(\cO_p) \\
         & =  &  \cup_{\sigma\in S_n} \B_n(K_p) \Iw \sigma \Iw = \cup_{\sigma\in S_n} \B_n(K_p) \N_n(\cO_p) \sigma \Iw \\
         & =  & \cup_{\sigma \in S_n} \B_n(K_p)\sigma \Iw
\end{eqnarray*}
the proposition is sufficient to calculate $W_\lambda$.    

\noindent {\bf Example:} For $\GL_2$ we have 
$$W_\lambda(\left(\scriptstyle\begin{array}{cc}\pi^{d_1} & \\ & \pi^{d_2}\end{array}\right))=\left\{\begin{array}{ll} q^{d_2-d_1}\lambda^{d_1+d_2} & \text{if } d_1\geq d_2\\ 0 & \text{otherwise,}\end{array}\right.$$ 
$$W_\lambda(\left(\scriptstyle\begin{array}{cc} & \pi^{d_1}  \\ \pi^{d_2} & \end{array}\right))=\left\{\begin{array}{ll} q^{d_2-d_1-1}\lambda^{d_1+d_2} & \text{if } d_1\geq d_2-1\\ 0 & \text{otherwise.}\end{array}\right.$$ 
This is the formula used in Drinfeld's article \cite{Drinfeld_galois_group}.

\subsection{Eigenvalues of some Hecke operators on the Steinberg representation}

To calculate $W_\lambda$, we need first to compute the eigenvalues of the Hecke operators on the $\Iw-$invariant vector in $\St_\lambda$. To this end we use the function $f_\Iw$. 
For an element $g\in \GL_n(K_p)$ we denote by $T_g$ the Hecke operator given by convolution with the 
characteristic function of the double coset $\Iw g \Iw$, i.e.
\begin{eqnarray*}
 T_g : \cC (\GL_n(K_p)/\Iw) & \to & \cC(\GL_n(K_p)/\Iw) \\
 f & \mapsto &  (T_g f)(x):= \sum_{h\in \Iw g\Iw/\Iw} f(x\cdot h).
\end{eqnarray*}
The Hecke operators given by the following particular matrices $t_{\leq i}$ will be very useful\footnote{In case $i=n$ the corresponding operation on parabolic bundles is the upper modification.}:
$$ t_{\leq i} (e_j)  = \left\{
\begin{array}{l}
\pi \cdot e_{i} \textrm{ if } j=1 \\
e_{j-1} \textrm{ if } 1<j\leq i\\
   e_j  \textrm{ if } j>i
\end{array}\right.
\textrm{ i.e. }
t_{\leq i} = 
\left(\begin{array}{cccc|c}
 & 1 & & & \\
 &   & \ddots & &\\
 &   &        & 1 & \\
 \pi &&        &   &  \\
\hline 
  &   &        &   &  \mathbf{1}_{n-i}
      \end{array}\right)
$$

The following --- presumably well known --- lemma gives the eigenvalues of some of the operators 
$T_g$ on $f_\Iw$:
\begin{lemma}\label{eigenwerte}\begin{enumerate}
\item $T_{(\diag(\pi^{d_1},\dots,\pi^{d_n}))}f_\Iw = \lambda^{\sum d_i} f_\Iw$ for all $d_1\geq \dots \geq d_n.$
\item $T_\sigma f_\Iw = \sign(\sigma) f_\Iw$ for all $\sigma \in S_n.$
\item $T_{t_{\leq i}} f_\Iw = (-1)^{i-1} \lambda f_\Iw.$
\end{enumerate}
\end{lemma}
Unfortunately, most of the results on Hecke algebras in the literature are formulated only for 
semi-simple groups. However, the Iwahori-Hecke algebra for $\GL_n$ differs from the one for
$\SL_n$ only by the additional element $T_{t_{\leq n}}$.

\noindent{\bf Proof:} First we note that Borel shows in \cite{Borel} that the eigenvalue of $T_\sigma$ is
$$(T_\sigma)f_\Iw= \sign(\sigma)f_\Iw \textrm{ for }\sigma\in S_n.$$
Further, we may assume $f_\Iw(1)=1$, since $f_\Iw(1)=0$ would imply that $f_\Iw$ is identically $0$ (see the calculations below), thus
$$f_\Iw(\left(\begin{array}{ccc} \pi^{d_1} &  & * \\  & \ddots &  \\  0 &  & \pi^{d_n}\end{array} \right))=  \lambda^{\sum d_i} q^{-\sum_{i<j} d_i-d_j} f_\Iw(1).$$
We apply this in the case $d_1\geq \dots \geq d_n$ to calculate\begin{eqnarray*} 
T_{\diag(\un{d})} f_\Iw(1) & = & \sum_{g\in \Iw \diag(\un{d})  \Iw/\Iw} f_\Iw(g)\\
                                         & = & \sum_{g\in \N_n(\cO) \diag(\un{d}) \Iw/\Iw} f_\Iw(g)\\
                                         & = & \vol \big( \N_n(\cO)\diag(\un{d})\Iw \big) \delta(\diag(\un{d})) f_\Iw(1) \\
                                         & = & q^{(\sum_{i<j} d_i-d_j)} \delta(\diag(\un{d})) f_\Iw(1) \\
                                         & = & \lambda^{(\sum d_i)} f_\Iw(1). 
\end{eqnarray*}
Here we used that an element of $\Iw$ is a product of an element in $\N_n(\cO)$ and a lower diagonal matrix contained in $\Iw$, and that for $d_1\geq \dots\geq d_n$
$$\left(\begin{array}{ccc}
        1 & & \\
        \cp  & \ddots & \\
        \cp & \cp & 1
\end{array}\right) \cdot \diag(\un{d}) = \diag(\un{d}) \cdot \left(\begin{array}{ccc}
        1 & & \\
        \pi^{d_1-d_2}\cp  & \ddots & \\
        \pi^{d_1-d_n}\cp & \pi^{d_2-d_n}\cp & 1
\end{array}\right)\in \diag(\un{d})\cdot \Iw.$$ 
Finally, to compute the eigenvalue of the operator $T_{t_\leq i}$ we first note that by (2):
$$ \sign (\sigma) = T_\sigma f_\Iw(1) = \sum_{g\in \Iw \sigma \Iw / \Iw} f_\Iw(g) = \vol(\Iw\sigma\Iw) f_\Iw(\sigma).$$
 
Further we need a description of the corresponding double coset $\Iw \cdot t_{\leq i}\cdot \Iw/\Iw$. 
Take an element $k\in \Iw$ and look at $k\cdot t_{\leq i}$:
\begin{align*}
\left(\begin{array}{cc} 
\Iw_i & \cO \\ 
\cp & \Iw_{n-i}
\end{array}\right)
\cdot t_{\leq i} 
& =  
\left(\begin{array}{c|c} 
\Iw_i\cdot  t & \cO \\ 
\hline
\underbrace{\pi\cdot\cp}_{1^{st} \textrm{\tiny column}} \,  \cp & \Iw_{n-i}
\end{array}\right) \\
& = 
t_{\leq_i} \cdot 
\left(\begin{array}{c|c} 
\Iw_i &   \begin{array}{l}
                        \pi^{-1} \cO \}{\textrm{\tiny $1^{st}$ line}}\\
                        \cO
          \end{array} \\
\hline
\underbrace{\pi\cdot\cp}_{1^{st} \textrm{\tiny column}} \, \cp & \Iw_{n-i}\end{array}\right).
\end{align*}
Thus, the matrices of the form  
$t_{\leq i}  
\left(\begin{array}{c|c} 
1_i & \begin{array}{c}
                \pi^{-1} v_1 \dots \pi^{-1} v_{n-i} \}{\textrm{\tiny $1^{st}$ line}}\\ 
                0 
      \end{array} \\
\hline
    & 1_{n-i}\end{array}\right) $, $v_i\in \bF_q$ form a set of representatives for $\Iw t_{\leq i} \Iw/\Iw$. 
Set 
$\sigma_{\leq i}^{-1}:= \left( \scriptsize\begin{array}{cccc|c}
0 &   &   & 1 & \\
1 &   &   &   &\\
  & \ddots& & & \\
  &   & 1 & 0 &  \\
\hline 
  &   &   &   &  \mathbf{1}_{n-i}
\end{array} \right),$
then we have 
\begin{eqnarray*}
T_{t_{\leq i}} f_\Iw(\sigma_{\leq i}^{-1})
 & = & 
\sum_{v\in \bF_q^{n-i}} 
f_\Iw(\sigma_{\leq i}^{-1} t_{\leq i}
\left(\scriptsize\begin{array}{ccc|c}
1 &       &     & \pi^{-1} v_1  \dots   \pi^{-1}v_{n-i}\\
 & \ddots &     & \\         
 &        & 1   & \\         
\hline
 &        &     &  \mathbf{1}_{n-i} 
\end{array}\right))\\
& = & \sum_{v\in \bF_q^{n-i}} 
f_\Iw(\left(\scriptsize \begin{array}{cccc|c}
\pi &   &   &  & v_1  \dots   v_{n-i}\\
    & 1 &   &  &                     \\         
    &   & \ddots & &                 \\
    &   &   & 1   &                 \\   
\hline
    &   &   &  &       \mathbf{1}_{n-i}       \\          
\end{array}\right))\\
& = & q^{n-i} \lambda q^{-(n-1)} = q^{-i+1}\lambda. \hspace{13em}\square_{\text{Lemma}}
\end{eqnarray*}

\subsection{The Whittaker function -- proof of the formula}

First we show the vanishing assertion. We know that $W_\lambda(u \cdot g\cdot \gamma)=\psi(u) W\lambda(g)$ for $\gamma\in \Iw$ and $u\in \N_n(K)$. We therefore compute {\scriptsize 
$$
{\underbrace{\left(\begin{array}{ccccc} 1 & u_{1} &  &         \\ 
                             & \ddots & \ddots &        \\
                             &   & 1 & u_{n-1}      \\
                                     &  &   & 1     \\
\end{array}\right)}_{\textstyle =:u}
\diag(\un{d})\sigma}  
 = \diag(\un{d}) \sigma 
\underbrace{\sigma^{-1}
\left(\begin{array}{ccccc} 
1 & \frac{\pi^{d_{2}}}{\pi^{d_1}} u_{1} &  &  &       \\ 
  & 1 &\frac{\pi^{d_{3}}}{\pi^{d_2}} u_{2} &  &      \\
  &   & \ddots & \ddots &     \\
  &   &        & 1 & \frac{\pi^{d_n}}{\pi^{d_{n-1}}}u_{n-1} \\
  &   &        &   & 1     \\
\end{array}\right)  \sigma}_{\textstyle =:u_\sigma \in \Iw \text{?}}.
$$
}
If $u_\sigma\in \Iw$, we must have either $\psi(u)=1$ or $W_\lambda(u\cdot \diag(\un{d})\cdot \sigma )=0$,
i.e. for $\sigma^{-1}(i)>\sigma^{-1}(i+1)$ our function $W_\lambda$ can be non-zero only if $$(\pi^{d_{i+1}-d_i})u_i\in \cp \Rightarrow \res(u_i)=0,$$ that is $d_i\geq d_{i+1}-1$ and if  $\sigma^{-1}(i)<\sigma^{-1}(i+1)$, we need $d_i\geq d_{i+1}$.
This gives the necessary condition for $W_\lambda\neq 0$ claimed in the lemma.

Next, we note that our formula holds for diagonal matrices with $d_1\geq d_2 \geq \dots \geq d_n$, because
\begin{multline*}
\lambda^{\sum d_i}\cdot W(1) = T_{\diag(\un{d})} W(1) =  
\sum_{g\in \Iw\cdot \diag(\un{d}) \cdot \Iw/\Iw} W(g)
              =  \sum_{g\in N(\cO) \diag(\un{d}) \cdot \Iw/\Iw} W(g) \\
                                         =  W(\diag(\un{d})) \cdot \vol(\Iw\cdot \diag(\un{d})\cdot \Iw)
 =   W(\diag(\un{d})) \cdot q^{\sum_{i<j} d_i-d_j}.
\end{multline*}
Now we proceed by descending induction on the number $i$ such that $\sigma(j)=j$ for all $j<i$: Assume that $\sigma(j)=j$ for all $j>i$  and $\sigma^{-1}(i)<i$.

We apply the Hecke operator $T_{t_{\leq i}}$ to express the value of $W(\diag(\un{d})\cdot \sigma)$ for elements $\sigma$ with $\sigma(i)=i-k$ in terms of the value of $W$ at points with $\sigma(i)=i-k+1$, which we know by induction:

Since $W_\lambda$ is an eigenfunction for $T_{t_{\leq i}}$, with eigenvalue $(-1)^{i-1}\lambda$, we get
$$ (-1)^{i-1}\lambda \cdot W(\sigma\cdot\diag(\un{d}))
 =   (T^1_{\leq i} W)(\sigma\cdot\diag(\un{d}))
 =  \sum_{k\in \Iw t_{\leq i} \Iw/\Iw} W(\sigma\cdot \diag(\un{d}) \cdot k).$$
Put $r:=\sigma(i)$ and write 
$\sigma\cdot \diag(\un{d}) =$
{\scriptsize $\left(\begin{array}{cccc|ccc}
 & D_1 &  & & & \\
 &     & & D_r   &&\\         
\vdots &     &       & &&\\         
 & & D_i & &          &\\
\hline
 & &     & & D_{i+1}  & \\
 & &     & &         & \ddots  \\
\end{array}\right)$},
then the above is equal to
{\scriptsize \begin{eqnarray*}
& = &\sum_{v\in \bF_q^{n-i}} W(
\left(\begin{array}{cccc|ccc}
 & & D_1   & & & \\
\pi D_r & & &  & &\\         
 & \vdots    & &  & & \\         
 & & & D_i  &         & \\
\hline
 & &     & & D_{i+1} &  \\
 & &     & &         & \ddots  \\
\end{array}\right)
\left(\begin{array}{c|ccc}
         &  \pi^{-1} v_1 & \dots  & \pi^{-1}v_{n-i}\\
1_{i\times i}    &               &        &\\         
     &               &        &\\         
\hline
           &  1       &        & \\
           &          & \ddots & \\
\end{array}\right)
) \\
& = & \sum W(
\left(\begin{array}{cccc|ccc}
 & & D_1   & & & &\\
\pi D_r & & &  & D_r v_1 & \dots  & D_r v_{n-i}\\         
 & \vdots    & &  & & &\\         
 & & & D_i  &         & &\\
\hline
 & &     & & D_{i+1} & & \\
 & &     & &         & \ddots & \\
 & &     & &         &         & D_n
\end{array}\right))\\
&  = & q^{n-i} W_\lambda(\diag(d_1,\dots,d_r+1,\dots,d_n)\sigma\circ \sigma_{\leq i}^{-1}. 
\end{eqnarray*}}
Here $\sigma_{\leq i}^{-1}$ is the cyclic permutation $(i,i-1,\dots 1)$ as in the proof of Lemma \ref{eigenwerte}.
Note that this gives the sufficient condition for $W_\lambda$ to be non-zero, because we know by induction 
that we must have $d_{\sigma(i)-1} \geq  d_{\sigma(i)} \geq d_{\sigma(i)+1}-1$, or equivalently $d_{\sigma(i)-1}+1 \leq d_{\sigma(i)}+1 \leq d_{\sigma(i)+1}$. To conclude we have to check that we get the right power of $q$ in the induction step:
  \begin{enumerate}
  \item $\vol(\sigma\circ (i,i-1,\dots,1)=:\sigma^\prime)=q^{\#\{k<j|\sigma^\prime(k)>\sigma^\prime(j)\}}$ and we have $${\#\{k<j|\sigma^\prime(k)>\sigma^\prime\}} = \vol(\sigma) - (i-\sigma(i)) + (\sigma(i)-1) = \vol(\sigma) - i + 2 \sigma(i)-1.$$
  \item Write $d^\prime_{\sigma(i)}:=d_{\sigma(i)}+1$ and $d^\prime_j:=d_j$ for $j\neq \sigma(i)$. Then $$q^{\sum_{k<j}d^\prime_k-d^\prime_j} =q^{(\sum_{k<j}d_k-d_j)+(n-\sigma(i)-(\sigma(i)-1)}.$$
  \end{enumerate}
   So these terms differ by a factor $q^{n-i}$, which is what we needed to show.
\hfill$\square$


\section{An analogue of Laumon's construction}

We fix an irreducible local system $\sf E$ of rank $n$ on our curve $C-S$,
ramified at a finite set of points $S\subset C$, 
such that the ramification group at any point $p\in S$ acts unipotently and indecomposably. We will state this condition as ``$\sfE$ has indecomposable unipotent ramification at $S$''.

We want to give a geometric construction for an irreducible perverse sheaf corresponding to the Fourier transform $\Phi(W_\sfE)$ of the Whittaker function $W_\sfE$, computed in the previous section. We will follow Laumon's construction closely, the only new ingredient needed for the construction being the notion of a coherent sheaf with parabolic structure. We will also need to prove generalizations of some results on vector bundles to the case of quasi-parabolic vector bundles.

\subsection{Parabolic vector bundles}

Denote by $\Bun_{n,S}^d$ the moduli space (algebraic stack) of vector bundles of rank $n$ and degree $d$ on $C$ with a full flag at the points of $S$, i.e.:
$$  \Bun_{n,S}^d(T) := \left\langle (\cE,\cE^{(i,p)})_{{i=1,\dots,n-1}\atop{p\in S\hfill }} \left| 
\begin{array}{l} \cE \textrm{ a vector bundle on } C\times T \\ 
\cE \subset \cE^{(1,p)}\subset \dots \subset \cE^{(n-1,p)} \subset \cE(p) \\
\length(\cE^{(i,p)}/\cE)=i \\
\rank(\cE)=n, \deg(\cE)=d\end{array}\right.\right\rangle. $$

Note that the sheaf $\cE^{(i,p)}/\cE$ is automatically flat over $T$ since it is of rank $0$ and its degree is constant.

\noindent{\bf Remark:} Usually one defines a vector bundle with full (quasi-)parabolic structure to be a vector bundle $\cE$
 together with a full flag $V_{1,p}\subsetneqq \dots \subsetneqq V_{n,p} =\cE\tensor k(p)$ of subspaces of the stalk of $\cE$
 at $p$. 
 
This is equivalent to the above definition --- set 
$$\cE^{i,p}:= \bigg(\Ker (\cE \to \cE\tensor k(p)/V_{i,p})\bigg)(p),$$
and conversely $$V_{i,p}:=\Ker \big(\cE\tensor k(p) \to \cE^{(i,p)}\tensor k(p)\big).$$
From this reformulation we get a description of the points of $\Bun_{n,S}^d$: Denote as before $\sfK:=\prod_{p\in(C-S)} \GL_n(\cO_p) \times \prod_{p\in S} \Iw_p$, then\footnote{Recall that given a vector bundle $\cE$ one can choose a trivialisation of $\cE$ at the generic point and at all complete local rings of $C$. The transition functions then give an element of $\GL_n(\bA)$, the double quotient is obtained by forgetting the trivialisations, keeping the flags at $S$.} 
$$\Bun_{n,S}^d(T)(\bF_q) = \GL_n(k(C))\backslash \GL_n(\bA)^{\text{norm}(\det)=d}/\sfK.$$

And the double quotient $\P_1(F)\backslash \GL_n(\bA)/ \sfK$ contains the points of the bundle
 $\Hom^{\text{inj}}(\cO,\cE)\to \Bun_{n,S}^d.$ 

\noindent{\bf Notations:}\begin{enumerate}\item We will write $\cE^{\bullet}:=(\cE,\cE^{(i,p)})_{i=1,\dots,n-1;p\in S}$.

\item Since $\cE \subset \cE^{(i,p)}\subset \cE(p)$ we also get $\cE(p) \subset \cE^{(i,p)}(p)$, thus a parabolic bundle is a chain of vector bundles
$$ \cE^{(i,p)} \subset \cE^{(i+1,p)} \subset \dots \subset \cE^{(n-1,p)} \subset \cE(p)  \subset \cE^{(1,p)}(p) \subset \dots, $$
where the cokernel of every inclusion is of length 1. For any integer $k\in \bZ$ we denote by $\cE^{(k\cdot n +i,p)}:= \cE^{(i,p)}(k\cdot p)$. 

Note furthermore that since the map $\cE\to\cE(p)$ is an isomorphism on $C-\{p\}$, for two distinct points $p,q\in S$ the vector 
bundle $\cE^{(i,p)}+\cE^{(j,q)}\subset \cE(p+q)$ is a vector bundle of degree $d+i+j$. We denote it by $\cE^{(i,p)+(j,q)}$. 
Analogously we define $\cE^{(i,S)}:=\cE^{\sum_{p \in S} (i,p)}$.

Thus we can shift the whole complex to obtain parabolic structures on the vector bundle $\cE^{(i,p)}$ for all $i$. This is called the $i-$th {\em upper modification} of $\cE$. 

\item $\cE(\frac{i}{n}p):=(\cE^{(i,p)},\cE^{(j,q)+(i,p)})_{j=1,\dots ,n-1,q\in S}$.  This notation might be justified, because $\cE^{(i,p)}$ is of degree $d+i=d+n(\frac{i}{n})$ and for $i=n$ we get the canonical parabolic structure on the vector bundle $\cE(p)$.
\end{enumerate}
We now want to mimic Laumon's construction of automorphic sheaves for unramified local systems. Consider for example the case of bundles of rank $2$.
We will view $\Phi(W_\sfE)$ as a function on vector bundles together with a meromorphic section of $\Omega$. At a point $\Omega\hookrightarrow \cE$ such that $\Omega \to \cE$ and $\Omega\to\cE^{(1,S)}$ are both maximal embeddings
$\Phi(W_\sfE)$ is defined as the sum over all sections of $\cE^0/\Omega$ with at most simple poles at $S$. But the line bundle $(\cE/\Omega)(S)\cong \cE^{(1,S)}/\Omega$, thus we might equivalently sum over all holomorphic sections of $\cE^{(1,S)}/\Omega$. 

To apply a similar consideration to bundles of larger rank, our calculation of $W_\sfE$ suggests that we need to consider quotients of $\cE^\bullet$ by subsheaves which are not maximal. We therefore look for a notion of coherent sheaves with parabolic structure\footnote{While I was thinking about this, Norbert Hoffmann explained to me that one can formally adjoin quotients of vector bundles with
 parabolic structure to the category  of such bundles to obtain an abelian
 category. The definition below may be viewed as a geometric interpretation of 
  these quotients. I would like to thank him for the helpful discussion.} which allows the operation $\cF^\bullet \mapsto \cF^\bullet (\frac{i}{n}S)$. This is easy with the above definition of parabolic structure:

\subsection{Parabolic coherent sheaves}
\begin{definition}
A {\em coherent sheaf on C with n-step parabolic structure at S} -- also called {\em parabolic sheaf} for short -- is a collection of coherent sheaves $\cF^{\bullet}:=(\cF=\cF^{(0,p)},\cF^{(i,p)})_{i=1,\dots,n-1;p\in S}$ together with morphisms $\phi^{(i,p)}:\cF^{(i,p)}\to \cF^{(i+1,p)}$ for $i=1,\dots,n$ and $p\in S$ (where $\cF^{(n,p)}:=\cF(p)$)
such that in the resulting sequence
$$\dots \map{\phi^{(n,p)}(-p)} \cF \map{\phi^{(1,p)}} \cF^{(1,p)} \map{\phi^{(2,p)}} \dots \map{\phi^{(n-1,p)}} \cF^{(n-1,p)} \map{\phi_n} \cF(p) \map{\phi^{(1,p)}(p)} \cF^{(1,p)}(p) \dots $$
the composition of $n$ maps $\xymatrix@C=14ex{{\cF}^{(i,p)}\ar[r]^{\phi^{(i-1,p)}(p)\circ \dots \circ \phi^{(i,p)}}& {\cF}^{(i,p)}(p)}$ is the natural morphism.
\end{definition}

\noindent{\bf Note:}\begin{enumerate}\item If the sheaf $\cF^{(i,p)}$ is not torsion free at $p$ for some $i$,
then  the natural map $\cF^{(i,p)} \to \cF^{(i,p)}(p)$ is not injective, so at least one of the 
$\phi^\bullet$'s is not injective (see the examples below).
\item The {\em degree} of $\cF^\bullet$ is defined as the collection $\deg(\cF^\bullet):=(\deg(\cF^{(i,p)}))_{0\leq i<n\atop p\in S}$.
\item Denote the algebraic stack of coherent sheaves of rank $r$ on $C$ with $n$-step parabolic structure at $S$ and
(multi-)degree $\underline{d}=(d^{(i,p)})_{0\leq i<n,p\in S}$ by $\Coh_{r,C,S}^{\underline{d}}$. 
Since we usually fix the curve $C$, we will omit it and write $\Coh_{r,S}^{\underline{d}}$ to shorten this lengthy notation.
\item As in the case of vector bundles we define $\cF^{(i,p)+(j,q)}:=(\cF^{(i,p)}\oplus\cF^{(j,q)})/\cF$ (for the diagonal embedding of $\cF$). Note that this quotient is the sheaf isomorphic to $\cF^{(i,p)}$ on $C-q$ and isomorphic to $\cF^{(j,q)}$ on $C-p$. These sheaves glue, since both are canonically isomorphic to $\cF$ on $C-\{p,q\}$. 
Analogously we define $\cF^{(i,S)}$. 
\item Again we define upper modifications as $\cF^\bullet(\frac{i}{n}p):=(\cF^{(i,p)},\cF^{(j,q)+(i,p)})_{0\leq j<n \atop q\in S}$.
\end{enumerate}

\noindent{\bf Example:} In our case, given a morphism $\Omega^{\tensor(n-1)} \to \cE$, we get an induced parabolic structure on the quotient $\cE/\Omega^{\tensor (n-1)}$. We only use that $\cE(p)/\Omega^{\tensor (n-1)}(p)=(\cE/\Omega^{\tensor (n-1)})(p)$ to get
{\small\begin{equation}\label{beispiel}
\vcenter{\xymatrix@C=1.5em{
\Omega^{\tensor (n-1)}\ar[d]\ar[r]^{Id} & \Omega^{\tensor (n-1)}\ar[d]\ar[r]^{Id} &  \dots\ar[r] & \Omega^{\tensor (n-1)}\ar[d]\ar[r] & \Omega^{\tensor (n-1)}(p)\ar[d] \\
{\cE} \ar[r]^{\phi^{(1,p)}}\ar[d]          & {\cE}^{(1,p)} \ar[r]^{\phi^{(2,p)}} \ar[d]  &  \dots\ar[r]&{\cE}^{(n-1,p)} \ar[r]^{\phi^{(n,p)}}\ar[d] & {\cE}(p) \ar[d]   \\
{\cE}/\Omega^{\tensor (n-1)} \ar[r]       & {\cE}^{(1,p)}/\Omega^{\tensor (n-1)} \ar[r] &  \dots\ar[r] &{\cE}^{(n-1,p)}/\Omega^{\tensor (n-1)} \ar[r] &(\cE/\Omega^{\tensor (n-1)})(p). \\
}}
\end{equation}}
Note that we can view $\Omega^{\tensor (n-1)}$ (or any coherent sheaf) as parabolic sheaf by defining $\Omega^{\tensor (n-1) (i,p)}:=\Omega^{\tensor (n-1)}$ for $i=0,\dots,n-1$. With this definition the above diagram is an extension of parabolic sheaves.

From this example we see that:
\begin{lemdef} 
The category of (quasi--)coherent sheaves with n-step parabolic structure is abelian.

We denote homomorphisms of parabolic sheaves by $\Hom_{\text{para}}(\underline{\quad},\underline{\quad})$, and the same for $\Ext^1_{\text{para}}$, etc.

The category of quasi--coherent sheaves has enough injectives.
\end{lemdef}
\noindent{\bf Proof:} The kernel and cokernel of a morphism can be defined componentwise. All compatibilities thus follow from the corresponding ones for coherent sheaves. 

Furthermore the above example shows that:
\begin{remark}\label{einbettung} The stack $\Coh_{k,C}^d$ classifying coherent sheaves of rank $k$
and degree $d$ can be embedded into the stack of parabolic sheaves:
\begin{eqnarray*}
j: \Coh_{k,C}^{d} &\hookrightarrow & \Coh_{k,S}^{(d,\dots,d)}\\
       \cF           & \mapsto        & \cF^\bullet:=(\cF,\cF^{(i,p)}:=\cF)
\end{eqnarray*}
For a coherent sheaf $\cF$ on $C$ we will write $(\cF)^\bullet$ for its image $j(\cF)$. The functors $(\underline{\quad})^\bullet$ and $(\underline{\quad})^{(0,S)}$ are adjoint functors:
$$ \Hom_{\text{para}}((\cF)^\bullet,\cG^\bullet) = \Hom_{\cO_C}(\cF,\cG^{(0,S)})$$
and
$$ \Hom_{\text{para}}(\cG^\bullet,(\cF)^\bullet)= \Hom_{\cO_C}(\cG^{(n-1,S)},\cF).$$
\end{remark}
For an injective sheaf $\cI$ the adjunction property yields 
$$\Hom_{\text{para}}(\cF^\bullet,(\cI)^\bullet)=\Hom_{\cO_C}(\cF^{(n-1,S)},\cI).$$ Since the functor $(\un{\quad})^{(n-1,S)}$ is exact we conclude that $\Hom_{\text{para}}(\un{\quad},(\cI)^\bullet)$ is exact. Thus choosing embeddings $\cG^{i,p}\hookrightarrow \cI_{i,p}$ into injective sheaves $\cI_{i,p}$ we get an embedding $\cG^\bullet \hookrightarrow \oplus (\cI_{i,p})^\bullet(\frac{n-1}{n}S+\frac{n-i-1}{n}p)$ of $\cG$ into an injective parabolic sheaf.
\hfill $\square_{\textrm{\tiny Lemma}}$ 

By the above we also have:
\begin{lemma}\label{extpara}
 The extensions of a parabolic sheaf $\cF^\bullet$ by a line bundle $\cL$  are classified by $\Ext^1_{\cO_C}(\cF^{(n-1,S)},\cL)$, i.e. 
$$\Ext^1_{\text{para}}(\cF^\bullet,(\cL)^\bullet)=\Ext_{\cO_C}^1(\cF^{(n-1,S)},\cL).$$
\end{lemma}
\noindent{\bf Proof:} By the above remark any injective resolution of $\cL$ defines an injective resolution of $(\cL)^\bullet$, and to such a resolution we may apply the adjunction formula. 

\hfill$\square_{\textrm{\tiny Lemma}}$

Note that we could give another proof of this lemma, calculating the Yoneda-Ext groups directly via the diagram (\ref{beispiel}). The only thing one has to check is that in this diagram we have $\cE^{(i,p)} \cong \cE^{(n-1,p)}\times_{(\cE^{(n-1,p)}/\Omega^{\tensor n-1})}(\cE^{(i,p)}/\Omega^{\tensor n-1})$.
 
\begin{kor}\label{serre-parabolisch} Let $\cF^\bullet$ be a parabolic sheaf and let $\cL$ be a line bundle on $C$. Then we
have by Serre duality:
$$ \Ext^1_{\text{para}}(\cF^\bullet, (\cL)^\bullet) = (\Hom_{\text{para}}((\cL \tensor \Omega^{-1})^{\bullet}(-\frac{n-1}{n}S), \cF^\bullet))^\vee.$$
\end{kor}
\noindent{\bf Proof:} This is just an application of the adjunction isomorphism to
$$ \Ext^1_{\cO_C}(\cF^{(n-1,S)},\cL) = \Hom_{\cO_C}(\cL\tensor\Omega^{-1},\cF^{(n-1,S)}).$$ 

\hfill$\square$

The above version of Serre duality (Corollary \ref{serre-parabolisch}) suggests to denote $\cO^\bullet:=(\cO)^\bullet$, 
$\Omega^\bullet:=(\Omega)^\bullet(\frac{n-1}{n}S)$ and 
analogously $\Omega^{\bullet,k}:= (\Omega^{\tensor k})^\bullet(k\frac{n-1}{n}S)$.
Then we can put $\cL:=\Omega^{\tensor k}$ to deduce from the corollary that:
$$ \Ext_{\text{para}}^1(\cF^\bullet,\Omega^{\bullet,k}) \cong \Hom_{\text{para}}(\Omega^{\bullet,k-1},\cF^\bullet)^\vee.$$

\subsection{The fundamental diagram}

Reformulating the preceding calculations for families of parabolic sheaves allows us to construct a variant of Laumon's ``fundamental diagram'' as follows: Denote by $\cE_{\univ}^\bullet$ (resp. $\cF^\bullet_{\univ}$) the universal parabolic sheaf on $\Bun_{n,S}^d\times C$ (resp. on $\Coh^{\un{d}}_{n,S}\times C$) and let $p_i$ be the projection to the $i-$th factor.

We can view the sheaf $p_{1,*}(\cH om(p_2^*\Omega^{\bullet,n-1},\cE^\bullet_{\univ}))$ as the classifying stack for parabolic vector bundles $\cE^\bullet$ together with a homomorphism $\Omega^{\bullet,n-1}\to \cE^\bullet$. Denote this stack by:
$$ \Hom_n(T):=\langle (\cF^\bullet,\Omega^{\bullet,n-1}\map{\phi}\cF^\bullet) \; | \; \cF^\bullet\in \Coh^{\un{d}}_{n,S} \rangle. $$
Write $\Hom_n^{\text{inj}}$ for the open substack of $\Hom_n$ where $\phi$ is injective.

Similarly write $\Ext^1_n$ for the stack classifying extensions of parabolic sheaves by $\Omega^{\bullet,n}$:
$$ \Ext_n^1(T):=\langle 0\to \Omega^{\bullet,n} \to \cF^\bullet_{n+1} \to \cF^{\bullet} \to 0 \; | \; \cF^\bullet\in\Coh^{\un{d}}_{n,S}\rangle. $$

Note that there are open substacks $\Bun_{n,S}^{d,\text{good}}\subset \Bun_{n,S}^{d}$ and 
$\Coh_{n,S}^{\underline{d},\text{good}}\subset \Coh_{n,S}^{\underline{d}}$ such that the restrictions of $\Hom_n$ and $\Ext^1_n$ to these substacks are vector bundles. 
More precisely we will call a coherent parabolic sheaf $\cF^\bullet$ {\em good} if 
$$ \Hom_{\text{para}}(\cF^\bullet, \Omega_{C/k}^{\bullet,n-i+1})=0 \textrm{ for all } 1\leq i \leq n-1.$$
By Serre duality this condition guarantees that 
$$\Ext^1_{\text{para}}(\Omega^{\bullet,n-i},\cF^\bullet)=\Ext^1_{\cO_C}(\Omega^{n-i},\cF^{(-(n-i)(n-1),S)})=0,$$
and moreover the same will be true for any quotient of $\cF^\bullet$. 

We define $\Coh_{n,S}^{\un{d},\text{good}}$ to be the stack of good parabolic bundles. 
Over these stacks the semi-continuity theorem tells us that the sheaves $p_* \cH om_{\text{para}}(p_2^*\Omega^{\bullet,n-i},\cF^\bullet_{\univ})$ and $\bR^1 p_{1,*}(\cH om_{\text{para}}(p_2^*\Omega^{\bullet,n},\cF^\bullet_{\univ}))$ are indeed vector bundles over $\Coh_{n,S}^{\un{d},\text{good}}$. Furthermore we have:
$$\Ext_n^1 = \bR^1 p_{1,*} (\cH om_{\text{para}}(p_2^*\Omega^{\bullet,n},\cF^\bullet_{\univ}))$$
since the corresponding $\bR^0 p_{1,*}$ vanishes for good sheaves. We will always consider the stacks $\Ext_n^1$ and $\Hom_n$ as stacks over $\Coh_{n,S}^{\text{good}}$.

As in Laumon's construction we have:
\begin{enumerate}
\item To give a short exact sequence $0 \to \Omega^{\bullet,n-1} \to \cF_n^\bullet \to \cF^\bullet \to 0$ it is sufficient to define the datum $0 \to \Omega^{\bullet, n-1} \to \cF^{\bullet}_{n}$. To restrict this remark to good parabolic sheaves we denote by $\Ext_n^{1,\text{good}}\subset \Ext_n^1$ the substack consisting of extensions in which the middle term is a good parabolic sheaf (therefore the right term is good as well). We then have an isomorphism
$$ I_n: \Hom_n \map{\cong} \Ext_n^{1,\text{good}}.$$
\item Over $\Coh^{\un{d}}_{n,S}$ the bundles $\Hom_n$ and $\Ext_n^1$ are dual vector bundles.
\end{enumerate}

Therefore we can define a fundamental diagram (which we split into several diagrams):

{\small
\begin{eqnarray}
&\vcenter{\xymatrix{
\Hom_n\ar[d] & \ar@{_(->}[l]_-{j_{\Hom}} \Hom^{\text{inj}}_n \ar[r]^{I_n}_{\cong}& \Ext^{1,\text{good}}_{n-1} \ar@{^(->}[r]^{j_{\Ext}} & \Ext^{1}_{n-1}\ar[dr]\\ 
\Coh_{n,S}^{\un{d},\text{good}} & & & & {\Coh^{\underline{d}_{n-1,S},\textrm{\tiny \text{good}}}_{n-1,S}}
}}& \dots \nonumber\\
&\vcenter{\xymatrix{
\Ext_{n-1}^1 \ar[dr]\ar@/^1.5pc/@{<..>}[rr]^{\textrm{\tiny dual bundles}} & & \Hom_{n-1}\ar[dl]\\
& {\Coh^{\underline{d}_{n-1},\text{good}}_{n-1,S}} & 
}}& \dots \label{Teil1}\\
\dots & \vcenter{\xymatrix{
& \Hom_{n-1}\ar[dl] & \ar@{_(->}[l]_-{j_{\Hom}} \Hom^{\text{inj}}_{n-1} \ar[r]^{I_{n-1}}_{\cong}& \Ext^{1,\text{good}}_{n-2} \ar@{^(->}[r]^{j_{\Ext}} & \Ext^{1}_{n-2}\ar[d] \\ 
\Coh_{n-1,S}^{\un{d},\text{good}} & & & & {\Coh^{\underline{d}_{n-2},\text{good}}_{n-2,S}}
}}& \dots \nonumber
\end{eqnarray}
}

Here the last line is the same as the first one with $n$ replaced by $n-1$. Thus we can continue this to end up with $\Coh_{0,S}^{\un{d}_0}$ (we drop the superscript ``good'', since all torsion sheaves are good) . We have to keep track of the degrees of the parabolic sheaves:
\begin{eqnarray*}
\underline{d}_{n-i}^{(*,p)} =(d_{n-i}-(n-i-1),d_{n-i}-(n-i-2),\dots,\underbrace{d_{n-i},\dots,d_{n-i}}_{i+1 - times})\\
\textrm{ with } d_{n-i}:=d-\sum_{j=1}^i ((n-j)(2g-2)+(n-j+1)).
\end{eqnarray*}
In particular, continuing the above diagram to the right, the last term will be $\Coh_{0,S}^{(d_0,\dots,d_0)}$. Laumon's construction started with a sheaf on $\Coh_{0}^{d_0}$ which corresponds to the Whittaker function for unramified local systems. This sheaf is pulled back to $\Ext_0^{1,\text{good}}$, then one applies the Fourier transform for the bundles in (\ref{Teil1}) and then continues with pull backs and intermediate extensions for the maps $j_{\Hom}$ and $j_{\Ext}$ in the upper line of the diagram  until one ends up with a sheaf on $\Hom_n$.

Thus, in our situation we need to find a sheaf on $\Coh_{0,S}^{\underline{d}_0}$ that corresponds to the Whittaker function as calculated in Section 1.

\subsection{The Whittaker sheaf $\cL_\sfE^d$}
As noted in Section \ref{einbettung}, there is an open embedding of torsion sheaves of degree $d_0$ on $C-S$ to parabolic torsion sheaves: 
\begin{eqnarray*}
j: \Coh^{d_0}_{0,C-S} &\hookrightarrow & \Coh_{0,S}^{\underline{d}_0}\\
       \cT           & \mapsto        & \cT^\bullet=(\cT,\cT^{(i,p)}:=\cT).
\end{eqnarray*}
  The map $j$ is open, since the condition that $\supp(\cT^{(0,S)})\subset C-S$ is open.
 
Furthermore we have Laumon's Whittaker sheaf $\cL_{\sfE|_{C-S}}^{d_0}$ on $\Coh_{0,C-S}^{d_0}$.  Recall the the definition of  $\cL_{\sfE|_{C-S}}^{d_0}$:  Let $\sfE|_{C-S}^{(d_0)}$ be the symmetric product of $\sfE$ restricted to the symmetric product $(C-S)^{(d_0)}$ of the curve $C-S$. Denote $j_{C-S}: (C-S)^{(d_0)} \map{D\mapsto \cO/\cO(-D)} \Coh_{0,C-S}^{d_0}$, which is almost an embedding (see \cite{Laumon_correspondance}). Then $\cL_{\sfE|_{C-S}}^{d_0}:= j_{C-S,!*}\sfE|_{C-S}^{(d_0)}.$

\begin{definition}\label{Whittakergarbe}
We define the {\em Whittaker sheaf} corresponding to $\sfE$ to be
$$\cL_{\sfE}^d := j_{!*} \cL_{\sfE|_{C-S}}^d = (j\circ j_{C-S})_{!*} (\sfE|_{C-S})^{(d)}.$$
\end{definition}

We will prove some properties of the Whittaker sheaf justifying its name in Section \ref{section_whittaker_sheaf}.

\subsection{Putting everything together: The Fourier transform of $\cL_\sfE^d$}
Now let $\quot:\Ext^1(\cF^\bullet_0,\cO^\bullet) \to  \Coh_{0,S}^{\underline{d}_0}$, ${(\cO^\bullet\hookrightarrow \cF_1^\bullet\tto\cF_0^\bullet)\mapsto \cF_0^\bullet}$ be the quotient map  and denote the Fourier transform by $\Four : D^b(\Hom_k)\to D^b(\Ext^1_k)$. Following the fundamental diagram (\ref{Teil1}) from right to left we define:

\begin{definition}\label{definition_FE} We inductively define the sheaves $\sfF_\sfE^k$ and $\sfF_{\sfE!}^k$ on $\Hom_k^{\text{inj}}$ as
\begin{eqnarray*}
\sfF_{\sfE}^1 & := & I_1^* j_{\Ext}^* \quot^* \cL_{\sfE}^{d_0}[d_0] =: \sfF_{\sfE,!}^1,\\
\sfF_{\sfE}^{k+1} & := & I_{k+1}^* j_{\Ext}^* \Four(j_{\Hom,!*}\sfF_{\sfE}^{k}),\\
\sfF_{\sfE,!}^{k+1} & := & I_{k+1}^* j_{\Ext}^* \Four(j_{\Hom,!}\sfF_{\sfE,!}^{k}).
\end{eqnarray*}
\end{definition}

The restriction of the sheaf $\sfF_\sfE^n$ to the stack of vector bundles with a section of $\Omega^{\bullet,n-1}$ will be our candidate to descend to an automorphic sheaf on $\Bun_{n,S}^{\underline{d}}$. 
By construction this is an irreducible perverse sheaf (because this property is preserved by $\Four$, $j_{\Hom,!*}$ and $j_{\Ext}^*$).

As in \cite{Laumon_premiere_construction} we also define the sheaves $\sfF_{\sfE,!}^{k}$, because it will be easy to prove that these have a Hecke eigensheaf property, and finally (in Section \ref{section_sauberkeit}) 
we will show that they are isomorphic to $\sfF_{\sfE}^k$ for $k\leq n\leq 3$.

To end this section we want to state our main theorem. To do this we need to define geometric Hecke operators for parabolic sheaves. We first give an example
indicating the relation between parabolic torsion sheaves and the Iwahori-Hecke algebra:

\subsection{Parabolic torsion sheaves and Hecke operators}
Assume for the moment that $n=2$, $S=\{p\}$, and consider the stack $\Coh_{0,p}^{1,1}$. Take any  $\cT^\bullet\in \Coh_{0,p}^{1,1}$. If $\supp(\cT)=q\neq p$, then $\cT^\bullet\cong(k(q))$, but if $\supp(\cT)=p$, then $\cT^\bullet$ is isomorphic to exactly one of the following sheaves:
\begin{enumerate}
\item $ \cT_0=k_p \map{id} \cT_1=k_p \map{0} \cT_0(p)=k_p \map{id}\dots $
\item $ \cT_0=k_p \map{0} \cT_1=k_p \map{id} \cT_0(p)=k_p \map{0}\dots $
\item $ \cT_0=k_p \map{0} \cT_1=k_p \map{0} \cT_0(p)=k_p \map{0}\dots $
\end{enumerate}
We want to relate these sheaves to some Hecke operators of the Iwahori-Hecke algebra at $p$, acting on parabolic vector bundles of rank $2$.
To do this, we consider torsion free extensions of a vector bundle $\cE^{\prime\bullet}$ by the first complex:
$$\xymatrix@C=2em@R=2ex{
\ar[r] & {\cE}^{\prime (0,p)} \ar[r]\ar[d] & {\cE}^{\prime (1,p)}\ar[r]\ar[d] & {\cE}^{\prime (0,p)}(p)\ar[r]\ar[d] & {\cE}^{\prime (1,p)}(p)\ar[d]\ar[r] & \dots\\
\ar[r] & {\cE}^{ (0,p)} \ar[r]^{\phi^{ 1}}\ar[d] & {\cE}^{(1,p)}\ar[r]^{\phi^{ 2}}\ar[d]\ar@{..>}[ur] & {\cE}^{(0,p)}(p)\ar[r]^{\phi^{ 1}(p)}\ar[d] & {\cE}^{(1,p)}(p)\ar[d] \ar[r]& \dots \\
\ar[r] & k_p \ar[r]^{id} &  k_p\ar[r]^{0} & k_p \ar[r]^{id} & k_p \ar[r] & \dots\\
}$$

The middle map in the lower sequence is $0$, therefore $\phi^{ 2}$ factors through $\cE^{(1,p)}\to\cE^{\prime (0,p)}(p)$. Since all the bundles $\cE^{(i,p)}$ are locally free this map is injective, and since the two bundles have the same degree it is an isomorphism.

The same argument shows that $\phi^{ 1}$ does not factor through $\cE^{ (1,p)}$, so the upper line is given by a parabolic structure on the vector bundle $\cE^{ (1,p)}$, different from the canonical structure $(\cE^{ (1,p)})^\bullet$.
Thus extensions of this type are the set indexing the summation of the Hecke-Operator $T_{\left({0\atop 1}{1 \atop 0}\right)}\circ T_{\left({0 \atop \pi}{1\atop 0}\right)}$.
According to Lemma \ref{eigenwerte} this operator should act with eigenvalue $\Trace(\Frob_p,{\sf E}_p)$ on the Whittaker sheaf.
Analogously we find that summing over extensions of parabolic bundles by the second torsion sheaf calculates $T_{\left({0 \atop \pi}{1 \atop 0}\right)}\circ T_{\left({0\atop 1}{1 \atop 0}\right)}.$ 
Finally the third torsion sheaf gives the Hecke operator $T_{\left({ 0\atop \pi }{ 1 \atop 0 }\right)}$ which acts with eigenvalue $-\Trace (\Frob_p,{\sf E}_p)$ on the Whittaker function. Note that this torsion sheaf is a 
point of codimension 2 in $\Coh_{0,p}^{1,1}$, and thus the perverse sheaf $\cL_\sfE$ will have some $H^1$ at this point. The minus sign of the eigenvalue will come from taking the trace of $\Frob$ on this cohomology group of odd degree (see Section 4).
Therefore we define {\em generalized Hecke operators} as follows:

Fix non negative degrees $\un{d}=\un{d_1}+\un{d_2}$, and let $\Hecke_n^{\underline{d_1},\underline{d_2}}$ be the stack classifying extensions of parabolic sheaves of degree $\un{d_2}$ by torsion sheaves of degree $\un{d_1}$, i.e.
$$ \Hecke_n^{\underline{d_1},\underline{d_2}} := \langle (0 \to \cF^{\prime\bullet} \to \cF^\bullet \to \cT^{\bullet} \to 0) | \cF^{\prime\bullet} \in \Coh_{n,S}^{\un{d_2}}, \cT^{\prime\bullet} \in \Coh_{0,S}^{\un{d_1}} \rangle.$$
The forgetful maps give rise to a correspondence
$$\xymatrix{
 & \Hecke_n^{\underline{d_1},\underline{d_2}}\ar[dl]_-{pr_{\text{big}}}\ar[dr]^-{pr_{\text{small}}\times quot} &\\
\Coh_{n,S}^{\underline{d_1}+\underline{d_2}} & & \Coh_{n,S}^{\underline{d_2}}\times \Coh_{0,S}^{\underline{d_1}}.
}$$
\begin{definition}\label{Heckeoperatoren} The {\em Hecke operator} $H^{\un{d_1},\un{d_2}}_n$ is defined as
\begin{eqnarray*}
H^{\underline{d_1},\underline{d_2}}_n: D^b(\Coh_{n,S}^{\un{d}}) & \to & D^b(\Coh_{n,S}^{\un{d_2}}\times \Coh_{0,S}^{\un{d_1}}) \\
 \sfK & \mapsto & H^{\underline{d_1},\underline{d_2}}_n\sfK :=\bR (pr_{\text{small}}\times quot)_! \circ pr_{\text{big}}^* \sfK.
\end{eqnarray*}
\end{definition}
To define operators on parabolic vector bundles which correspond to the action of the Iwahori-Hecke algebra we have to introduce for every 
$(0,\dots,0)\leq \un{\epsilon} \leq (1,\dots,1)$ the stack 
$$\overline{\Coh}_{0,S}^{\un{\epsilon}}:=\Coh_{0,S}/\textrm{central $\bG_m$-automorphisms}.$$ 
And we define more Hecke correspondences:
$$\xymatrix@C=2ex{ 
& {\left\langle {\cE}^{\prime\bullet} \subset {\cE}^\bullet \left| \begin{array}{l} {\cE}^\bullet \in \Bun_{n,S}^{\un{d}} \\ {\cE}^{\prime\bullet}\in \Bun_{n,S}^{\un{d}-\un{\epsilon}} \end{array}\right.\right\rangle} \ar[dl]_-{pr_{\text{big}}}\ar[dr]_-{pr_{\text{small}}\times quot}\ar[drr]^-{pr_{\text{small}}\times pr_C} & & \\
\Bun_{n,S}^{d} & & \Bun_{n,S}^{d-\un{\epsilon}} \times \overline{\Coh}_{0,S}^{\un{\epsilon}}\ar[r] & \Bun_{n,S}^{\un{d}-\un{\epsilon}} \times C.
}$$
Here $\Bun_{n,S}^{d-\un{\epsilon}}$ denotes the moduli space of torsion free parabolic sheaves of degree $d^{i,p}=d+i-\epsilon^{i,p}$, i.e. the moduli space of vector bundles with partial parabolic structure at $S$.

\begin{definition}
    The {\em Hecke operator} $H^{\un{\epsilon}}$ is defined by:
\begin{eqnarray*}
H^{\un{\epsilon}}: D^b(\Bun_{n,S}^{d}) & \to & D^b(\Bun_{n,S}^{d-\un{\epsilon}}\times \overline{\Coh}_{0,S}^{\un{\epsilon}}) \\
 \sfK & \mapsto & H^{\underline{\epsilon}}\sfK :=\bR (pr_{\text{small}}\times quot)_! \circ pr_{\text{big}}^* \sfK.
\end{eqnarray*} 
For $\un{\epsilon}=(1,\dots,1)$ we set:
\begin{eqnarray*}
H^1_C : D^b(\Bun_{n,S}^{d}) & \to & D^b(\Bun_{n,S}^{d-1} \times C) \\
\sfK & \mapsto & H^{1}_C\sfK :=\bR (pr_{\text{small}}\times pr_C)_! \circ pr_{\text{big}}^* \sfK.
\end{eqnarray*}
\end{definition}

Finally note that the sheaf $\cL_\sfE^1$ descends to a sheaf $\overline{\cL}_\sfE^1$ on $\overline{\Coh}_{0,S}^{\un{1}}$. Denote by $j_C:C-S \hookrightarrow C$ the inclusion. We will prove the following:

\begin{Satz}\label{ergebnis} Let $\sfE$ be a local system on the curve $C-S$ with indecomposable unipotent ramification at $S$ and assume $n=rank(\sfE)\leq3$. 
Then 
\begin{enumerate}
\item $\sfF_\sfE^n\cong \sfF_{\sfE,!}^n.$ 
\item $\sfF_\sfE^n$ descends to
a nonzero perverse sheaf $\sfA^{good}_\sfE$ on $\Bun_{n,S}^{\un{d},good}$. 
\item $\sfA^{good}_\sfE$ extends to a {\em Hecke eigensheaf} $\sfA_\sfE$ on $\Bun_{n,S}$, i.e. there is a unique extension $\sfA_\sfE$ of $\sfA_\sfE^{good}$ to $\Bun_{n,S}$ such that: 

\begin{eqnarray*}
H^{\un{1}} \sfA_\sfE & \cong & \sfA_\sfE \boxtimes \overline{\cL}_\sfE^1\\
H^{\un{\epsilon}} \sfA_\sfE  & = & 0 \textrm{\hspace{3em} for } \un{0} < \un{\epsilon} < \un{1}\\
H^{1}_C \sfA_\sfE & \cong & \sfA_\sfE \boxtimes j_{C,*}\sfE\\
\textrm{and the isomorphism}\\
H^1_C \circ H^1_C \sfA_\sfE &\cong&  \sfA_\sfE \boxtimes j_{C,*}\sfE \boxtimes j_{C,*}\sfE \textrm{ is } S_2\textrm{--equivariant.} 
\end{eqnarray*}
\end{enumerate}
\end{Satz}

Furthermore, we will show that this implies that the function $\tr_{\sfA_\sfE}$ is an eigenfunction for 
the action of the Iwahori-Hecke algebra. Indeed, by the example given above we have already seen that the points of $\overline{\Coh}_{0,S}^{\un{1}}$ give a set of generators for the Iwahori-Hecke algebra
(the invertible element corresponding to $\tensor \cO(\frac{1}{n}p)$ and the operators corresponding
to the transpositions in $S_n$ generate the algebra).


\section{Some properties of parabolic sheaves}

This section is an attempt to clarify the notion of parabolic sheaves. First we give a description of the isomorphism classes
of parabolic torsion sheaves, then we prove some lemmata concerning homological algebra of parabolic sheaves. 
At the end of this section we give an explicit description of the moduli space of torsion sheaves on $\bA^1$ with parabolic structure at $0$.
All these results are simple, but for completeness they are collected in this paragraph.

\subsection{The structure of parabolic torsion sheaves}

The structure theorem for modules over principal ideal domains shows that any torsion sheaves on a curve $C/k$ is a direct sum of sheaves of the form $\cO/(\cp^d)=:\cO_{dp}$ for some prime ideals $\cp$. We will prove a similar result for parabolic torsion sheaves. 
The constituents of a sheaf $\cT^\bullet$ supported in $p\in S$ will be of the form (we only give the sheaves in degree $(*,p)$)
$$ \dots \to \cO/ \cp^d \to \dots \to \cO/ \cp^d \tto \cO/\cp^{(d-1)} \to \dots \to \cO/\cp^{(d-1)} \hookrightarrow \cO/\cp^d \to \dots; $$
more precisely these are isomorphic to $\cO^\bullet_{\frac{k}{n}p}(\frac{i}{n}p):=\cO^\bullet(\frac{i}{n}p)/\cO^\bullet(\frac{i-k}{n}p)$ for some $0\leq k < i \in \bN$, and in the sequence above $d=\lceil\frac{k}{n}\rceil$ is the smallest integer bigger than $\frac{k}{n}$. The key step is to prove:
\begin{lemma}\label{maximale_summanden}
Let $\cT^\bullet$ be a parabolic torsion sheaf supported in $p\in S$, and let further $\cO^\bullet_{\frac{k}{n}p}(\frac{i}{n}p) \incl{\psi} \cT^\bullet$ be an inclusion such that
the sum of the degrees of the torsion sheaves occurring in $\cO^\bullet_{\frac{k}{n}p}(\frac{i}{n}p)$ is
maximal. Then there is a splitting of $\psi$.
\end{lemma}
\noindent{\bf Proof:} Since $\cT^\bullet$ is supported at $p$, we may assume that $S=\{p\}$. We choose a local parameter $\pi$ at $p$. Shifting our sequences we may further assume that $i=0$. 

Note that any inclusion $\cO_{d^\prime p}\hookrightarrow \cT^{(m,p)}$ gives rise to an inclusion of some $\cO^\bullet_{\frac{k^\prime}{n}}(\frac{m}{n})\hookrightarrow \cT^\bullet$ with $d^\prime=\lceil\frac{k^\prime}{n}\rceil$. Thus for a maximal embedding $\psi$ we know that $d=\lceil\frac{k}{n}\rceil$ is the maximal length of the torsion sheaves occurring in $\cT^\bullet$. In particular, any inclusion $\cO_{dp}\hookrightarrow \cT^{(m,p)}$ splits. Thus we have
$$\xymatrix{
 {\cO_{(d-1)p}} \ar[r]^-{\cong}\ar@{^(->}[d]^{\psi_0}  &  \dots \ar[r]^-\cong & {\cO_{(d-1)p}} \ar@{^(->}[r]\ar@{^(->}[d] & {\cO_{dp}} \ar[r]^-\cong\ar[d]^{\psi_r} & \ar@{^(->}[d]^{\psi_{n-1}}\ar@{->>}[r]{\cO_{dp}} & \ar@{^(->}[d]{\cO_{(d-1)p}}\\
 {\cT}^0 \ar[r] & \dots \ar[r] & {\cT}^{r-1}\ar[r] & {\cT}^r \ar[r]^{\phi_{n-1}}& {\cT}^{n-1}(p) \ar[r] & {\cT}^0(p) \\
}$$
(where the number of submodules of the form $\cO_{(d-1)p}$ might be zero). We know that there is a splitting for $\psi_{n-1}$, and this induces compatible splittings for $\psi_l$ for $r\leq l < n-1$. In particular, if $\frac{k}{n}\in \bN$ (i.e. in the upper line of the diagram all terms are of the form $\cO_{dp})$, then any splitting of $\psi_{n-1}$ induces a compatible splitting for $\psi$. We may therefore assume that $\psi_0:\cO_{(d-1)p} \hookrightarrow \cT^0$ and that $d-1\neq 0$.

{\em Claim:} There is a splitting of $\psi_{r-1}$.

Otherwise $\psi_{r-1}(1)=\pi \cdot e_{r-1}$ for some $e_{r-1}\in \cT^{r-1}$ with $\pi^{d-1}e_{r-1}\neq 0$. Then $\pi^{d-1}\phi_{n-1}(e_{r-1})=\pi^{d-2}(\phi_{n-1}(\pi e_{r-1})= \pi^{n-1}\psi_{n-1}(1)\neq 0$ contradicting the maximality of $\frac{k}{n}$.
 
Thus we only need to find a compatible splitting of $\psi_{r-1}$ and $\psi_{n-1}$.  To do this, we may replace $\cT^{r-1}$ by its image in $\cT^{n-1}$, since the above argument still works for this replacement. We then have $\cT^{r-1} \hookrightarrow \cT^{n-1}$, and the cokernel of $\phi_n$ is of length $1$. In this case any splitting of $\psi_{r-1}$ can be extended to one of $\psi_{n-1}$. (Choose a decomposition $\cT^{r-1}\cong \cO_{(d-1)p}\psi_{r-1}(1)\oplus \bigoplus_i \cO_{d_ip}e_i$. Then either $\phi_{n-1}(e_i)$ generates a direct summand of $\cT^{n-1}$, or $\phi_{n-1}(e_i)=\pi e_i^\prime$, and $e_i^\prime$ generates a direct summand. Completing this to a decomposition of $\cT^{n-1}$ into indecomposable sheaves we can define an extension of $\psi_{r-1}$.) 
\hfill $\square$

From this lemma we get:
\begin{satz}{(Structure of parabolic torsion sheaves)}\label{struktur}\begin{enumerate}
\item Any parabolic torsion sheaf is a direct sum of sheaves of the form
$$ \cO_{\frac{j}{n}p}({\textstyle \frac{i}{n}}p)^\bullet = \cO^\bullet({\textstyle \frac{i}{n}}p)/\cO^\bullet({\textstyle \frac{i-j}{n}}p), \; i,j\in \bN, \; p\in S$$
and sheaves supported outside $S$.
\item Any parabolic torsion sheaf $\cT^\bullet$  has a filtration $\cT^\bullet_j\subset \cT^\bullet_{j+1}\subset \dots \subset \cT^\bullet$ such that the filtration quotients $\cT_{j+1}^\bullet/\cT_j^\bullet$ are isomorphic to one of the following:
  \begin{enumerate}
   \item $ \cT_{j+1}^\bullet/\cT_j^\bullet \cong (k(q))^\bullet$ and  $q\not\in S$
   \item There is a $p_0\in S$ and $0\leq i_0<n$ such that
$$\cT_{j+1}^{(i,p)}/\cT_j^{(i,p)} = \left\{\begin{array}{ll} k(p_0) & i=i_0,p=p_0\in S \\ 0 & \textrm{else.}\end{array}\right.$$
  \end{enumerate}
\item Any parabolic torsion sheaf $\cT^\bullet$ of constant degree $deg{\cT^\bullet}=(d,\dots,d)$
has a filtration $\cT_1^\bullet\subset\cdots\subset\cT_i^\bullet \subset\cdots\subset \cT^\bullet$ such that $deg(\cT_i^\bullet)=(i,\dots,i)$.
\end{enumerate}
\end{satz} 
\noindent{\bf Proof:} Since for any torsion sheaf $\cT$ we have a canonical decomposition $\cT = \oplus_{q\in supp(\cT)} \cT\tensor\cO_{C,p}$, we may assume that $\cT^\bullet$ is a parabolic torsion sheaf concentrated in a single point $q$, i.e. $\supp(\cT^{(i,p)})=q$ for any $(i,p)$. 

If $q\not\in S$, we know that all the $\cT^{(i,p)}$ are isomorphic because the functor $\tensor \cO_C(S)$ is the identity functor on sheaves supported in $C-S$. Hence $\cT^\bullet= (\cT^{(0,p)})^\bullet$ and for torsion sheaves without extra structure the lemma holds.

For torsion sheaves supported in $S$ the previous lemma implies (1) and the sheaves $\cO_{\frac{k}{n}p}(\frac{j}{n}p)$ have a filtration satisfying (2) Counting degrees we also get (3). \hfill $\square$

Finally note that for an arbitrary parabolic sheaf the torsion subsheaves are always a direct summand:
\begin{remark}
Let $\cF^\bullet$ be a parabolic sheaf on $C/k$. Then $\cF^\bullet=\cE^\bullet\oplus \cT^\bullet$, where $\cT^\bullet$ is
a parabolic torsion sheaf and all $\cE^{i,p}$ are torsion free.
\end{remark}
\noindent{\bf Proof:} We know that $\cT^\bullet:= \text{torsion}(\cF^\bullet)\subset \cF^\bullet$ is a parabolic 
torsion sheaf and $\cF^{(0,S)}\cong \cT^{(0,S)} \oplus \cE^{(0,S)}$.
And since the $\phi^{(i,p)}$ are isomorphisms over the generic fibre of $C$ the images 
$\phi_{i,p}(\cE^{(0,p)})$ can be used to define maximal torsion free subsheaves of $\cF^{(i,p)}$,
 these define the desired decomposition. \hfill $\square$

\subsection{Homological algebra of parabolic sheaves}

\begin{lemma}
For coherent parabolic sheaves on $C_{/k}$ the functors $\Ext^i_{\text{para}}$ vanish for $i>1$. 
\end{lemma}
\noindent{\bf Proof:} Let $\cF^\bullet$ be a parabolic sheaf. We prove that $\Ext_{\text{para}}^i(\;\;,\cF^\bullet)=0$ 
for $i>1$ by descending induction on the rank and degree of $\cF^\bullet$.

For a line bundle $\cL$ on $C$ the functor $\Hom_{\text{para}}(\;\; , (\cL)^{\bullet}(\frac{i}{n}S))$ 
coincides with a $\Hom-$functor on coherent sheaves, and for $\Ext_{\cO_C}^i$ the lemma holds. By induction, we may therefore assume that
$\cF^\bullet$ is a parabolic torsion sheaf. 
By Lemma \ref{struktur} giving the structure of parabolic torsion sheaves, we may further 
assume that $\cF^\bullet$ is a quotient of two line bundles of arbitrarily high degree, 
which establishes the claim. \hfill $\square_{\textrm{\tiny Lemma}}$

\begin{lemma}\label{dim_Hom_und_Ext}
Let $\cT^\bullet$ be a parabolic torsion sheaf and $\cE^\bullet$ a parabolic vector bundle. Then:\begin{enumerate}
\item $\dim(\Ext^1_{\text{para}}(\cT^\bullet,\cE^\bullet))$ and $\dim(\Hom_{\text{para}}(\cE^\bullet,\cT^\bullet))$ only depend on $\rank(\cE^\bullet)$, $\deg(\cE^\bullet)$ and $\deg(\cT^\bullet)$. 

\item More precisely, for $\cT^{(i,p)}=\left\{\begin{array}{ll} k(p_0) & i=i_0,p=p_0 \\ 0 & \textrm{else.}\end{array}\right.$ we have $$\dim(\Ext^1_{\text{para}}(\cT^\bullet,\cE^\bullet))= \deg(\cE^{(i+1,p)})-\deg(\cE^{(i,p)}).$$
$$\dim(\Hom_{\text{para}}(\cE^\bullet,\cT^\bullet))= \deg(\cE^{(i,p)})-\deg(\cE^{(i-1,p)}).$$

\item If $\deg(\cT^\bullet)=(d)$ is constant, we get 
$$\dim(\Ext^1_{\text{para}}(\cT^\bullet,\cE^\bullet))= d \cdot \rank(\cE) = \dim(\Hom_{\text{para}}(\cE^\bullet,\cT^\bullet)).$$
\end{enumerate}
\end{lemma}
\noindent{\bf Proof:} We give a proof of the statements on $\Ext^1_{\text{para}}$, the case of homomorphisms is even simpler. Since $\cE^\bullet$ is torsion free, $\Hom_{\text{para}}(\cT^\bullet,\cE^\bullet)=0$. Thus for any exact sequence $0\to \cT^{\prime \bullet} \to \cT^\bullet \to \cT^{\pprime\bullet}$ the sequence
$$ 0 \to Ext^1_{\text{para}}(\cT^{\prime\bullet},\cE^\bullet) \to Ext^1_{\text{para}}(\cT^{\bullet},\cE^\bullet) \to Ext^1_{\text{para}}(\cT^{\pprime\bullet},\cE^\bullet) \to 0$$
is exact as well.

To prove the lemma, apply this remark to the filtration $\cT_i^\bullet \subset \cT^\bullet$ constructed in Lemma \ref{struktur} (1) and reduce to the case $\cT^{(i,p)}=\left\{\begin{array}{ll} k(p_0) & i=i_0,p=p_0 \\ 0 & \textrm{otherwise}\end{array}\right.$. 
We may shift $\cE^\bullet,\cT^\bullet$ and assume that $i_0=0$. Write $\cT= (\cL)^\bullet / (\cL)^\bullet(-\frac{1}{n}p_0)$ for some line bundle $\cL$ and for simplicity choose $\deg(\cL)<<0$  such that $\Ext^1_{\cO_C}(\cL,\cE^{i,p})=0$ for all $p\in S$,$-1\leq i\leq n$. Then
\begin{eqnarray*}
\dim(\Ext^1_{\text{para}}(\cT^\bullet,\cE^\bullet)) & = & \chi(\bR \Hom_{\text{para}}(\cL^\bullet(-\frac{1}{n}p_0),\cE^\bullet)) - \chi(\bR \Hom_{\text{para}}(\cL^\bullet,\cE^\bullet))\\
& = & \chi(\bR \Hom_{\cO_C}(\cL,\cE^{(1,p_0)})) - \chi(\bR \Hom_{\cO_C}(\cL,\cE^{(0,p_0)}))\\
& = & \dim ( \Hom_{\cO} (\cL,\cE^{(1,p_0)})) - \dim(\Hom_{\cO} (\cL,\cE^{(0,p_0)}))\\
& = & \deg(\cE^{(1,p_0)}) - \deg(\cE^{(0,p_0)})   
\end{eqnarray*}\hfill $\square_{\textrm{\tiny Lemma}}$

\subsection{The moduli stack of parabolic torsion sheaves}

First let us consider the moduli stack of torsion sheaves on $\bA^1$ with parabolic structure at $p=0$ as an example:

This stack classifies sequences of torsion sheaves\footnote{I suppress the upper index $p$ since we have assumed that $S=\{p\}=\{0\}$}:
$$ \cdots \to \cT^0 \map{\phi_1} \cT^1 \map{\phi_2} \cT^2 \map{\phi_3} \dots \map{\phi_{n-1}} \cT^{n-1} \map{\phi_n} \cT^0(p) \map{\phi_1(p)} \cdots $$
with the property that the induced maps $\cT_i \to \cT_i(p)$ are the natural ones.

Recall that a single torsion sheaf $\cT$ on $\bA^1$ can be described by giving its vector space of
global sections $H^0(\bA^1,\cT)$ together with the endomorphism given by multiplication by the 
coordinate $t$ of $\bA^1=\Spec(k[t])$. Hence we get a presentation of the moduli space of 
torsion sheaves of degree $d$ on $\bA^1$:
$$ \Coh_{0,\bA^1}^{d} \cong [\Mat_{d,d}/\GL_d], $$
where $\GL_d$ acts on $\Mat_{d,d}$ by conjugation. (Under this identification the support of a sheaf is given by the eigenvalues of the corresponding matrix, and the length of the indecomposable summands is given by the Jordan decomposition.)

For torsion sheaves with parabolic structure we can define a similar presentation as follows:

Note that the natural map $\cT^i\to\cT^i(p)$ is given by the multiplication by the coordinate $t$.
Thus for any collection $(\phi_i:k^{\oplus d_{i-1}}\to k^{\oplus d_i})_{i=1}^n$ we may define $\cT_j$ by
 $(k^{\oplus d_j}, \phi_j\circ \phi_{j-1} \cdots \phi_1 \circ \phi_n \circ \cdots \circ \phi_{j+1})$ and with this definition the 
 $\phi_i$ automatically define homomorphisms $\cT_{i-1}\to \cT_i$ of $\cO_{\bA^1}$-modules.  
This proves:
\begin{lemma}
$$\Coh_{0,\{p\}}^{d_0,\dots,d_{n-1}}\cong [\Mat_{d_1,d_0}\times \Mat_{d_2,d_1} \times \dots \times \Mat_{d_{0},d_{n-1}}/(\GL_{d_0}\times\dots\times \GL_{d_{n-1}})],$$
where an element $(g_0,\dots,g_{n-1})\in \GL_{d_0}\times\dots\times \GL_{d_{n-1}}$ operates on 
$(\phi_1,\dots,\phi_n)\in \Mat_{d_1,d_0}\times \Mat_{d_2,d_1} \times \dots \times \Mat_{d_{0},d_{n-1}}$ as
$$ (g_0,\dots,g_{n-1})\cdot (\phi_1,\dots,\phi_n) := (g_1\phi_1 g_0^{-1},g_2\phi_2 g_1^{-1},\dots,g_{0} \phi_n g_{n-1}^{-1}).$$ \hfill$\square_{\textrm{\em\tiny Lemma}}$
\end{lemma}

\begin{kor}\label{coh-glatt}
For any smooth curve $C$ and any finite set $S\subset C(k)$ the stack $\Coh_{0,S}^{\underline{d}}$ is a smooth stack of relative dimension $0$. 
\end{kor}
\noindent{\bf Proof:} To show the lifting property for smoothness at a point $\cT^\bullet\in \Coh_{0,S}^{\un{d}}$, we only need to consider sheaves on $\Spec(\Pi_{q\in supp(\cT)} \widehat{\cO}_{C,q})$. But for a smooth curve we know that $\widehat{\cO}_{C,q}\cong k[[t]] \cong \widehat{\cO}_{\bA^1,0}$, and therefore it is sufficient to prove the corollary in case $C=\bA^1$ and $S=\{0\}$, which is proven in the previous lemma. \hfill $\square_{\textrm{\tiny Corollary}}$

In case one does not want to consider deformations of parabolic sheaves one could use the above lemma
and the fundamental diagram to get smooth presentations of the stacks $\Coh_{n,S}^{\underline{d}}$:

\begin{kor} 
For any smooth curve $C$ and any finite set $S\subset C(k)$ the stacks $\Coh_{n,S}^{\underline{d}}$ 
are smooth algebraic stacks.
\end{kor} 


\section{Properties of the Whittaker sheaf $\cL_{\sfE}^d$}\label{section_whittaker_sheaf}

Our main goal in this section is to prove a Hecke property of the sheaf $\cL_\sfE^d$ defined in \ref{Whittakergarbe} (Proposition \ref{hecketorsion}). In the case of unramified local systems Laumon proved this in two steps: First he introduced a small resolution of the stack of torsion sheaves, defined as the stack classifying torsion sheaves, together with a full flag of subsheaves. Thereby he obtained  a geometric description of the Whittaker sheaf, which he then used to prove the Hecke property. 

Translating this into our situation we encounter two problems. The first one is that $\cL_\sfE^1$ is already a complex of sheaves. The second problem is that the analogue of Laumon's resolution is not small in the case of parabolic torsion sheaves. 

Since $\cL_\sfE^d$ is a perverse sheaf on the moduli stack of parabolic torsion sheaves and most of the questions are local in the \'etale topology we will often be able to reduce to the case that our curve is $\bA^1$ and our local system is ramified only at the point $0$. Our first aim is therefore to calculate $\cL_\sfE^1$ in this case. After translating these results into the general situation we then proceed with Laumon's strategy as described above. Here we simultaneously prove that the Hecke property of $\cL_\sfE^d$ holds and that we can give a geometric description of $\cL_\sfE^d$.

\subsection{Calculation of the sheaf $j_{!*}\sfE$ on $\Coh^{\un{1}}_{0,\bA^1,0}$}\label{A1}

Consider the case $C=\bA^1$ and $S=\{0\}$. Let $\sfE_n$ be the $n-$dimensional local system
 on $\bG_m$, ramified at $0$, such that the ramification group acts unipotently and indecomposably
 --- i.e. the invariants under the ramification group are $1$-dimensional --- constructed as follows:
We have $\Ext_{\bG_m}^1(\bQ_\ell(-1),\bQ_\ell)=H^1(\bG_m,\bQ_\ell(1))= H^1(\bG_m,\bQ_\ell)(1)= \bQ_\ell$ and 
therefore there is a unique nontrivial extension $\sfE_2$ of the sheaf $\bQ_\ell(-1)$ by 
the constant sheaf $\bQ_\ell$. The long exact cohomology sequence corresponding to this extension gives 
$H^1(\bG_m,\sfE_2)=\bQ_\ell(-2)$, thus we can repeat this argument to define $\sfE_n$, filtered by 
$\bQ_\ell=\sfE_1\subset\sfE_2\subset\dots\subset\sfE_{n-1}\subset\sfE_n$ with subquotients 
$\sfE_i/\sfE_{i-1}\cong \bQ_\ell(-i+1)$. 
Alternatively we could describe this as $Sym^{n-1} (\sfE_2)$. 

Since $\Coh_{0,\bG_m}^1$ -- the stack of torsion sheaves of length $1$ on $\bG_m$ -- is 
isomorphic to $[\bG_m/\bG_m]$ for the trivial action of $\bG_m$ on $\bG_m$, the sheaf 
$\sfE_n$ descends to a sheaf on $\Coh_{0,\bG_m}^1$ which we denote again by $\sfE_n$.

We want to calculate the middle extension $j_{!*}\sfE_n$ with respect to the inclusion 
$ j: \Coh_{0,\bG_m}^1 \to \Coh_{0,\bA^1, \{0\}}^{1,\dots,1} $ 
(where $\Coh_{0,\bA^1, \{0\}}^{1,\dots,1}$ is the stack of torsion sheaves with $k-$step parabolic 
structure -- in this section we allow $n\neq k$).
Because of the theorem on smooth base change it is sufficient to do this on a smooth 
representation of these stacks, for example:

$$\xymatrix{ 
& {\bG_m}^k \ar[r]^{j} \ar@{->>}[dd]\ar[dl]_{m}\ar@{}[ddr]|\square & \bA^k \ar@{->>}[dd] \\
\bG_m\ar[dr]_{/ \bG_m-\text{trivial}}  & & \\
&  \Coh_{0,\bG_m}^{1} \ar[r]^{j} & \Coh_{0,\bA^1,0}^{1,\dots,1}.
}$$

So we are left calculating $j_{!*} m^*\sfE_n$ on $\bA^k$.

Denote $D_i:=\{x_i=0\}\subset \bA^k$ and for
a subset $I\subset\{1,\dots,k\}$ define $D_I:=\cap_{i\in I} D_i$. This stratification of the complement of $\bG_m^k$ gives rise to open immersions $j_i$: 
$${\bG_m}^{k}=\bA^k-\cup D_i \incl{j_1} \bA^k-\cup D_{ij} \incl{j_2} \dots \incl{j_{k-1}} \bA^k-(0,\dots,0) \incl{j_k} \bA^k $$ 

And $j_{!*}m^*\sfE_n = \tau^{<k}\bR j_{k,*} \tau^{<k-1} \bR j_{k-1} \cdots \tau^{<1} \bR j_{1,*}m^*\sfE_n$ 
(this is a definition in Intersection Homology II \cite{GoreskyMacPherson} and a proposition 
(2.1.11) in Faisceaux Pervers \cite{faisceaux_pervers}). 

To simplify notation, let $U_i:=\bA^k-\cup_{\# I=i} D_I$ and denote $\sfE_n:=m^*\sfE_n$.

For $k=1$ we have $R^p j_* \sfE_n|_0 = \left\{ \begin{array}{ll} \bQ_\ell & p=0 \\ \bQ_\ell(-n) & p=1 
\end{array}\right.$ on $\bA^1$. 

Therefore on $\bA^k$ we get that the stalk at $0$ is: 
$R^pj_* \sfE_n|_0 \stackrel{(*)}{\cong} H^p(\bA^{k}_{\overline{\bF}_q},\bR j_*\sfE_n)=H^p(\bG^k_{m,\Fbar_q},\sfE_n )$ 
and this isomorphism is compatible with the action of the Galois group. 
The equality $(*)$ holds, because $\sfE_n$ is an extension of constant sheaves, 
for which the two cohomology groups are canonically isomorphic. 

To calculate the latter cohomology group, we can factor $m$ into an isomorphism 
$\bG_m^k \map{ (a_i)\mapsto (\Pi a_i, a_2,\dots,a_n) } \bG_m^k $ 
followed by the projection onto the first factor, to obtain:
\begin{eqnarray*}
H^*(\bG_m^k,m^* \sfE_n ) & \cong & H^*(\bG_m^k,pr_1^* \sfE_n ) \cong H^*(\bG_m,\sfE_n) \tensor H^*(\bG_m^{k-1},\bQ_\ell)\\
                         & = &\left\{\begin{array}{ll} \bQ_\ell & *=0 \\ \bQ_\ell(-1)^{\oplus k-1} \oplus \bQ_\ell(-n) & * =1 
                         \\ \textrm{etc.}  & \\\end{array}\right.
\end{eqnarray*}
Analogously we get a formula for the stalk of $R^pj_*\sfE_n$  at a point lying on exactly $r$ of 
the divisors:  
$$R^*j_* \sfE_n|_{D_{(i_1,\dots,i_r)}-\cup_{j\not\in\underline{i}}D_{(i_1,\dots,i_r,j)}}\cong H^*(\bG_m,\sfE_n)\tensor H^*(\bG_m^{r-1}\times \bA^{k-r},\bQ_\ell).$$
If the terms of weight $\geq 2n$ did not appear, then the truncation functors $\tau^{<i}$ used 
in the definition of $j_{!*}\sfE_n$ would be trivial and $\bR j_* \sfE_n$ would ``be an irreducible perverse sheaf''. 
But these terms do disappear if we pass to the inductive limit of all the 
$\sfE_n\hookrightarrow \sfE_{n+1} \hookrightarrow \cdots$. Therefore define $\sfE_\infty := \varinjlim \sfE_n$. We have the following proposition:

\begin{satz}\label{e_unendlich}  For $n\geq k$ there is an exact triangle of complexes on $\bA^k$:
$$ \to j_{!*} \sfE_n \to  \bR j_* \sfE_\infty \to j_! \sfE_\infty (-n) \map{[1]} $$
\end{satz}

\noindent{\bf Proof:} (inductively calculating $\tau^{<i}\bR j_{i,*}$)  Use the shorthand $j_{i\dots 1}:=j_i\circ \cdots \circ j_1$.

We start with the exact sequence of sheaves on $\bG_m^k$:
$$0 \to \sfE_n \to \sfE_\infty \to \sfE_\infty (-n) \to 0.$$

Applying $\tau^{<1}\bR j_{1,*}=j_{1,*}$ we get
$$ 0 \to j_{1,*} \sfE_n \to j_{1,*} \sfE_\infty \to j_{1,!}\sfE_\infty(-n) \to 0.$$

Using the previous calculation, that $j_{1,*}\sfE_\infty=\bR j_{1,*}\sfE_\infty$  
we get an exact triangle of complexes on $U_2$
$$ \to j_{1,!*} \sfE_n \to \bR j_{1,*} \sfE_\infty  \to j_{1,!} \sfE_\infty (-n) \map{[1]}.$$
Because $j_{i\dots 1,!} \sfE_\infty(-n)$ is an extension of $j_{i\dots 1,!}\bQ_\ell(-n-r)$ with 
$r\geq 0$, we will need to calculate $\bR j_{i+1,*} j_{i\dots 1,!} \bQ_\ell$ via the sequence
$$ 0 \to j_{i\dots 1,!} \bQ_\ell \to \bQ_\ell \to \bQ_\ell|_{\cup_{j=1}^k D_j} \to 0,$$
Where $\bQ_\ell|_{\cup_{j=1}^k D_j}$ is the constant sheaf on the union of the divisors. (We will
often use this shorthand: For a closed subscheme $Z\incl{i}X$ and a sheaf $\sfK$ on $X$ we write
$\sfK|_Z := i_*i^* \sfK$.)

For the last term --- viewed as a sheaf on the whole of $\bA^k$ --- we have a resolution
$$ 0 \to \bQ_\ell|_{\cup_{j=1}^k D_j} \to \oplus_{j=1}^k \bQ_\ell|_{D_j} \to \oplus_{1\leq i_1<i_2\leq k} \bQ_\ell|_{D_{i_1,i_2}} \to \dots \to \bQ_\ell|_{D_{1,\dots,k}} \to 0.$$
Restricting this resolution to $U_i=\bA^k-\cup_{1\leq k_1< \dots < k_i\leq k} D_{k_1\dots k_i}$ 
all terms $\bQ_\ell|_{D_I}$ with $|I|\geq i$ disappear and thus on $U_{i+1}$ we have a resolution
$$  0 \to j_{i\dots 1,!} \bQ_\ell \to \bQ_\ell \to \oplus_{j=1}^k \bQ_\ell|_{D_j} \to \dots \to \oplus_{I\subset\{1,\dots,k\} \atop |I|=i} \bQ_\ell|_{D_I} \to 0.$$

\begin{lemma}\label{eindimlemma}
For any $m>i\geq 0$ the complex $j_{m\dots i+1,!*}j_{i\dots 1,!} \bQ_\ell$ is quasi-isomorphic to
$$ j_{m\dots i+1,*}(0\to \bQ_\ell \to \oplus_{j=1}^k \bQ_\ell|_{D_j} \to \dots \to \oplus_{I\subset\{1,\dots,k\} \atop |I|=i} \bQ_\ell|_{D_I} \to 0).$$ 
\end{lemma}
\noindent{\bf Proof:} For $i=0$ there is nothing to prove, so we may assume $i>0$. 
Note that for $|I|=c$ we have 
$$\bR^p j_{i+1,*} (\bQ_\ell|_{D_I}) = \left\{ \begin{array}{ll} \bQ_\ell|_{D_I} & p=0 \\ \oplus_{I^\prime \supset I; |I^\prime|=k+1} \bQ_\ell(k-i-1)|_{D_{I^\prime}} & p = 2(i+1-c)-1\\ 0 & \textrm{otherwise.}\end{array}\right.$$
because $j_{i+1}$ adds a smooth boundary of codimension $i+1-c$ to $D_I$.
Therefore looking at the spectral sequence calculating $\bR j_{i+1,*}j_{i\dots 1,!}\bQ_\ell$ via
our resolution of $j_{i\dots 1,!}\bQ$ we see that 
the only terms appearing in cohomological dimension $<i+1$ are as claimed. 
This proves the lemma for $m=i+1$. Inductively we may apply the same argument for $m$ to see that
in our spectral sequence the cohomology in degrees 
$p$ with $i+1\leq p < 2(m-i)-1+i=2m-(i+1)$ vanishes.
\hfill $\square_{\textrm{\tiny Lemma}}$  

We will use the last statement of the above proof again:
\begin{lemma}\label{eindimplus} 
For $m>i+1$ we we have
$$ j_{m,!*} (j_{m-1\dots i+1,!*}j_{i\dots 1,!} \bQ_\ell) \cong \tau^{<2m-(i+1)} \bR j_{m,*} (j_{m-1\dots i+1,!*}j_{i\dots 1,!} \bQ_\ell)$$
\end{lemma}\hfill$\square$

To finish the proof of Proposition \ref{e_unendlich}, we still have to calculate 
$$\tau^{<i+1}\bR j_{i+1,*} \Big(  \to j_{i\dots 1,!*} \sfE_n \to \bR j_{i\dots 1,*} \sfE_\infty \to j_{i\dots 1,!} \sfE_\infty(-n)\map{[1]} \Big).$$ 
By our calculation (Lemma \ref{eindimlemma}) of $j_{i+1,!*} j_{i\dots 1,!}\bQ_\ell=\tau^{< i+1} \bR j_{i+1,*} j_{i\dots 1,!} \bQ_\ell$ we know that
$$ \bR^p j_{i+1,*} j_{i\dots 1,!} \sfE_\infty (-n) = \left\{ \begin{array}{ll}
 j_{i+1\dots 1,!} \sfE_\infty (-n) & p = 0 \\
 0                                 & 0<p<i \\
 \textrm{of weight }\geq 2n        & p= i \end{array}\right. .$$
Considering the long exact cohomology sequence for $$\bR j_{i+1,*} j_{i\dots 1,!*} \sfE_n \to \bR j_{i+1\dots 1,*} \sfE_\infty \to \bR j_{i+1,*} j_{i\dots 1,!}\sfE_\infty (-n) \to$$ this calculation implies that the map 
$$R^i j_{i+1\dots 1,*} \sfE_\infty \to R^i j_{i+1,*} j_{i\dots 1,!} \sfE_\infty(-n)$$
must be zero because the weights of the two sheaves are distinct. Thus we have proven the proposition.
  \hfill $\square_{\text{Proposition \ref{e_unendlich}}}$


Later we will need the following description of $(j_{!*}{\sfE}_\infty)|_{D_I}$, which is 
implicit in the above:

\begin{lemma}\label{restriction_to_D}
The morphism $j_{1,*}{\sfE}_\infty \to i_{D_1,*}i_{D_1}^* j_{1,*}{\sfE}_\infty = i_{D_1,*}(\bQ_\ell|_{D_1\cap U_2})$
on $U_2=\bA^k-(\cup_{i\neq j} D_{ij})$ induces an isomorphism 
\begin{align*}
(\bR j_* {\sfE}_\infty)|_{D_1}&  \map{\cong} \bR j_{k\dots 2,*} \bQ_\ell|_{D_1\cap U_2},\\
\intertext{therefore}
(j_{!*} {\sfE}_n)|_{D_1}& \map{\cong} \bR j_{k\dots 2,*} \bQ_\ell|_{D_1\cap U_2},\\
\intertext{and more generally}
 (j_{!*} {\sfE}_n)|_{D_{1\dots l}}&  \map{\cong} \bR j_{k\dots l,*} \big((\bR j_{l-1\dots 2,*} \bQ_\ell|_{D_1\cap U_2})|_{D_{1\dots l-1}}\big).
\end{align*}
\end{lemma}

\noindent{\bf Proof:} For the first statement consider $j_{D_1}: U_1=\bA^k-\cup_i D_i \hookrightarrow \bA^k-\cup_{i>1} D_i$ and $j_1^\prime:  \bA^k-\cup_{i>1} D_i \hookrightarrow U_2=\bA^k-(\cup_{i\neq j} D_{ij})$. 
This induces an exact sequence
$$ 0 \to j_{D_1,!} {\sfE}_\infty \to j_{D_1,*} {\sfE}_\infty \to  i_{D_1,*}\bQ_\ell|_{D_1} \to 0.$$
We therefore have to show that $\big(\bR (j_{k\dots 2}\circ j^\prime_1)_* j_{D_1,!}\sfE_\infty\big)|_{D_1}=0$. Again we first show that the stalk at $0$ vanishes. We know that $R^p(j_{k\dots 2}\circ j^\prime_{1})_* j_{D_1,!}\sfE_\infty|_0 = H^p(\bA^k-\cup_{i>1} D_i,j_{D_1,!}{\sfE}_\infty)$, because this is true for the other two sheaves in the sequence above 
(for the middle term we proved this to calculate $\bR j_* \sfE_\infty$).

The cartesian diagram
$$\xymatrix{ 
{\bG}_m^k \ar@{^(->}[r]^-{j_{D_1}}\ar[d]_{(a_i)\mapsto(\Pi a_i,a_2,\dots a_k)}^{\cong} & {\bA}^1\times \bG_m^{k-1} \ar[d]^{(a_i)\mapsto(\Pi a_i,a_2,\dots a_k)} _{\cong} \\
{\bG}_m^k \ar@{^(->}[r]^-{j_{D_1}}\ar[d]^{pr_1} &  {\bA}^1\times \bG_m^{k-1} \ar[d]^{pr_1} \\
{\bG_m} \ar@{^(->}[r]^{j_{\bG_m}} & {\bA}^1.
}$$
shows that 
\begin{multline*}
H^*(\bA^1\times \bG_m^{k-1},j_{D_1,!} {\sfE}_\infty) = H^*(\bA^1 \times \bG_m^{k-1}, pr_1^* j_{\bG_m,!} {\sfE}_\infty) \\  = H^*(\bA^1,j_{\bG_m,!} {\sfE}_\infty)\tensor H^*(\bG_m^{k-1},\Qbar_l) = 0,
\end{multline*}
 because  $H^*(\bA^1,j_{\bG_m,!} \sfE_\infty)=0$ (we know that $H^*(\bA^1,j_{\bG_m,!}\Qbar_l) =0$ and that $\sfE_\infty$ is an extension of constant sheaves). 

Analogously we get that the fibre of the above complex at a point lying on $D_1$ and exactly $c$ other divisors is isomorphic to $H^*(\bA^1,j_{\bG_m,!}\sfE_\infty) \tensor H^*(\bG_m^{k-c-1},\bQ_\ell)=0$. So we have proven, the first part of the lemma.

This can be used to give an analogous description of $(\bR j_* \sfE_\infty)|_{D_I}$, because $D_1\cong \bA^{k-1}$ and so we can apply the same reasoning again. Consider the divisor $D_1^\prime\subset\bA^{k-1}\cong D_1$ and $j_{D_1^\prime}:\bG_m^{k-1}\hookrightarrow \bA^1\times \bG_m^{k-2}$.  Look at 
$$  \to j_{D_1^\prime,!} \bQ_\ell\to \bR j_{D_1^\prime,*} \bQ_\ell \to (\bR  j_{D_1^\prime,*} \bQ_\ell)|_{D_1^\prime} \map{[1]}.$$ 
Again $H^*(\bA^1\times \bG_m^{k-1},j_{D_1^\prime,!} \bQ_\ell)=0$, such that inductively
$$ (\bR j_* \bQ_\ell)|_{D_{1\dots l}} \cong \bR j_{l\dots k,*} \big( (\bR j_{1\dots l-1,*} \bQ_\ell)|_{D_{1\dots l-1}}\big) .$$ \hfill $\square$

\begin{kor}\label{spur-der-whittaker-garbe} For an arbitrary curve $C$, let $\sfE$ be a rank $n$ local system with indecomposable unipotent ramification at a finite set of points $S\subset C$.  Let  $I\subset \{1,\dots,k\}$ and let $\overline{D}_{I,p}\subset \Coh_{0,S}^{\underline{1}}$ be the substack defined by $(\phi^{i,p}=0)_{i\in I}$ (i.e. for $C=\bA^1$ this is the substack defined by $D_I\subset \bA^1$) and denote by $D^\circ_{I,p}$ the substack defined by $\phi^{(j,p)}\neq 0$ for $j\not\in I$. Then the following holds:\begin{enumerate}
\item For $0<|I|<k$ we have $H^*_c(\overline{D}_{I,p}, j_{!*}\sfE|_{\overline{D}_I})=0$.
\item Let $pr: \Coh_{0,S}^{\un{1}} \map{\cT^\bullet \mapsto \cT^{(0,S)}} \Coh_{0}^1$ be the projection. Then we have a canonical isomorphism $\bR pr_! j_{!*}\sfE \cong \sfE$.
\item $\bR pr_! (j_{!*}\sfE|_{D_{I,p}^\circ})\cong \sfE_{pr(D_{I,p})} [|I|-1]$ for any $I$ and again the isomorphism is canonical.
\end{enumerate}
\end{kor}
\noindent{\bf Proof:} In the special case $(C,S)=(\bA^1,\{0\})$ the corollary follows from Lemma \ref{restriction_to_D} which shows:
\begin{align}
j_{!*}\sfE_n|_{D_{1,\dots,l}^\circ}& \cong \sfE_p \tensor H^*(\bG_m,\Qbar_l)^{\tensor l-1} && \text{is constant and}\label{gm}\\
j_{!*}\sfE_n |_{\overline{D}_{1,\dots,l}}& = \bR j_{D\hookrightarrow \overline{D},*}j_{!*}\sfE_n|_{D_{1,\dots,l}}.
\end{align}
The first assertion of the corollary now follows from the K\"unneth formula and the fact that for $j:\bG_m \to \bA^1$ we have $H^*_c(\bA^1,\bR j_* \Qbar_l)=0$. This shows that $H^*_c(D_{I},m^*\sfE_n)=0$, and by the spectral sequence calculating the cohomology of a stack
from the cohomology of a presentation this gives the result for $H^*_c(\overline{D}_{I},\sfE_n)$.

The second assertion follows because by (1) the cohomology of $j_{!*}\sfE_n$ restricted to the complement of the section $\cT \mapsto \cT^\bullet$ vanishes. 
More precisely we get that in this case the canonical morphism $\sfE \to \bR pr_! j_{!*} \sfE$ given by the section $\Coh_0^1 \to \Coh_{0,S}^{\un{1}}$
is an isomorphism.

To prove (3) we note that $H^*(\bG_m,\Qbar_\ell)\cong H^*_c(\bG_m,\Qbar_\ell)[1]$ and compare (\ref{gm}) with the Leray spectral sequence for 
$$\xymatrix{ [pt/\bG_m] \ar[r] \ar[dr]^{id} & [pt/\bG_m^{l}]\cong D_{I,p}^\circ \ar[d]^{pr|_P} \\
 & [pt/\bG_m]=pr(D_{I,p}^\circ).} $$
Thus again the claimed isomorphism is obtained from a morphism defined by the
canonical section $s:pr(D_I) \to D_I$. 

The general case follows from these calculations, because the statements are local in the \'etale
topology on  $\Coh^{\un{1}}_{0,C,S}$. Therefore it is a problem which is local in the \'etale 
topology on $C$, thus to check that the morphisms given above are isomorphisms we may assume that
$C=C^{sh}_p$ is strictly henselian and $S=\{p\}$, i.e. $(C,p)\cong (\bA^{1,h}_{\overline{\bF}_q},0)$.
In this case any irreducibly ramified sheaf on $\bA^{1,h}_{\overline{\bF}_q}$ is isomorphic to our sheaf $\sfE_n$. \hfill $\square_{\textrm{\tiny Corollary}}$   

\noindent{\bf Remark:} This corollary implies that for any parabolic torsion sheaf $\cT^\bullet$ of degree one $\tr(\Frob_\cT,j_{!*}\sfE)$ is the eigenvalue of the Hecke operator corresponding
to $\cT$ applied to the Whittaker function $W_\sfE$, i.e. $(-1)^{\text{codim}(p)} \tr(\Frob_p,j_*E_p)$ 
 
To end this section we will prove two more corollaries of the above calculations. First we
give a description of $j_{!*} \sfE_k$ for $k<n$:

\begin{kor}
\begin{enumerate}
\item For any $0 < k < m$ we have a distinguished triangle:
$$ \to j_{k,!*} \sfE_k \to j_{k,!*} \sfE_m \to j_{k,!} \sfE_{m-k}(-k) \map{[1]} $$
\item $j_{!*} \sfE_k \cong \tau^{<k}\bR j_{k+1\dots n,*}j_{k,!*}\sfE_k$ 
\end{enumerate}
\end{kor}
\noindent{\bf Proof:} The first part of the corollary has been proven above. We may also recover it
by comparing the triangle from Proposition \ref{e_unendlich} for $\sfE_k$ with the one for $\sfE_m$.

The second part follows from Lemma \ref{eindimlemma} by induction on $k$: We have an exact triangle
$$  \to j_{k-1,!*} \sfE_{k-1} \to j_{k-1,!*} \sfE_k \to j_{k-1 \dots 1,!} \bQ_\ell(-k+1) \to $$
If we apply $j_{!*,k}$ to this, we get that $j_{k!*} \sfE_k \to \tau^{<k}\bR j_{k,*}j_{k-1\dots 1,!}\bQ_\ell(-k+1)$
induces a surjective map on the $(k-1)$-th cohomology sheaves, because by induction the $(k-1)$-th 
cohomology of $j_{k,!*} \sfE_{k-1}$ vanishes and the cohomology of the other two terms are 
isomorphic.
Now we can apply $\tau^{< k}\bR j_{k+1 \dots n,*}$ to the resulting triangle and by the Lemma 
\ref{eindimplus} we get again that the induced map to the $k-1$-th cohomology of the right term is 
surjective.
\hfill $\square$  

Finally we note that there is a -- perhaps surprising --  analogue of Corollary \ref{spur-der-whittaker-garbe} for
the tensor product $j_{!*}\sfE_n \tensor j_{!*} \sfE_{n+k}$ which will be needed later on:

\begin{kor}\label{pr-zwei} Let $(C,S)$ be a curve together with a finite set of points and let
$E_n, E_{n+k}$ be local systems of rank $n$ and $n+k$ on $C-S$, with indecomposable unipotent ramification
at all points in $S$. 

Let $pr: \Coh_{0,S}^{\un{1}} \to \Coh_{0}^1$ be the map 
forgetting the $n-$step parabolic structure of the torsion sheaves and denote by
$j:\Coh_{0,C-S}^1 \to \Coh_{0,C}^1$ the inclusion.

Then $$\bR pr_! (\cL_{\sfE_n}^1 \tensor \cL_{\sfE_{n+k}}^1) = j_*(\sfE_n\tensor \sfE_{n+k}).$$
\end{kor}

\noindent{\bf Proof:} To prove the corollary we have to show that $\bR pr_!(\cL_{\sfE_n}\tensor \cL_{\sfE_{n+k}})$
is the (middle) extension of its restriction to $\Coh_{0,C-S}^1$. This is a local problem on $C$, thus we
may assume as before that $(C,S)=(\bA^1,\{0\})$ and that $\sfE_n$ and $\sfE_{n+k}$ are the unipotently
ramified sheaves on $\bG_m$ defined at the beginning of this section. 

In this case note that 
$$j_*(\sfE_n\tensor \sfE_{n+k})\cong \oplus_{i=0}^{n-1} j_*\sfE_{2(n-i)+k}(-i).$$ 
This is just the Jordan decomposition for a tensor product (see for example \cite{Fulton-Harris}, Exercise 11.11).

We use the filtration 
$$j_{!*} \sfE_{n-1} \to j_{!*} \sfE_{n} \to j_{!*}j_{n-1\dots 1,!} \bQ_\ell(-n+1)$$
and prove by induction that the corollary also holds if we replace $\sfE_n$ by $\sfE_{n-i}$.
 
By Lemma \ref{eindimlemma} $j_{!*}j_{n-1\dots 1,!} \bQ_\ell(-n+1)$ is represented by the complex
$$ \bQ_\ell \to \oplus \bQ_\ell|_{D_i} \to \dots \to \oplus_{|I|=n-1} \bQ_\ell|_{D_I}. $$
And by Corollary \ref{spur-der-whittaker-garbe} we know 
$\bR pr_! ((j_{!*} \sfE_{n+k}) \tensor \bQ_\ell|_{D_I}) = 0$ for $0< |I| < n.$

Therefore 
\begin{eqnarray*}
\bR pr_! ( j_{!*} \sfE_{n+k} \tensor j_{!*}j_{n-1\dots 1,!} \bQ_\ell(-n+1)) & =& \bR pr_! (j_{!*} \sfE_{n+k} \tensor \bQ_\ell(-n+1))\\
 & = & j_* \sfE_{n+k}(n-1). 
\end{eqnarray*} 

Using induction we thus get an exact triangle
$$ \to j_*(\sfE_{n-1}\tensor \sfE_{n+k}) \to \bR pr_! (\cL_{\sfE_{n}}\tensor \cL_{\sfE_{n+k}}) \to j_* \sfE_{n+k}(-n+1) \map{[1]} $$
proving that the middle term is a perverse sheaf, which is the middle extension of its restriction to $\Coh_{0,C_S}^1$.
\hfill $\square$

\subsection{A Hecke property on $\Coh^{\un{d}}_{0,S}$}\label{section_hecketorsion}

Consider as before $S\subset C$ and a rank $n$ local system $\sfE$ on $C-S$ with indecomposable unipotent ramification at $S$. To reduce the number of
constants we will assume that we are looking at $n$-step parabolic sheaves 
(it would be sufficient to assume that $\rank(\sfE)\geq \textrm{length of structure}$).
 
Using Definition \ref{Heckeoperatoren} of the generalized {\em Hecke operators} the aim of this section is to prove:
\begin{satz}\label{hecketorsion} $\cL_\sfE^d$ is a {\em Hecke eigensheaf} on $\Coh_{0,S}^{d,\dots,d}$, i.e. for non-negative degrees $\un{d}=\un{d}^\prime+\un{d}^\pprime$ we have 
$$H^{\underline{d}^\pprime,\underline{d}^\prime}_0 \cL_\sfE^d = \left\{ \begin{array}{ll} \cL_\sfE^{d_1}\boxtimes \cL_\sfE^{d_2} & \textrm{if } \un{d}^\prime,\un{d}^\pprime \textrm{ are constant} \\ 0 & \textrm{otherwise.} \end{array}\right.$$
\end{satz}

To prove this, we need an analogue of Laumon's description of the Whittaker sheaf $\cL_\sfE^d$:  

Let $\widetilde{\Coh}^{\underline{d}}_{0,S}$ be the stack classifying parabolic torsion sheaves on $(C,S)$ together with a complete flag of subsheaves: 
$$\widetilde{\Coh}^{\underline{d}}_{0,S}(T):=\langle \cT_d^\bullet \supset \cT_{d-1}^\bullet\supset \dots \supset \cT_1^\bullet | \cT^\bullet_i \in \Coh^{\underline{i}}_{0,S}(T) \rangle$$
We use the diagram
$$\xymatrix{
{\widetilde{\Coh}}^{\underline{d}}_{0,\bA^1,0}\ar[d]_{\forget_{\Flag}}\ar[r]^-{\gr} & \Pi_{i=1}^d \Coh^{\underline{1}}_{0,\bA^1,0} \\
 \Coh^{\underline{d}}_{0,\bA^1,0}}$$ 
to define the sheaf $\tilde{\cL}_{\sfE} := \bR \forget_{\Flag,*} \gr^* j_{!*}\sfE^{\boxtimes d}$ on $\Coh^{\underline{d}}_{0,S}$.
Note that the map $\forget_{\Flag}$ is projective but not small (nor semi-small) in general.

\begin{satz} For any decomposition $\underline{d}=\underline{d}^\prime+\underline{d}^{\prime\prime}$ we have
$$  H_0^{\underline{d}^\pprime,\underline{d}^{\prime}} \tilde{\cL}_{\sfE}^d = \left\{\begin{array}{cl} \oplus_{S_d/(S_{d^\prime}\times S_{d^{\prime\prime}})} \tilde{\cL}_{\sfE}^{d^\prime} \boxtimes \tilde{\cL}_{\sfE}^{d^{\prime\prime}} & \textrm{ if } \underline{d}^\prime = (d^\prime,\dots,d^\prime) \textrm{ is constant}\\ 0 & otherwise. \end{array}\right.$$  
\end{satz}

\noindent{\bf Proof:}  Extend the diagram used to define the Hecke-operators as follows:
{\small $$\xymatrix@C=2ex{
& & {\Hecke}^{\un{d}^\prime,\un{d}^\pprime}\times_{\Coh_{0,S}^{\underline{d}}} \widetilde{\Coh}_{0,S}^{\underline{d}}\ar[r]\ar[dl]\ar@/^3.6pc/[rrd]^{\tilde{\gr}_{\Ext}} \ar@/_2pc/[dll]_{\tilde{\gr}_{\Flag}} & {\Hecke}^{\un{d}^\prime,\un{d}^\pprime}  \ar[dl]_{\pi_{big}}\ar[dr]^{\pi_{\text{small}}\times quot} & \\
(\Coh_{0,S}^{\underline{1}})^d &\ar[l]\widetilde{\Coh}_{0,S}^{\underline{d}}\ar[r]^{\forget_{\Flag}}& \Coh_{0,S}^{\underline{d}} & & \Coh_{0,S}^{\underline{d}^\prime}\times \Coh_{0,S}^{\underline{d}^{\prime\prime}}
}$$}
Using the base change theorem for the proper map $\forget_{\Flag}$ , we see that $$H_0^{\underline{d}^\prime,\underline{d}^{\prime\prime}} \tilde{\cL}_{\sfE}^d=\bR \widetilde{\gr}_{\Ext,!} \tilde{\gr}_{\Flag}^* \cL_{\sfE}^{1,\boxtimes d}.$$

The fibre product ${\Hecke}^{\un{d}^\prime,\un{d}^\pprime}\times_{\Coh_{0,S}^{\underline{d}}} \widetilde{\Coh}_{0,S}^{\underline{d}}$ classifies
$$\langle \cT^{\prime\bullet} \subset \cT^{\bullet} \to \cT^{\pprime\bullet} \; , \; \cT_1^{\bullet} \subset \dots \subset \cT_{d-1}^{\bullet} \subset \cT^{\bullet} \rangle.$$

For every such collection of torsion sheaves we can pull back the filtration of $\cT^{\bullet}$ 
to $\cT^{\prime\bullet}$, and by fixing the degrees $\underline{d}^\prime_i$  of the resulting torsion sheaves we obtain a stratification of the above stack
$$ {\Hecke}^{\un{d}^\prime,\un{d}^\pprime}\times_{\Coh_{0,S}^{\underline{d}}} \widetilde{\Coh}_{0,S}^{\underline{d}} = \cup_{\underline{d}^\prime_i} ({\Hecke}^{\un{d}^\prime,\un{d}^\pprime}\times_{\Coh_{0,S}^{\underline{d}}} \widetilde{\Coh}_{0,S}^{\underline{d}})^{\underline{d}^\prime_i}, $$
where the substacks of ${\Hecke}^{\un{d}^\prime,\un{d}^\pprime}\times_{\Coh_{0,S}^{\underline{d}}} \widetilde{\Coh}_{0,S}^{\underline{d}}$ are defined as 
$$({\Hecke}^{\un{d}^\prime,\un{d}^\pprime}\times_{\Coh_{0,S}^{\underline{d}}} \widetilde{\Coh}_{0,S}^{\underline{d}})^{\underline{d}^\prime_i}
:=\left\langle \left( {\begin{array}{l} \cT^{\prime\bullet}\to \cT^\bullet \to \cT^{\pprime\bullet} \\ \cT_i^{\bullet}\subset\cT^\bullet \end{array} }\right) \; | \; \deg(\cT^{\prime\bullet}\cap \cT_i^{\bullet}) = \underline{d}^\prime_i \right\rangle .$$

{\em 1$^{st}$ case: ${\underline{d}^\prime_i}=(d_i^\prime,\dots,d_i^\prime)$ is constant for all $i$.}
In this case we have a commutative diagram:
$$\xymatrix@C=2ex{
({\Hecke}^{\un{d}^\prime,\un{d}^\pprime}\times_{\Coh_{0,S}^{\underline{d}}} \widetilde{\Coh}_{0,S}^{\underline{d}})^{\underline{d}^\prime_i} \ar[r]^-{\widetilde{\forget}_{\Ext}}\ar[d] & {\widetilde{\Coh}}_{0,S}^{\underline{d}^\prime}\times \widetilde{\Coh}_{0,S}^{\underline{d}^\pprime} \ar[dr]^{\forget_{\Flag}^\prime\times\forget_{\Flag}^\pprime}\ar[d] & \\
(\Coh_{0,S}^{\underline{1}})^{\times d} \ar[r]^-{\cong} & (\Coh_{0,S}^{\underline{1}})^{\times d^\prime} \times (\Coh_{0,S}^{\underline{1}})^{\times d^\pprime} &\Coh_{0,S}^{\underline{d}^\prime} \times \Coh_{0,S}^{\underline{d}^\pprime}. 
}$$
By Lemma \ref{exakte_sequenzen} the map $\widetilde{\forget}_{\Ext}$ is smooth, the fibres being generalized affine spaces.
These are of dimension $0$, since both stacks are smooth of dimension $0$, thus
$$\bR \tilde{\gr}_{\Ext,!}
(\tilde{\gr}_{\Flag}^* (\cL_{\sfE}^1)^{\boxtimes{d}})|_
{({\Hecke}^{\un{d}^\prime,\un{d}^\pprime}\times_{\Coh_{0,S}^{\underline{d}}} \widetilde{\Coh}_{0,S}^{\underline{d}})^{\underline{d}^\prime_i}} = \tilde{\cL}_{\sfE}^{\underline{d}^\prime} \boxtimes \tilde{\cL}_{\sfE}^{\underline{d}^\pprime}. $$

{\em 2$^{nd}$ case: ${\underline{d}^\prime_i}$ not a constant sequence for some $i$.}

Let $\text{Flag}^{(\underline{d}_i^\prime)}$ be the stack, classifying torsion sheaves with a flag of subsheaves of degree $(\underline{d}_i^\prime)$. Then we can still factor the restriction of $\tilde{\gr}_{\Ext}$ to the corresponding stratum into
$$({\Hecke}^{\un{d}^\prime,\un{d}^\pprime}\times_{\Coh_{0,S}^{\underline{d}}} \widetilde{\Coh}_{0,S}^{\underline{d}})^{\underline{d}^\prime_i} \map{\widetilde{\forget}_{\Ext}} \text{Flag}^{(\underline{d}_i^\prime)}\times \text{Flag}^{(\underline{d}_i^\pprime)} \to \Coh_{0,S}^{\underline{d}^\prime} \times \Coh_{0,S}^{\underline{d}^\pprime}
$$
{\bf Claim:} $\bR \widetilde{\forget}_{\Ext,!} \tilde{\gr}_{\Flag}^* \cL_{\sfE}^{1,\boxtimes d} = 0.$ 

As in the first case the map $\widetilde{\forget}_{\Ext}$ is smooth, and the fibres are generalized affine spaces:
For a fixed point $(\cT^{\underline{d}_i^\prime},\cT^{\underline{d}_i^\pprime}) \in \text{Flag}^{(\underline{d}_i^\prime)}\times \text{Flag}^{(\underline{d}_i^\pprime)}$ the fibre of $\widetilde{\forget}_{\Ext}$ over this point consists of extensions
$$\xymatrix{
  {\cT}_{1}^{\prime\bullet} \ar@{^(->}[r]\ar@{^(->}[d] & \dots \ar@{^(->}[r] & {\cT}_{d-1}^{\prime\bullet}\ar@{^(->}[d] \ar@{^(->}[r] & {\cT}^{\prime\bullet}\ar@{^(->}[d] \\
  {\cT}^{\bullet}_{1} \ar@{^(->}[r]\ar@{->>}[d] & \dots \ar@{^(->}[r] & {\cT}_{d-1}^\bullet \ar@{->>}[d] \ar@{^(->}[r] & {\cT}^{\bullet}\ar@{->>}[d]\\
  {\cT}_{1}^{\pprime\bullet} \ar@{^(->}[r]& \dots \ar@{^(->}[r] & {\cT}_{d-1}^{\pprime\bullet}\ar@{^(->}[r] & {\cT}^{\pprime\bullet}
}$$
Let $\gr_i \cT^{\prime\bullet} := \cT^{\prime\bullet}_i /\cT^{\prime\bullet}_{i-1}$. Then we may
factor $\widetilde{\forget}_{\Ext}$ into
$$({\Hecke}^{\un{d}^\prime,\un{d}^\pprime}\times_{\Coh_{0,S}^{\underline{d}}} \widetilde{\Coh}_{0,S}^{\underline{d}})^{\underline{d}^\prime_i} \map{\gr_{\Ext}} \Pi_{i=1}^d \Ext(\gr_i \cT^{\pprime\bullet},\gr_i \cT^{\prime\bullet})\to \text{Flag}^{(\un{d}_i^\prime)}\times \text{Flag}^{\un{d}_i^\pprime},$$
where $\Ext(\gr_i \cT^{\pprime\bullet},\gr_i \cT^{\prime\bullet})$ is the generalized vector bundle
over $\text{Flag}^{(\underline{d}_i^\prime)}\times \text{Flag}^{(\underline{d}_i^\pprime)}$ classifying extensions
of the filtration quotients.
Furthermore Lemma \ref{exakte_sequenzen} shows that $\gr_{\Ext}$ is 
a generalized affine space bundle, which can be factored into maps with 
fibres $\Ext(\gr_i \cT^{\pprime \bullet},\cT_{i-1}^{\prime\bullet})$.

Since $\widetilde{\gr}_{\Flag}$ also factors through $\gr_{\Ext}$, the sheaf $\widetilde{\gr}_{\Flag}^*\tilde{\cL}_\sfE^d$ is 
constant on the fibres of $\gr_{\Ext}$ and thus by 
the K\"unneth formula it is sufficient to prove that for $d=1$ and any non-trivial decomposition 
$$\underline{d}=(1,\dots,1)= \underbrace{(\epsilon_1,\dots,\epsilon_n)}_{=:\underline{d}^\prime} + \underbrace{(1-\epsilon_1,\dots,1-\epsilon_n)}_{=:\underline{d}^\pprime}$$
 we have $H_0^{\underline{d}^\prime,\underline{d}^\pprime} \cL^1_{\sfE} = 0$. But here we can apply the calculation of $\cL^1_{\sfE} |_{D_I}$ given in Corollary \ref{restriction_to_D} to establish the claim. 

Now we have shown that $H_0^{\underline{d}^\prime,\underline{d}^{\prime\prime}} \tilde{\cL}_{\sfE}^d$ 
has a filtration such that the subquotients are isomorphic to the sheaves 
$\tilde{\cL}_{\sfE}^{d^\prime}\boxtimes\tilde{\cL}_\sfE^{d^\pprime}$. 
Furthermore we know that over the substack where 
$supp(\cT^{\prime\bullet})\cup supp(\cT^{\pprime\bullet})$ consists of $d$ distinct points, 
this extension splits. 
The proof of the following lemma will only use this fact to show that all these sheaves are perverse sheaves which are 
the middle extension of their restrictions to any open subset. 
Therefore the filtration splits globally.
\hfill $\square$

\begin{lemma}\label{symmetrische_gruppe_operiert}
The complex $\tilde{\cL}_{\sfE}^d=\bR \forget_! \gr^* ((\cL_{\sfE}^1)^{\boxtimes d})$ is a perverse sheaf which is the intermediate extension of its restriction to $\Coh_{0,C-S}^d$:
$$  \tilde{\cL}_{\sfE}^d=\bR \forget_! \gr^* ((\cL_{\sfE}^1)^{\boxtimes d} = j_{!*} \tilde{\cL}_{\sfE}^d|_{\Coh_{0,C-S}^d}.$$
In particular, it carries a natural action of the symmetric group $S_d$ and
$$\cL_{\sfE}^d=(\tilde{\cL}_{\sfE}^{d})^{S_d}.$$
\end{lemma}
Again we denoted by $j:\Coh_{0,C-S}^1\hookrightarrow \Coh_{0,S}^{\un{1}}$ the inclusion.

\noindent{\bf Proof of Lemma \ref{symmetrische_gruppe_operiert}:}
By Laumon's results \cite{Laumon_premiere_construction} we know that the restriction of $\tilde{\cL}_{\sfE}$ to $\Coh_{0,C-S}^d$ is indeed a perverse sheaf which is the middle extension of its restriction to every open subset. 

Since the question is local on $\Coh_{0,S}^{\underline{d}}$ we may assume that our local system 
$\sfE$ is pure. Then $\cL_{\sfE}^1$ is pure (it is irreducible and perverse) and thus, by Deligne's theorem (\cite{Deligne_Weil_II}, 6.2.6) $\widetilde{\cL}_{\sfE}^d$ is 
also pure. Therefore we may apply the Decomposition Theorem (\cite{faisceaux_pervers}, 5.4.6) 
to decompose $\tilde{\cL}_{\sfE}^d = j_{!*} j^*\tilde{\cL}_\sfE^d \oplus \text{rest}^d$. 

We prove the lemma by induction on $d$.  Assume that $\text{rest}^k=0$ for all $k<d$. (By definition of $\cL_{\sfE}^1$ the statement is true for $d=1$.) 

By the induction hypothesis and the fact that the restriction of $\widetilde{\cL}_\sfE^d$ to $\Coh_{0,C-S}$ is perverse we furthermore know that $\supp(\text{rest}^d)\subset \langle \cT^\bullet |\supp(\cT)= p\in S\rangle$.
The preceding proposition shows a Hecke property of $\tilde{\cL}_\sfE^d$ and this implies in particular that $H_0^{\un{i},(d)-\un{i}} rest^d=0$ for all $\un{i}>0$. 

Choose $\cT^\bullet\in \supp(\text{rest}^d)$ such that the degree of a maximal indecomposable summand of $\cT^\bullet$ is maximal.
And write $\cT^\bullet= \cO^\bullet_{\frac{i}{n}p}(\frac{j}{n}p) \oplus \cT^{\prime\bullet}$, such that $\cO^\bullet_{\frac{i}{n}p}(\frac{j}{n}p)$ is a direct summand of maximal degree (this is possible by Lemma \ref{maximale_summanden}). Note that $\cT^\bullet\not\cong \cO^\bullet_{dp}$ since the latter sheaf has a unique filtration. Now define $\un{d}^\prime := \deg(\cT^{\prime \bullet})$ and look at the fibre $F$ of the 
Hecke-correspondence $\Hecke_{0}^{\un{d}^\prime,(d)-\un{d}^\prime}$ over the point $(\cT^{\prime\bullet},\cO_{\frac{i}{n}p}(\frac{j}{n}p))\in \Coh_{0,S}^{\un{d}^\prime} \times \Coh_{0,S}^{\un{d}-\un{d}^\prime}$. 
Then $\cT^\bullet$ is the only sheaf contained in $\supp(\text{rest}^d)\cap F$, because every non-trivial extension of the two sheaves contradicts our maximality assumption (again by Lemma \ref{maximale_summanden}). 

Therefore if $\text{rest}^d|_{\cT^\bullet}\neq 0$ then $H_0^{\un{i},\un{d}-\un{i}} \text{rest}^d \neq 0$, contradicting our 
assumption that all the $\tilde{\cL}^k_{\sfE}$ are irreducible perverse sheaves for $k<d$.\hfill $\square$

\noindent{\bf Proof of Proposition \ref{hecketorsion}:} This now follows from the above lemma by taking $S_d-$invariants in the Hecke property of $\widetilde{\cL}_{\sfE}^d$. \hfill$\square_{\text{Proposition}}$


\section{The sheaf $\sfF_{\sfE,!}^{n}$ corresponds to the function $\Phi(W_\sfE)$}

The aim of this section is to explain the relation between the function $\tr_{\sfF_{\sfE,!}^n}$ and Shalika's definition of $\Phi(W_\sfE)$.
As in the case of unramified local systems the problem to compare the two functions stems from the fact that the interpretation of Laumon's diagram in terms of adeles does not immediately correspond to the definition of $\Phi$. The main ingredient needed to solve this problem is to introduce an analogue of Drinfeld's compactification as defined in \cite{FGV3}. This moduli space is on the one hand related to the fundamental diagram and on the other hand its points have a simple adelic description. All this follows easily from \cite{FGV3}.

However, to prove that the function
$\tr_{\sfF_{\sfE,!}^n}$ is indeed a non-zero multiple of the function $\Phi(W_{\sfE})$, we can not copy the proof of \cite{FGKV}, since this argument uses results on the affine Grassmannian for which we do not know analogous statements for the affine flag manifold. We will use an elementary approach instead. This yields an inductive argument to calculate the function $\tr_{\sfF_{\sfE,!}^n}$ on a subset which is sufficiently big to conclude the proof of our main theorem once we have calculated this function for $n-1$.
We will then give a calculation for $n\leq 2$. 

\subsection{An analogue of Drinfeld's compactification}

First we rewrite the inductive definition of $\sfF_{\sfE,!}^n$ as in the appendix of \cite{Laumon_premiere_construction} and \cite{FGV3}:

Denote by $\langle\OmegaExt\subset \cE^\bullet\rangle$ the stack classifying
$$ \langle \OmegaExt\subset\cE^\bullet\rangle(T):= 
\left\langle \begin{array}{l} 
\cE^\bullet \in Bun_{n,S}^{\un{d},\text{good}}(T), \cJ_i^\bullet\in Bun_{i,S}(T) \\
\cJ^\bullet_1\subset \cJ^\bullet_2\subset\dots\subset\cJ^{\bullet}_n\subset \cE \\
\cJ^\bullet_i/\cJ^\bullet_{i-1} \map{\cong} \Omega^{\bullet,n-i}  
\end{array}\right\rangle$$

We may define maps
\begin{align*} 
 \quot:\langle\OmegaExt\subset \cE\rangle & \to Coh_{0,S}^{\un{d}}\\
(\cE^\bullet,\cJ_i^\bullet) & \mapsto \cE^\bullet/\cJ^\bullet_n\\
\ext : \langle\OmegaExt\subset \cE\rangle & \to  \Pi_{i=1}^{n-1} \Ext_{\text{para}}(\Omega^{n-1-i,\bullet},\Omega^{n-i,\bullet}) \map{\sum} \bA^1\\
(\cE^\bullet,\cJ_i^\bullet) & \mapsto  \sum_{i=1}^{n-1}(\Omega^{\bullet,n-i} \to \cJ^\bullet_{i+1}/\cJ^\bullet_{i-1} \to \Omega^{\bullet,n-i-1})\\
\forget : \langle\OmegaExt\subset \cE\rangle & \to \Hom^{\text{inj}}_n\\
 (\cE^\bullet,\cJ_i^\bullet)& \mapsto (\cJ^\bullet_1 \hookrightarrow \cE^\bullet).
\end{align*}

Then by definition of $\sfF_{\sfE,!}^n$ we have 
$$\sfF_{\sfE,!}^n = \bR \forget_! (\quot^*\cL_\sfE \tensor \ext^* \psi)[c],$$
where $c$ is the dimension of the fibres of $\forget$.

\noindent\bem We have an adelic description of the points of the stack $\langle\OmegaExt\subset \cE\rangle$:
$$ \langle\OmegaExt\subset \cE\rangle(\bF_q) \subset \sfN_n(k(C)) \backslash \sfN_n(\bA) \times_{\sfN(\cO)} \GL_n(\bA) /(\GL_n(\cO_{C-S})\times \Iw_{S}) $$ 
We will not need this (it is the same as in \cite{FGKV}, Section 3) but note, that this is not the set
which is used in the definition of the function $\Phi(W_\sfE)$. 

To define a moduli space whose points will be a subset of $$\N_n(k(C))\backslash \GL_n(\bA) / \GL_n(\cO_{C-S}\times \Iw_S)$$ we argue as in \cite{FGV3} and define a moduli space classifying
parabolic vector bundles together with a full flag of subspaces of the generic fibre of the bundle, satisfying
some regularity condition: 

For a parabolic vector bundle $\cE^\bullet$ we denote by $\bigwedge^k \cE^\bullet$ its $k-$th exterior power, which is defined as the sequence of vector bundles 
$$ \dots \to \bigwedge^k \cE^{i,p} \to \bigwedge^k \cE^{i+1,p} \to \dots$$ 
Analogously denote for parabolic bundles $\cE_1^\bullet\tensor \cE_2^\bullet$ the tensor product
taken componentwise, together with the natural maps.

\begin{definition}{(Drinfeld's compactification) } 
The stack $\Omega\text{--Pl\"ucker}$ classifies:
$$ \OmegaP(T):= \left\langle \begin{array}{l}
\cE^\bullet\in Bun_{n,S}^{\un{d}}(T),\\
s_{1}:\Omega^{\bullet,n-1} \hookrightarrow \cE^\bullet,\\
s_{i}:\Omega^{\bullet,n-1} \tensor \dots \tensor \Omega^{\bullet,n-i} \hookrightarrow \wedge^{i} \cE^\bullet,\\
s_n:  \Omega^{\bullet,n-1} \tensor \dots \tensor \Omega^\bullet \tensor \cO^\bullet \hookrightarrow \wedge^n \cE^\bullet\\
\textrm{s. th. the } s_i \textrm{ satisfy the Pl\"ucker relations}
\end{array} 
\right\rangle.$$
\end{definition}
Recall that the Pl\"ucker relations are given by the condition that over the generic point of $C$ the maps $s_i$ define a full flag of subspaces of one (or equivalently all)  $\cE^{(i,p)}$. 
 
If all the $s_i$ are maximal embeddings (i.e. if the $s_i$ are maximal
embeddings in every degree $(j,p)$), then the $s_i$ define a 
full flag of $\cE^\bullet$ at every point of the curve, i.e. the $s_i$ define a
full flag of subbundles of $\cE^\bullet$.  

Therefore the points of this stack have a simple description:

\begin{lemdef}
The stack $\Omega\text{--Pl\"ucker}$ has a stratification by locally closed substacks indexed
by  degrees of parabolic divisors $\un{d_1},\dots,\un{d_n}$ (i.e. the coefficients of points in $S$ are allowed
to lie in $\frac{1}{n}\bZ$). 
The strata are given by:
\begin{align*}
(\OmegaP)_{(\un{d_1},\dots,\un{d_n})}(T) &:=  
\left\langle\begin{array}{l} 
\cE^\bullet\in Bun_{n,S}^{\un{d}}(T), D_i\in Div_{C,S}^{\un{d_i}},\\
s_{1}:\Omega^{\bullet,n-1}(D_1) \hookrightarrow \cE^\bullet,\\
s_{i}:\Omega^{\bullet,n-1}(D_1) \tensor \dots \tensor \Omega^{\bullet,n-i}(D_{n-i}) \hookrightarrow \wedge^{i} \cE^\bullet,\\
s_n:  \Omega^{\bullet,n-1}(D_1) \tensor \dots \tensor \cO^\bullet(D_n) \hookrightarrow \wedge^n \cE^\bullet\\
\textrm{such that the } s_i \textrm{ are maximal embeddings}\\
\textrm{and satisfy the Pl\"ucker relations,}\\
\textrm{and } \sum_{i=1}^k D_i \textrm{ is effective for all } 1\leq k \leq n.
\end{array}\right\rangle\\ 
&\cong \left\langle\begin{array}{l} 
\cE^\bullet \in Bun_{n,S}^{\un{d}}, \cJ_i^\bullet\in Bun_{i,S} \\
\cJ^\bullet_1\subset \cJ^\bullet_2\subset\dots\subset\cJ^{\bullet}_n=\cE \\
\cJ^\bullet_i/\cJ^\bullet_{i-1} \map{\cong} (\Omega^{n-i})^\bullet(D_i)  
\end{array}\right\rangle 
\end{align*}
For fixed parabolic divisors $D_1,\dots,D_n$ denote by $\OmegaP_{D_1,\dots,D_n}$ the corresponding 
substack of the above stack.
\end{lemdef}
\hfill $\square$

\begin{remark}\label{punkte-adelisch} The points of the stack $\OmegaP$ can be described as a subset:
$$ \OmegaP(\bF_q) \subset \N_n(k(C)) \backslash \GL_n^\Omega(\bA) / (\GL_n^\Omega(\cO_{C-S}) \times \Iw_S).$$
\end{remark}
\noindent {\bf Proof of Remark \ref{punkte-adelisch}:} This is the same as Weil's description of vector bundles (see also \cite{FGKV}). However to compare the function $W_\sfE$ with a sheaf on $\OmegaP$ we will need a precise form of the inclusion, therefore we will recall the construction of the map. 

Given a point $(\cE^{\bullet},s_i,D_i)\in\OmegaP_{D_1,\dots,D_n}$ we define an element of $\GL_n^\Omega(\bA)$ as follows: Let $N:=-(n-1)^2$ be the shift in the definition of $\Omega^{\bullet,n-1}$. 
(Note that if all $D_i=0$ then the bundle $\cE^{(N,S)}$ is equipped with a filtration with subquotients $\Omega^{\tensor n-i}(-(i-1)S)$.)

Recall that in \ref{GLomega} we have chosen an identification of $\GL_n(\bA)$ with $\GL_n^\Omega(\bA)$, i.e. we decided to use $\oplus_{i=0}^{n-1}\Omega^i$ as standard bundle instead of the trivial one. 
 
Denote by $\eta$ the generic point of $C$ and choose an isomorphism $f_\eta: \oplus_{i=0}^{n-1} \Omega^{\tensor i}_\eta \map{\cong} \cE^{(n-1,S)}_\eta$ such that the image of $\oplus_{i=n-j}^{n-1} \Omega^{\tensor i}_\eta$ is the subspace defined by $(s_i)_{i\leq j}$. 

Further, for $p\in C-S$ choose a trivialization $f_p: \oplus_{i=0}^{n-1} \Omega^{i}\map{\cong} \cE^{N,S}\tensor \widehat{\cO}_p$ again compatible with the filtration induced by the $s_i$. Then $f_p^{-1}\circ f_\eta\in\GL_n^{\Omega}(K_p)$   will be an element of the form $N_p \cdot \diag(d_{n,p},\dots,d_{1,p})$, where $N_p$ is a unipotent upper triangular matrix and the second term is a diagonal matrix such that the valuations of the entries are given by the $p-$part of the divisors $D_i$. 

For $p\in S$ we have to choose an isomorphism $f_p: \oplus_{i=0}^{n-1}\Omega^{\tensor i} \tensor \widehat{\cO}_p\map{\cong} \cE^{(N,S)}\tensor \widehat{\cO}_p$ compatible with the filtration of the stalk $\cE^{(N,S)}\tensor k(p)$. Thus we have to choose $f_p$ such that the induced map $\oplus_{i=0}^j \Omega^{\tensor i}\tensor k(p) \to \cE^{(N,S)}\tensor k(p)$ factors through $\ker\big(\cE^{(N,S)}\tensor k(p) \to (\cE^{(N+j,S)}\tensor k(p))\big)$. 

Again define $f_p^{-1}\circ f_\eta\in \GL_n(K_p)$. To describe this element, let $D_i = (d_{i} + \frac{k_i}{n}) p + D_i^\prime$ with $p\not \in \supp(D_i^\prime)$ and $0\leq k_i < n$, and choose a local parameter $\pi_p$ at $p$.  Then $f_p^{-1}\circ f_\eta (\Omega^{\tensor n-1})$ is contained in the $\widehat{\cO}_p$-submodule $\pi_p^{d_1+1}(\oplus_{j=0}^{k_1-1} \Omega^{\tensor j}) \oplus \pi_p^{d_1} (\oplus_{j=k_1}^{n-1} \Omega^j)$. Analogously the image of $f_p^{-1}\circ f_\eta (\Omega^{\tensor n-2})$ is contained in the subspace generated by $\pi_p^{d_2}(\oplus_{j=0}^{k_2} \Omega^j) \oplus \pi_p^{-1+d_2} (\oplus_{j=k_2+1}^{n-1}\Omega^j)$, etc. (We will only need this for $n=2$.)

Note that in this way we get an element of $\GL_n(K_p)$ for which we have calculated the value of the Whittaker function in Proposition \ref{formel}. In particular the shift in the definition of $\Omega^{i,\bullet}$ assures that the support of the Whittaker function is the subset of $\OmegaP$ where $D_1\leq D_2 \leq \dots \leq D_{n-1}$.   
\hfill $\square$

Note furthermore that we can factorize the forgetful map:
$$ \forget:\langle \OmegaExt\subset \cE\rangle \map{\forget_{Tor}} \Omega\textrm{--Pl\"ucker} \to \Hom(\Omega^{n-1,\bullet},\cE). $$

Therefore to prove that $\tr_{\sfF_{\sfE,!}^n}=\Phi(W_\sfE)$ is the correct function (up to a scalar), it is sufficient to prove that $\tr_{\bR \forget_{Tor,!} (\quot^*\cL_{\sfE}\tensor ext^* \psi)} =W_\sfE$ (up to a scalar). 
Our first aim is to show that the left hand side of the last equation is an element of the space of Whittaker functions (Proposition \ref{whittaker_eigenschaft}).

We denote by $\OmegaExt_{D_1,\dots,D_n}$ the preimage $\forget_{Tor}^{-1}(\OmegaP_{D_1,\dots,D_n})$.

Note that whenever we have $0\leq D_1\leq D_2 \leq \dots \leq D_n$, we can define a sheaf $\psi_{D_1,\dots,D_n}$ on $\OmegaP_{D_1,\dots,D_n}$ via
\begin{multline*} \ext_{D_1,\dots,D_n} : \OmegaP_{D_1,\dots, D_n} \to \Pi_{i=1}^{n-1} \Ext_{\text{para}} (\Omega^{\bullet,n-i-1}(D_{i+1}),\Omega^{\bullet,n-i}(D_i))\\  \to \Pi \Ext_{\text{para}} (\Omega^{\bullet,n-i-1},\Omega^{\bullet,n-i}) \to \bA^1 
\end{multline*}
$$ \psi_{D_1,\dots,D_n} := \ext_{D_1,\dots,D_n}^* \psi.$$
Let $d_j:=\min_{i,p} \un{d_j}^{(i,p)}$ and denote by $\lfloor D_j \rfloor$ the biggest divisor smaller than the parabolic divisor $D_j$. Then we also have a map
$$ div: \OmegaP_{D_1\dots D_n} \to C^{(d_1)} \times C^{(d_2-d_1)} \times \dots \times C^{(d_n-d_{n-1})},$$
sending $(D_1,\dots,D_n)$ to $(\lfloor D_1 \rfloor, \lfloor D_2 \rfloor - \lfloor D_1 \rfloor, \dots , \lfloor D_n\rfloor - \lfloor D_n-1 \rfloor)$.

The aim of this section is to prove:

\begin{satz}\label{whittaker_eigenschaft}
Let $D_1,\dots,D_n$ be (parabolic)-divisors and assume $0\leq D_1$. Then:\begin{enumerate}
\item If $D_i\not\leq D_{i+1}$ for some $i$, then 
$$\bR\forget_{Tor,!} (quot^*\cL_\sfE\tensor ext^*\sfL_\psi))|_{\Omega\textrm{--P\"ucker}_{D_1,\dots,D_n}}  = 0.$$
\item If $0\leq D_1\leq D_2\leq \dots \leq D_n$, then there is a sheaf $\sfW_\sfE$ on $C^{(d_1)} \times \dots \times C^{(d_n-d_1)}$ and a constant $c$ such that
$$\bR\forget_{Tor,!} (quot^*\cL_\sfE\tensor ext^*\sfL_\psi))|_{\Omega\textrm{--P\"ucker}_{D_1,\dots,D_n}} = \Psi_{D_1,\dots,D_n} \tensor div^* \sfW_\sfE[-2c](-c).$$ 
The sheaf $\sfW_\sfE$ and the constant $c$ depend on the parabolic degrees of the $(D_j)$ and will be defined explicitly in the proof.
                 \end{enumerate}
\end{satz}

\noindent{\bf Proof:} We may assume that all $D_i$'s are effective, since otherwise the fibres of $\forget_{Tor}$ above $\OmegaP_{D_1,\dots,D_n}$ are empty.

The vanishing assertion follows essentially from the following simple observation, which simply is the geometric reformulation of the corresponding statement proven in \ref{formel}.
Assume that $D_1\not\leq D_{2}$, and let $D$ be the effective part of $(D_1-D_{2})$.
Let furthermore $\cE^\bullet_{n-1}(=\cE^\bullet/\Omega^{\bullet,n-1}(D_1))$ be a parabolic vector bundle of arbitrary degree and
assume that we are given a commutative diagram
$$\xymatrix{
  \Omega^{\bullet,n-1} \ar@{^(->}[d]\ar@{^(->}[r]^-{s} & {\cE}^\bullet_{n-1} \oplus \Omega^{\bullet,n-1}_{D_1}\ar@{->>}[d]\\
  \Omega^{\bullet,n-1}(D_{2}) \ar@{^(->}[r]^-{\overline{s}}_-{\text{maximal}} & {\cE}^\bullet_{n-1}
}$$
where the vertical maps are the canonical ones, i.e. we are given a lift $s$ of 
the maximal embedding $\overline{s}$.

\noindent{\bf Claim:} For any $\phi:\Omega^{\bullet,n-2} \to \Omega^{\bullet,n-1}_D \hookrightarrow \Omega^{\bullet,n-1}_{D_1}$
the two extensions
$$ 0 \to \Omega^{\bullet,n-2} \map{s} \cE^\bullet_{n-1} \oplus \Omega^{\bullet,n-1}_{D_1} \to \cF^\bullet \to 0$$
and
$$ 0 \to \Omega^{\bullet,n-2} \map{s+\phi} \cE^\bullet_{n-1} \oplus \Omega^{\bullet,n-1}_{D_1} \to \cF^{\prime\bullet} \to 0$$
are isomorphic.

This holds because locally at $\supp(D)$ we can lift $\phi$:
$$\xymatrix{
\Omega^{\bullet,n-2}(D_{2})/\Omega^{\bullet,n-2}(-D)\ar@{..>}[r]^-{\tilde{\phi}} & \Omega^{\bullet,n-1}_{D_{1}}\\
\Omega^{\bullet,n-2}/\Omega^{\bullet,n-2}(-D) \ar[r]^-{\phi}\ar@{^(->}[u] & \Omega^{\bullet,n-1}_D\ar@{^(->}[u]
}$$

And --- choosing a splitting of $\overline{s}$ at $\supp(D)$ --- this can be used to define an automorphism
of $\cE^\bullet_{n-1}\oplus \Omega^{\bullet,n-1}_{D_1}$ mapping $s$ to $s+\phi$. \hfill$\square_{\textrm{\tiny Claim}}$

This tells us that for $D_1\not\leq D_2$ the fibres of $\forget_{Tor}$ over $\OmegaP_{D_1,\dots,D_n}$ carry an action of the additive group $\Hom(\Omega^{\bullet,n-1},\Omega^{\bullet,n-1}_D)$ and on the orbits of this action the sheaf $\ext^*\sfL_\psi$ is non-trivial whereas $\quot^*\cL_\sfE^d$ is constant. Therefore the cohomology of the fibres with values in $\ext^* \sfL_\psi\times quot^*\cL_\sfE^d$ vanishes. 

To state this more precisely and to apply this argument inductively we need some notation:  

The forgetful map: $\OmegaExt_{D_1,\dots,D_n} \to \OmegaP_{D_1,\dots,D_n}$
factors through the stacks $\Omega_{D_1,\dots,D_k} \Ext_{D_{k+1},\dots,D_n}$, classifying
$$ \left\langle\begin{array}{l}
\cJ_1^\bullet \subset \dots \subset \cJ_n^\bullet\subset \cE^\bullet\\
\cJ_i^\bullet/\cJ^\bullet_{i-1} \map{\cong} \Omega^{n-i,\bullet} (D_i) \textrm{ for } i \leq k\\
\cJ_i^\bullet/\cJ^\bullet_{i-1} \map{\cong} \Omega^{n-i,\bullet}       \textrm{ for } i > k\\
\textrm{such that } \cJ_k^\bullet\subset\cE^\bullet \textrm{ is a maximal embedding} 
\end{array}\right\rangle.$$

And let $\forget_{D_1\dots D_i}=\forget_{D_i}\circ \dots \circ \forget_{D_1}$ be the corresponding forgetful maps.
In case $k=1$, the above stack classifies maximal embeddings $\Omega^{\bullet,n-1}(D_1)\hookrightarrow \cE^\bullet$
together with an iterated extension of $\cO^\bullet,\dots,\Omega^{\bullet,n-2}$ (namely $\cJ_n^\bullet/\cJ^\bullet_1$) 
contained in the quotient bundle $\cE^\bullet/\Omega^{\bullet,n-1}(D_1)$.
Therefore the fibre of $\forget_{D_1}$ consists of the different lifts of the iterated extensions to 
$\cE^\bullet/\Omega^{\bullet,n-1}\cong \cE^{\bullet}/\Omega^{\bullet,n-1}(D_1)\oplus\Omega^{\bullet,n-1}_{D_1}$,
more precisely, $\forget_{D_1}$ is a $\Hom(\cJ_n^\bullet/\cJ_1^\bullet,\Omega^{\bullet,n-1}_{D_1})$-torsor.
 
In particular $\Hom(\Omega^{\bullet,n-2},\Omega^{\bullet,n-1}_{D_1})$ acts on the fibres
of $\forget_{D_1}$, and by the above claim for any point in a fibre, the extension
$$\Omega^{\bullet,n-2} \to \cE^\bullet/\Omega^{\bullet,n-1} \to (\cE^\bullet/\Omega^{\bullet,n-1})/\Omega^{\bullet,n-2}$$
does not change if we modify it by $s\in \Hom(\Omega^{\bullet,n-2},\Omega^{\bullet,n-1}_{D})$.
On the other hand the induced extension of $\Omega^{\bullet,n-2}$ by $\Omega^{\bullet,n-1}$ 
does change by the residue of $s$. 
Therefore for $i=1$ we get that $\bR \forget_{D_1,!}(\Psi\tensor quot^*\cL_E^d)=0$ if 
$D_1 \not \leq D_2$.

Assume now that $0\leq D_1 \leq D_2$. Note that for a point $(\cJ_i^\bullet\subset\cE^\bullet)\in \OmegaExt_{D_1,\dots,D_n}$ the quotient $\cE^\bullet/\cJ^\bullet_N$ is also equipped with a filtration and morphisms identifying the subquotients with $\Omega^{\bullet,n-i}_{D_i}$. Denote by $\Ext(\cO^\bullet_{D_n},\dots,\Omega^{\bullet,n-1}_{D_1})$ the stack classifying torsion sheaves with such a filtration. Then we get a diagram
$$\xymatrix{
\Omega\Ext_{D_1,\dots,D_n} \ar[d]^{pr_{\Fib}}\ar[drr]^-{quot}\ar[dr]^-{qext}                & &\\
\textrm{Fib}               \ar[d]^{pr_1}\ar[r]_-{pr_2}          & \Ext(\cO^\bullet_{D_n},\dots,\Omega^{\bullet,n-1}_{D_1})\ar[d]^{gr_{D_1}} \ar[r]_-{pr_{big}} & Coh_{0,S}^{\un{d}}\\
\Omega_{D_1}\Ext_{D_2,\dots,D_n} \ar[r]^-{qext_{n-2}}          & C^{(d_1)} \times \Ext(\cO^\bullet_{D_n},\dots,\Omega^{\bullet,n-2}_{D_2}). &
}$$
Here $\Fib$ is the fibre product making the lower square cartesian and 
the occurring maps are the natural maps.
Note that $pr_{\Fib}$ is a vector bundle map and that its fibres are isomorphic to $\Hom(\cE^\bullet/(\Omega^{\bullet,n-1}(D_1)),\Omega^{\bullet,n-1}_{D_1})$. In particular, the map $ext$ factors through $pr_{\Fib}$ and thus $ext^*\sfL_\psi\tensor quot^*\cL_\sfE^d$ descends to $\Fib$.
Moreover we can again define a map $ext_{D_1}:\Omega_{D_1}\Ext_{D_2,\dots,D_n}\to \bA^1$ as the composition of the natural maps:
{\small $$\Omega_{D_1}\Ext_{D_2,\dots,D_n}  \to \Ext^1(\Omega^{\bullet,n-2}(D_2),\Omega^{\bullet,n-1}(D_1))\times \Pi_{i=2}^{n-1} \Ext^1(\Omega^{\bullet,n-i-1},\Omega^{\bullet,n-i}) \to \bA^1$$}
And finally since $D_1\leq D_2$ we get an exact sequence
$$ \Hom(\Omega^{\bullet,n-2}(D_2),\Omega^{\bullet,n-1}_{D_1}) \map{0} \Hom(\Omega^{\bullet,n-2},\Omega^{\bullet,n-1}_{D_1}) \tto \Ext^1(\Omega^{\bullet,n-2}_{D_2},\Omega^{\bullet,n-1}_{D_1}) \to 0,$$
and therefore the residue map $\res:\Hom(\Omega^{\bullet,n-2},\Omega^{\bullet,n-1}_{D_1})\to \bA^1$
factors through $\Ext(\Omega^{\bullet,n-2}_{D_2},\Omega^{\bullet,n-1}_{D_1})$. Thus $\res^*\sfL_\psi$
defines a sheaf on the latter space, which we can pull back to a sheaf $\Psi_{12}$ on
$\Ext(\cO^\bullet_{D_n},\dots,\Omega^{\bullet,n-1}_{D_1})$.

By the additivity of $\Psi$ on $\bA^1$ we get  
$$ext^*\sfL_\psi\cong \forget_{D_1}^* ext_{D_1}^* \sfL_\psi \tensor \qext^* \Psi_{12}.$$
Thus we can calculate our sheaf $\bR \forget_{D_1,!} (\ext^*\sfL_\psi \tensor \quot^*\cL_\sfE^d)$ as follows:
\begin{multline*}
\bR \forget_{D_1,!} (\ext^*\sfL_\psi \tensor \quot^*\cL_\sfE^d)  \\ \cong  \bR \forget_{D_1,!} (\forget_{D_1}^* \ext_{D_1}^* \sfL_\psi \tensor \qext^* \Psi_{12} \tensor \quot^*\cL_\sfE^d)\\
           \cong  \ext_{D_1}^* \sfL_\psi \tensor \qext_{n-1}^* (\bR \gr_{D_1,!} (\cL_\sfE^d\tensor \Psi_{12}))[-2 c_1](-c_1), 
\end{multline*}
where $c_1=\dim(\Hom(\cE^\bullet/(\Omega^{\bullet,n-1}(D_1)),\Omega^{\bullet,n-1}_{D_1}).$ 

We can inductively apply the same considerations to the maps $\forget_{D_i}$ to get:
\begin{enumerate}
\item $\bR \forget_{Torsion,!} (ext^* \sfL_\Psi \tensor quot^*{\cL_\sfE^d}) = 0$ unless $0 \leq D_1 \dots \leq D_n$.
\item If we have $0\leq D_1 \leq \dots \leq D_n$, then we may define a sheaf $\Psi_{Tor}$ on 
the stack $\Ext(\cO^\bullet_{D_n},\Omega^\bullet_{D_{n-1}},\dots,\Omega^{\bullet,n-1}_{D_1})$ 
as the tensor product of the pull-backs of the sheaves $\Psi_{i,{i+1}}$ on $\Ext(\Omega^{\bullet,n-i-1}_{D_{i+1}},\Omega^{\bullet,n-i}_{D_i})$.
\item Denote by $gr$ the natural map
$$ gr: \Ext(\cO^\bullet_{D_n},\Omega^\bullet_{D_{n-1}},\dots,\Omega^{\bullet,n-1}_{D_1})\to C^{(d_1)}\times \dots \times C^{(d_n-d_{n-1})}$$
and define  $\sfW_\sfE := \bR gr_! (\cL_\sfE^d\tensor \Psi_{Tor})$. Then
 $$\bR \forget_{Tor,!} (ext^*\sfL_\psi \tensor quot^*\cL_{\sfE}^d) \cong \Psi_{D_1,\dots,D_n} \tensor div^* \sfW_\sfE[-2c](-c),$$
where $c=\sum_{i=1}^{n-1} c_i$ and $c_i=\dim \Hom(\cE^\bullet/\cJ^\bullet_{i}, \Omega_{D_i}^{\bullet,n-i})$.
\end{enumerate} 
\hfill $\square_{\textrm{\tiny Proposition}}$
\pagebreak[2]
 
To compare the trace function of $\bR \forget_{Tor,!} \ext^*\sfL_\psi \tensor \quot^*\cL_\sfE^d$ and $W_\sfE$ we therefore only need to calculate the trace function of $\sfW_\sfE$. By construction it is sufficient to do this in the case that all $D_i$'s are supported at a single point $p$ since the trace of $\sfW_\sfE$ will be a product of these. We may also assume that $p\in S$, because for $p\not\in S$ we can use the calculations for unramified local systems \cite{FGV3} (note however that a calculation similar to the one we do below (Lemma \ref{W2}) could be applied for $p\not\in S$ as well). 

\subsection{Calculation in the case rank $=2$}

\begin{lemma}\label{W2}
Consider sheaves with $2-$step parabolic structure at $S=\{p\}\in C$. Denote by $\lambda_\sfE:= \tr(Frob_p,j_*\sfE)$. 
Then for any $d\in \bN$ and $k\in \frac{1}{2}\bN$ with $0\leq k \leq d-k$ we have
$$ \sum_{e\in \Ext^1(\cO^\bullet_{(d-k)p},\Omega^\bullet_{k p})} Tr(Frob_e,\Psi_{Tor}\tensor \cL_\sfE^d) = \left\{ \begin{array}{rl} q^{2k} \lambda_\sfE^d & \textrm{for } k\in\bN\\[1ex]
-q^{2k} \lambda_E^{d} & \textrm{for } k \in \frac{1}{2}+\bN.                           \end{array} \right.$$  
\end{lemma}

\noindent{\bf Remark:} The Hecke property of $\cL_\sfE^d$ (Lemma \ref{hecketorsion}) implies that
$$ \sum_{e\in \Ext^1(\cO^\bullet_{(d-k)p},\Omega^\bullet_{k p})} \tr(Frob_e,\cL_\sfE^d) = \left\{ \begin{array}{rl} 
                    q^{k} \lambda_\sfE^d & \text{for } k\in\bN\\[1ex]
                             0           & \text{for } k \in \frac{1}{2}+\bN. 
        \end{array} \right.
$$
(Note that the set $\Ext^1$ used above differs from the stack $\underline\Ext$ by some automorphisms, whereby we obtain the factor $q^k$ in the above formula.)
Note further that we have $\Ext^1(\cO^\bullet_{(d-k)p},\Omega^\bullet_{k p})\supset \Ext^1(\cO^\bullet_{(d-k)p},\Omega^\bullet_{(k-1) p})\supset \dots \supset 0$, where the inclusions are given by push outs. Therefore for any $e$ which lies in the subset $\Ext^1(\cO^\bullet_{(d-k)p},\Omega^\bullet_{(k-l) p})-\Ext^1(\cO^\bullet_{(d-k)p},\Omega^\bullet_{(k-l-1)p})$ the corresponding parabolic torsion sheaf $\cT^\bullet$ is isomorphic to
$$ \cT^\bullet \cong \left\{ \begin{array}{cl} \Omega^\bullet_{(d-l)p} \oplus \cO^\bullet_{l p} & \textrm{for } 0\leq l < k \textrm{ and } k\in \frac{1}{2}+\bN \\[1ex]
\Omega^\bullet_{(d-l-\frac{1}{2})p} \oplus \cO^\bullet_{(l+\frac{1}{2})} & \textrm{for }0\leq l < k  \textrm{ and } k\in\bN_{>0}\\[1ex]
\cO^\bullet_{(d-k)p} \oplus \Omega^\bullet_{k p} & \textrm{if } e=0. \end{array}\right.
$$ 
It might be helpful to write this out in the simplest cases:
If $d=1, k=\frac{1}{2}$ then $\Omega^\bullet_{\frac{1}{2}p} = (\to \Omega_p \to 0 \to)$ and $\cO^\bullet_{\frac{1}{2}p}= ( \to 0 \to \cO_p \to)$, thus the extensions are of the form $(\to \Omega_p \map{0} \cO_p \map{\phi_1})$, and $\phi_1$ is nonzero if the extension is non trivial.
Similarly for $d=2,k=1$ the nontrivial extension is of the form $(\Omega_{2p} \to \Omega_p \oplus \cO_p \to)$.

As in the unramified situation $\tr(\Frob_{\cO^\bullet_{dp}},\cL_\sfE^d)=\lambda_\sfE^d = \tr(\Frob_{\Omega_{dp}^\bullet},\cL_\sfE^d)$. Therefore the above gives a recursion relation for $\tr(\Frob_{\cO^\bullet_{d-k}\oplus \Omega^{\bullet}_{k}},\cL_{\sfE}^d) =: L_\sfE^d(k)$. (Note that $L^d_{\sfE}(k)=L^d_{\sfE}(d-k)$, since the two torsion sheaves differ only by a shift.)
\begin{equation}\label{rekursion} L^d_\sfE(k) = \left\{ \begin{array}{rl} 
               q^k \lambda_\sfE^d - (q-1) \sum_{i=0}\limits^{k-1} q^i L^d_\sfE(k-\frac{1}{2}-i) & \textrm{for } k\in \bN \\
               -(q-1) \sum\limits_{i=0}^{k-\frac{1}{2}} q^i L_\sfE(k-\frac{1}{2}-i) & \textrm{for } k\in \frac{1}{2} + \bN
                           \end{array} \right.
\end{equation}
Note further that this recursion relation does not depend on the rank of $\sfE$.

\noindent{\bf Proof of Lemma \ref{W2}:} By induction on $k$ (for $k=0$ there is nothing to show):
Note that in our formula most of the summands cancel out since $\sum_{x\in \bA^1} \Psi(x) =0$:
\begin{eqnarray*}
\lefteqn{\sum_{e\in \Ext^1(d-k,k)} \Psi(e) L_\sfE(e)  =  L^d_\sfE(k) - L^d_\sfE(k-\frac{1}{2})}\\
&\!\!\!\! \eqweil{(\ref{rekursion})} &\!\!
   \left\{ \begin{array}{rl}
      q^k \lambda_\sfE^d - q L^d_\sfE(k-\frac{1}{2}) - (q-1) \sum\limits_{i=1}^{k-1} q^i L^d_\sfE(k-\frac{1}{2}-i) & \textrm{for } k\in \bN\\
      -q L^d_{\sfE}(k-\frac{1}{2}) -(q-1) \sum\limits_{i=1}^{k-\frac{1}{2}} q^i L^d_\sfE(k-\frac{1}{2}-i) & \textrm{for } k\in \1halb +\bN \\
   \end{array}\right.\\
&\!\!\!\! = & \!\! \left\{\begin{array}{ll}
       \bigg( q^k\lambda_\sfE^d - \sum\limits_{i=0}^{k-1} q^{i+1}\big(L^d_\sfE(k-\1halb-i) - L^d_\sfE(k-1-i)\big) & \\
       - \sum\limits_{i=0}^{k-2} q^{i+1}\big(L^d_\sfE(k-1-i) - L^d_\sfE(k-\frac{3}{2}-i)) - q^k L^d_\sfE(0)\big) \bigg)& \!\!\!\!\text{for } k\in \bN\\[2.5ex]
       -\sum\limits_{i=0}^{k-\1halb} q^{i+1} L_\sfE^d(k-\1halb-i)+\sum\limits_{i=0}^{k-\frac{3}{2}} q^{i+1} L_\sfE^d(k-\frac{3}{2}-i) & \!\!\!\!\text{for }k\in \1halb+\bN 
    \end{array}\right.\\
 &\!\!\!\!\eqweil{induct.} & \!\! \left\{\begin{array}{ll}
-\sum\limits_{i=0}^{k-1} q^{i+1}(-q^{2(k-i)-1} \lambda^d_\sfE) - \sum\limits_{i=0}^{k-2} q^{i+1}(q^{2(k-1)-2}\lambda_\sfE^d) & \textrm{for }k\in \bN\\[1ex]
-\sum\limits_{i=0}^{k-\1halb}q^{i+1}(q^{2(k-i)-1} \lambda_\sfE^d) - \sum\limits_{i=0}^{k-\frac{3}{2}} q^{i+1}(-q^{2(k-i-1)} \lambda_\sfE^d) & \textrm{for } k\in \1halb +\bN\end{array}\right. \\
&\!\!\!\!= & \!\!\left\{ \begin{array}{rl}
                 q^{2k} \lambda_\sfE^d & \textrm{for } k\in \bN\\[1ex]
                 -q^{2k} \lambda_\sfE^d & \textrm{for } k\in \1halb +\bN
                \end{array} \right.
\end{eqnarray*}
\hfill$\square_{\textrm{\tiny Lemma}}$
     
\begin{kor}\label{spurFE}
Let $\sfE$ be a local system on $C-S$ with indecomposable unipotent ramification at $S$ and denote by $\lambda_\sfE := \prod_{p\in S} \tr(Frob_p,j_*\sfE)$. \begin{enumerate}
\item If $\sfE$ is of rank {\bf $2$}, then for any point $x\in \OmegaP$ we have
$$ \tr(\Frob_x, \bR \forget_{Tor,!} (ext^*\sfL_\psi\tensor quot^*\cL_\sfE^d))=\lambda_\sfE \cdot q^{|S|} \cdot W_\sfE(x).$$
In particular for any point $\overline{x}\in \Hom_2^{\text{inj}}$ we have
$$ \tr(\Frob_{\overline{x}},\sfF_{\sfE,!}^2) = \lambda_\sfE \cdot q^{|S|}\cdot \Phi(W_\sfE)(\overline{x}).$$
\item If $\sfE$ is of rank {\bf $3$}, then for any point $x\in \OmegaP$ with $D_1=0$ we have
$$ \tr(\Frob_x, \bR \forget_{Tor,!} (\ext^*\sfL_\psi\tensor \quot^*\cL_\sfE^d))=\lambda_\sfE^{3} q^{3|S|} W_\sfE(x).$$
In particular, for any point $\overline{x}\in \Hom_3^{\text{inj}}$ corresponding to a maximal embedding $\Omega^{\bullet,2}\hookrightarrow \cE^\bullet$ we have
$$ tr(Frob_{\overline{x}},\sfF_{\sfE,!}^3) = \lambda_\sfE \cdot q^{3|S|} \cdot \Phi(W_\sfE)(\overline{x}).$$    \end{enumerate}
\end{kor}
\noindent{\bf Proof:} Comparing the above lemma with the calculation of $W_\sfE$ we get the first assertion. Note that since the power of $\lambda_\sfE$ appearing on either side of the equation depends only on the degree, we just have to compare these for the trivial bundle.

Similarily the power of $q$ only depends on the difference $D_1-D_2$. 

For the second assertion note that for a maximal embedding, the quotient sheaf $\cE^\bullet/\Omega^{\bullet,n-1}$ may be viewed as a bundle with $(n-1)$-step parabolic structure since the $n-$th morphism in the parabolic structure is an isomorphism. Thus for rank $3$ bundles we may apply the calculation given above.
\hfill$\square_{\textrm{\tiny Corollary}}$


\section{Constructing $\sfA_\sfE$ under the assumption $\sfF_\sfE^n=\sfF_{\sfE,!}^n$}
In this section we give a proof of the main Theorem \ref{ergebnis} under the additional assumption that $\sfF_\sfE^n=\sfF_{\sfE,!}^n$. Here the proofs are almost identical to the ones in the case of unramified local systems: First we show that the Hecke property for $\cL_\sfE^d$ implies that $\sfF_{\sfE,!}$ is a Hecke eigensheaf as well. The second step is to deduce from Lafforgue's theorem and the calculation of the previous section, that the function $t_{\sfF_{\sfE,!}^n}$ descends to a function on $Bun_{n,S}^d$. Therefore we can  argue as in \cite{FGV3} that the sheaf $\sfF_\sfE$ also descends to the space of parabolic vector bundles. The resulting sheaf $\sfA_\sfE$ inherits the Hecke property from $\sfF_{\sfE,!}$, and we show that this property implies the one stated in the theorem. 

\subsection{The Hecke operators on the ``fundamental diagram''}

We want to check that Laumon's arguments in \cite{Laumon_premiere_construction} carry over to our situation. 
We define operators analogous to the operators $H_k^{\un{i}}$ on the spaces occurring in the fundamental diagram (\ref{Teil1}). We start with $\Hom^{\text{inj}}_k$:  
$$
\xymatrix@C=1.5ex{ 
  & \left\langle { \begin{array}{c} {\cF}^{\prime\bullet} \subset {\cF}^\bullet \\ \Omega^{\bullet,k-1} \hookrightarrow {\cF}^\bullet \end{array}} \right\rangle \ar[dl]_{\pi_{\text{big}}} & 
\langle \Omega^{\bullet,k-1} \hookrightarrow \cF^{\prime\bullet} \subset \cF^\bullet \rangle \ar@{_(->}[l]_-{\iota} \ar[dr]^{\pi_{\text{small}}\times \quot}\\
\Hom^{\text{inj}}_k & & & \Hom^{\text{inj}}_k \times \Coh_{0,S}^{\underline{i}}
}$$
\begin{eqnarray*}
H^{\un{i}}_{k,\Hom^{\text{inj}}} : D^b(\Hom^{\text{inj}}_k) & \to & D^b(\Hom^{\text{inj}}_k\times \Coh_{0,S}^{\underline{i}}) \\
 \sfK & \mapsto & \bR (\pi_{\text{small}}\times \quot)_! \iota^* \pi_{\text{big}}^* \sfK 
\end{eqnarray*}                                                      
Analogously we define an operator 
$$H^{\un{i}}_{k,\Hom}: D^b(\Hom(\Omega^{\bullet,k-1}, \cF^\bullet)) \to D^b(\Hom(\Omega^{\bullet,k-1},\cF^{\prime,\bullet})\times \Coh_{0,S}^{\un{i}}).$$

We have the same on $\Ext^1_k$ (and $\Ext_k^{1,\text{good}}$):
$$\xymatrix@C=1.6ex{
 & \left\langle {\begin{array}{c} \cF^{\prime\bullet}\subset \cF^\bullet \\ \Omega^{\bullet,k-1}\to \cF_{k}^\bullet \to \cF^\bullet \end{array}}\right\rangle \ar[dl]_{\pi_{\text{big}}} \ar@{->>}[r]^-{p} & 
\left\langle {\begin{array}{c} \cF^{\prime\bullet}\subset \cF^\bullet \\ \Omega^{\bullet,k-1}\to \cF_{k}^{\prime\bullet} \to \cF^{\prime\bullet} \end{array}}\right\rangle
\ar[dr]^{\pi_{\text{small}}\times \quot} & \\
\Ext^1_k & & & \Ext^1_k\times \Coh_{0,S}^{\underline{i}}
}$$
\begin{eqnarray*}
H^{\un{i}}_{k,\Ext^1} : D^b(\Ext^1_k) & \to & D^b(\Ext^1_k\times \Coh_{0,S}^{\underline{i}}) \\
 \sfK & \mapsto & \bR  (\pi_{\text{small}}\times \quot)_! \bR p_! \pi_{\text{big}}^* \sfK 
\end{eqnarray*}

\subsection{The Hecke property of $\sfF_{\sfE,!}^k$}
We want to show that these Hecke operators commute with the functors used to construct $\sfF_{\sfE,!}^k$.
Let $i_{k,S}:= deg(\cT^{k,S})$ for any $\cT^\bullet \in \Coh_{0,S}^{\un{i}}$. And denote by $pr_{\Coh_{0,S}^{\un{d}}}:\Ext^1_1 \to \Coh_{0,S}^{\un{d}}$ the projection. 

\begin{satz}\label{hecke_und_konstruktion}
  For any $\un{d},\un{i}$ as above and any complex $\sfK$ we have:\begin{enumerate}
\item $ H^{\un{i}}_{0,\Ext^1} pr_{\Coh_{0,S}^{\un{d}}}^* \sfK = (pr_{\Coh_{0,S}^{\un{d}-\un{i}}}^* \times Id_{\Coh_{0,S}^{\un{i}}})^* H_0^{\un{i},\un{d}-\un{i}} \sfK[-2i_{0,S}](-i_{0,S})$
\item $ H^{\un{i}}_{k,\Ext^{1,\text{good}}} j_{\Ext}^* \sfK = j_{\Ext}^* H^{\un{i}}_{k,\Ext^1} \sfK.$
\item $ H^{\un{i}}_{k,\Hom} j_{\Hom,!} \sfK = j_{\Hom,!} H^{\un{i}}_{k,\Hom^{\text{inj}}}.$
\item $ H^{\un{i}}_{k,\Hom^{\text{inj}}} I^*\sfK = I^* H^{\un{i}}_{k-1,\Ext^1} \sfK.$
\item $ H^{i}_{k,\Ext^1} \circ \Four\sfK = \Four \circ H^{i}_{k,\Hom^{\text{inj}}} \sfK [-i_{k(n-1),S}](-i_{k(n-1),S}) $ 
\end{enumerate}  
\end{satz}

\noindent{\bf Proof:}\begin{enumerate}\item Write down the definition of the correspondences:
$$\xymatrix@C=1.5ex{
&  \left\langle {\begin{array}{c} \cT^{\prime\bullet}\subset \cT^\bullet \\ \cO^\bullet \to \cF_{1}^\bullet \to \cT^\bullet \end{array}}\right\rangle
\ar@{->>}[r]^-{p}\ar[dl]_{\pi_{\text{big}}} \ar[dr]^{pr_{\text{left}}}& \ar[d]_{pr_{\text{right}}}
\left\langle {\begin{array}{c} \cT^{\prime\bullet}\subset \cT^\bullet \\ \cO^{\bullet}\to \cF_{1}^{\prime\bullet} \to \cT^{\prime\bullet} \end{array}}\right\rangle
\ar[dr]^-{\pi_{\text{small}}\times \quot} & \\
\Ext^1_1\ar[d]^-{pr_{\Coh_{0,S}^{\un{d}}}} & &\langle \cT^{\bullet\prime}\subset \cT^\bullet \rangle \ar[dll]_{\forget}\ar[dr]^{gr}  &  \Ext^1_1\times \Coh_{0,S}^{\un{i}} \ar[d]^-{pr_{\Coh_{0,S}^{\un{d}-\un{i}}}\times Id}\\
\Coh_{0,S}^{\un{d}} & & & \Coh_{0,S}^{\un{d}-\un{i}}\times \Coh_{0,S}^{\un{i}}
}$$
The left- and right-hand squares are cartesian and $p$ is a vector bundle projection, therefore we get our claim:
\begin{align*}
\begin{split}
H_{0,\Ext^1}^{\un{i}} \pr_{\Coh_{0,S}^{\un{d}}}^* \sfK = \bR (\pi_{\text{small}}\times \quot)_! \bR p_! (\pi_{\text{big}}\circ pr_{\Coh_{0,S}})^* \sfK \\
\end{split}\\
\quad\quad &=  \bR (\pi_{\text{small}}\times \quot)_! \bR p_! ( \forget \circ \pr_{\text{right}} \circ p)^* \sfK \\
 &= \bR (\pi_{\text{small}}\times \quot)_! (\forget \circ \pr_{\text{right}})^* \sfK[-i_{0,S}](-i_{0,S}) && \text{as } p\text{ is a bundle}\\
 &= (pr_{\Coh_{0,S}^{\un{d}-\un{i}}}\times Id) ^* ( \bR gr_! \forget^* \sfK )[-i_{0,S}](-i_{0,S})&&\text{by base-change}
\end{align*}

\item This holds, because extensions of good sheaves by torsion sheaves are good.
\item By definition.

\item This is again true by definition, because there is a canonical isomorphism
$$ \langle \Ext^1(\cF_{k-1}^\bullet,\Omega^{\bullet,k-1}),\cF_{k-1}^{\prime\bullet} \rangle  \cong \langle \Omega^{k-1} \hookrightarrow \cF_k^{\prime\bullet} \subset\cF^\bullet \rangle.$$
And thus we get an isomorphism of the diagrams defining the two Hecke functors.
\item Again Laumon's proof can be copied word by word, the only thing used is the compatibility 
of the Fourier transform with bundle maps: $\Four (\iota^* \sfK) = \bR p_! \Four \sfK [i_{k(n-1),S}](i_{k(n-1),S})$ (see \cite{Laumon_transformation} Thm 1.2.2.1 and 1.2.2.4). \hfill $\square_{\textrm{\tiny Proposition}}$
\end{enumerate}

\begin{kor}\label{hecke_eigenschaft}
The sheaf $\sfF_{\sfE!}^k$ is a Hecke eigensheaf on $\Hom^{\text{inj}}_k$, i.e.:
$$ H^{\un{i}}_{k,\Hom^{\text{inj}}} \sfF_{\sfE!}^k = \left\{ \begin{array}{ll} \sfF_{\sfE!}^k \boxtimes \cL_{\sfE}^i[-ki](-ki) & \textrm{if } \un{i} \textrm{ is constant}\\
0 & \textrm{otherwise.} \end{array}\right.$$     
\end{kor}
\noindent{\bf Proof:} By the above Proposition \ref{hecke_und_konstruktion} this follows from the Hecke property of $\cL_\sfE^i$ (Proposition \ref{hecketorsion}). 
\hfill $\square_{\textrm{\tiny Corollary}}$

\subsection{Comparison of the Hecke operators and the generalized Hecke operators}

In the same way as in \cite{FGV3}, Proposition 8.4 we want to show that for some sheaves on $\Bun_{n,S}$  the eigensheaf property with respect to $H_n^{\un{d}}$ implies the eigensheaf property for $H^d_n$. To do this we need to note some general properties of the maps $\pi_{\text{small}}$ and $\pi_{\text{big}}$ used in the definition of the operators $H^{\un{d}}_k$. 

Fix a degree $\underline{d}=(d^{(j,p)})$ of parabolic sheaves, and let $\underline{i}$ some positive degree. We have defined a diagram
$$\xymatrix{ 
  & Hecke_{k}^{\underline{d},\underline{i}} \ar[dl]_{\pi_{\text{big}}}\ar[dr]^{\pi_{\text{small}}\times \quot} & \\
\Coh_{k,S}^{\un{d}} & & \Coh_{k,S}^{\underline{d}-\un{i}} \times \Coh_{0,S}^{\underline{i}} 
}$$
Denote further $\Coh_{k,S}^{\underline{d},\leq \un{i}}:=\langle \cF^\bullet\in \Coh_{k,S}^{\underline{d}}\; | \; \length(\torsion(\cF^\bullet)\leq \un{i} \rangle$. Then we have:

\begin{remark}\label{hecke_abbildungs_eigenschaft}\begin{enumerate}
\item The map $\pi_{\text{big}}$ is representable and projective.
\item The restriction of $\pi_{\text{big}}$ to the pre-image of $\Coh_{n,S}^{\underline{d},\leq \un{i}}$ is smooth.
\item The map $\pi_{\text{small}}\times \quot$ is a generalized vector bundle, in particular it is smooth.
\item The map $\pi_{\text{small}}$ is smooth. 
\item 1., 4. and the second part of 3. are true for the analogous maps defined by replacing $\Coh_{n,S}^{\un{d}}$ and  $\Coh_{k,S}^{\un{d}-\un{i}}$ by $\Bun_{k,S}^{\un{d}}$ and $\Bun_{k,S}^{\un{d}-\un{i}}$ respectively. 
\end{enumerate}\end{remark}
\noindent{\bf Proof:} \begin{enumerate}
\item The fibres of $\pi_{\text{big}}$ are closed subschemes in the scheme $\prod Quot_{\rank\; n \hfill \atop \deg\; d^{(j,p)}}(\cF^{(j,p)})$
which is projective (see \cite{FGA}). 
\item As in \cite{FGV3} we want to show that for every fibre the tangent space at a point $(\cE^\bullet \hookrightarrow \cF^\bullet \tto \cT^\bullet:=\cF^\bullet/\cE^\bullet)$
is isomorphic to $\Hom_{para}(\cE^\bullet,\cT^\bullet)$. In Lemma \ref{dim_Hom_und_Ext} we have 
shown that this space is of constant dimension. And therefore, since $\pi_{\text{big}}$ is representable and 
projective, it must be smooth.

To see that the tangent space at a point of the fibre over $\cF^\bullet_0$ is indeed $\Hom(\cE^\bullet_0,\cT^\bullet_0)$, note that a point in the tangent space is a deformation to $k[\epsilon]/(\epsilon^2)$:
$$ \xymatrix{ {\cE}^\bullet \ar[r]^-{\psi}\ar[d] &  {\cF}_0^\bullet\tensor_k k[\epsilon]/\epsilon^2 \ar[r] \ar[d]& {\cT}^\bullet \ar[d] \\
{\cE}_0^\bullet \ar[r] & {\cF}^\bullet_0 \ar[r] & {\cT}_0^\bullet }$$
But in this case $\cE^\bullet \cong \cE_0 \times_{\cF^\bullet_0} \cF^\bullet_0 \tensor k[\epsilon]/\epsilon^2 \cong \cE_0 \tensor k[\epsilon]/\epsilon^2$. And therefore the choices of $\psi$ are given by $\Hom(\cE_0^\bullet,\cT^\bullet_0)$, as claimed.

\item Furthermore the map $\pi_{\text{small}}\times \quot$ is the projection from the generalized vector bundle 
$$\bV ( \bR pr_{12,*} \Hom(pr_{23}^*\cT_{univ}^\bullet,pr_{13}^*\cF_{univ}^\bullet)) \to \Coh_{n,S}^{\underline{d}-\un{i}} \times \Coh_{0,S}^{\underline{i}},$$
 where $pr_{jl}$ are the projections from $\Coh_{k,S}^{\un{d}-\un{i}}\times \Coh_{0,S}^{\un{i}}\times C$ on the $j$ and $l$-th factors, and $\cT^\bullet_{univ}$ and $\cF_{univ}^{\bullet}$ are the universal bundles on $\Coh_{0,S}^{\un{i}}\times C$ and $\Coh_{k,S}^{\un{d}-\un{i}}\times C$ respectively. 
\item By 3. we only need to note that $\Coh_{0,S}^{\underline{i}}$ is a smooth stack (Lemma \ref{coh-glatt}).
\item Since $\Bun_{k,S}^{\un{d}}\subset \Coh_{n,S}^{\un{d}}$ is open the maps are still smooth. The restriction of $\pi_{\text{big}}$ is still projective because 
subsheaves of vector bundles on curves are automatically vector bundles.
\end{enumerate}\hfill $\square_{\text{Remark}}$

Recall that $\overline{\Coh}_{0,S}^{\underline{d}}:=\Coh_{0,S}^{\un{d}}/(\textrm{central }\bG_m\textrm{-automorphisms})$. And the diagram defining the Hecke operators $H^{\un{1}}$ can be written as
$$\xymatrix{
& \bP \Ext^1(\cT^\bullet,\cE^{\prime\bullet}) \ar[dl]_{\pi_{\text{big}}}\ar[dr]^{\pi_{\text{small}}\times \overline{\quot}} & \\
\Bun_{0,S}^{\underline{d}}\incl{j_{\Bun}}\Coh_{0,S}^{\underline{d}} & & \Bun_{0,S}^{\underline{d}} \times \overline{\Coh}_{0,S}^{\underline{1}}  
}$$
 
\begin{satz}\label{normale-hecke-eigenschaft}
Assume that $A_\sfE$ is a Hecke eigensheaf for $\sfE$ on $\Coh_{n,S}^{\underline{d}}$, such that $\bD A_\sfE$ is a Hecke eigensheaf for $\bD \sfE=: \sfE^\vee$.
Then $A_\sfE|_{\Bun_{n,S}^{\un{d}}}$ is an eigensheaf for $H^{\un{1}}$, i.e.
$$H^{\un{1}}_n A_\sfE|_{\Bun_{n,S}} = A_\sfE|_{\Bun_{n,S}} \boxtimes \overline{\cL}_\sfE^1[-n+1](-n+1)$$
\end{satz}
\noindent{\bf Proof:} 
Look at the Hecke correspondence restricted to $\Bun_{n,S}^{\underline{d}-\un{1}}$:
$$\xymatrix{
 & \Ext^1(\cT^\bullet,\cE^{\prime\bullet}) \ar[dl]_{\pi_{\text{big}}}\ar[dr]^{\pi_{\text{small}}\times \quot} & \\
\Coh_{n,S}^{\underline{d},\leq \un{1}} & & \Bun_{n,S}^{\underline{d}-\un{1}} \times \Coh_{0,S}^{\un{1}}
}$$
We know by Lemma \ref{hecke_abbildungs_eigenschaft} that $\pi_{\text{big}}$ is smooth over $\Coh_{n,S}^{\underline{d},\leq \un{1}}$, therefore:
\begin{eqnarray*}
 A_\sfE \boxtimes \sfE [-n+1](-n+1) & \eqweil{\bD A_\sfE \text{ eigensheaf}} & \bD H^{\un{1}} \bD A_\sfE|_{\Coh_{n,S}^{\underline{d},\leq \un{1}}}\\
 & \eqweil{\text{def. of } H} & \bD \bR (\pi_{\text{small}}\times \quot)_! \pi_{\text{big}}^* \bD A_\sfE\\
 & \eqweil{\pi_{\text{big}} \text{ smooth}} & \bD \bR (\pi_{\text{small}}\times \quot)_! \pi_{\text{big}}^! \bD A_{\sfE} [-2n](-n)\\
 & = & \bD \bR (\pi_{\text{small}}\times \quot)_! \bD \pi_{\text{big}}^* A_\sfE [-2n](-n)\\
 & = & \bR (\pi_{\text{small}}\times \quot)_* \pi_{\text{big}}^* A_\sfE[-2n](-n).
\end{eqnarray*}
In other words, for $A_\sfE$ we can replace $\bR (\pi_{\text{small}}\times \quot)_!$ by $\bR (\pi_{\text{small}}\times \quot)_*$ in the definition of the Hecke operators.

Let $\un{\epsilon}$ be some degree with entries $\epsilon^{(i,p)}\in\{ 0,1\}$, such that not all entries are equal to $0$ and not all entries are equal to $1$.

Then we know that $H^{\un{\epsilon}} A_\sfE = 0$, and this helps to prove:

\begin{lemma}Under the assumptions of \ref{normale-hecke-eigenschaft} 
the restriction of the sheaf $\sfA_\sfE$ to the stack $\Coh_{n,S}^{\un{d},\leq \un{\epsilon}}-\Bun_{n,S}^{\un{d}}$ 
is zero.
\end{lemma}

\noindent{\bf Proof:} Since $\pi_{\text{big}}^* A_\sfE$ is $\bG_m$-equivariant, we can apply Lemma \ref{brylinski-lemma} to the generalized vector bundle $\Ext^1(\cT^\bullet,\cE^\bullet)\to \Bun_{n,S}^{\un{d}-\un{\epsilon}} \times \Coh_{0,S}^{\un{\epsilon}}$
and get that $s_0^* \pi_{\text{big}}^* A_\sfE =\bR (\pi_{\text{small}}\times \quot)_* \pi_{\text{big}}^* A_\sfE = 0$. 
\hfill $\square_{\textrm{\tiny Lemma}}$

Now we can apply Lemma 8.5. of \cite{FGV3} to get that
$$ \bR (\pi_{\text{small}}\times \overline{\quot})_! \pi_{\text{big}}^* A_\sfE = A_\sfE \boxtimes \overline{\cL}_\sfE^1[-n+1](-n+1).$$

By the above lemma this establishes the proposition, because the stalk of $\sfA_\sfE$ is zero at sheaves $\cF^\bullet$ with $\un{0} < \deg(\torsion(\cF^\bullet)) < \un{1}$.
\hfill $\square_{\textrm{\tiny Proposition}}$
  
\begin{kor}\label{FEnichtnull}
Assume that $A_\sfE$ is a Hecke eigensheaf for $\sfE$ on $\Coh_{n,S}^{\un{d}}$, such that 
$\bD A_\sfE$ is a Hecke eigensheaf for $\bD \sfE:= \sfE^\vee$. Then the corresponding function 
$t_{A_\sfE}$ on $\Bun_{n,S}^{\un{d}}$ is an eigenfunction for the Iwahori-Hecke algebra.
\end{kor} 

\noindent{\bf Proof:} We just have proven the Hecke-property of the restriction of $A_\sfE$ to $\Bun_{n,S}$.
Therefore we only need to compare the result with the computation of $\cL_\sfE^1$ on 
$\overline{\Coh}_{0,S}^{\un{1}}$  (Lemma \ref{spur-der-whittaker-garbe}) and note that the 
Iwahori-Hecke-algebra at $S$ is generated by elements corresponding to the points of $\Coh_{0,S}^{\un{1}}$. 
\hfill$\square$
\subsection{Descent of the sheaf $\sfF_\sfE^n$}
\begin{satz}
Assume that we know that $\sfF_\sfE^n = \sfF_{\sfE,!}^n$, then Lafforgue's theorem implies
that $\sfF_\sfE^n$ descends to a Hecke eigensheaf on $\Bun_{n,S}^{\text{good}}$, 
and this sheaf can be extended to a non zero Hecke eigensheaf $\sfA_\sfE$ on $\Bun_{n,S}$.
\end{satz}

\noindent{\bf Proof:} By definition $\sfF_{\sfE}^n$ is an irreducible perverse sheaf and by our
assumption $\sfF_\sfE^n=\sfF_{\sfE,!}^n$ is a Hecke eigensheaf (by Corollary \ref{FEnichtnull}).

We first want to explain why the function $\Phi(W_\sfE)$ does not depend on the section $\Omega^{\bullet,n-1}\hookrightarrow \cE^\bullet$. On the one hand by Lafforgue's theorem \cite{lafforgue} there is a (cuspidal) Hecke eigenfunction on $\Bun_{n,S}(\bF_q)$ with eigenvalues given by $\tr_{\cL_\sfE}$ (\ref{formel}). On the other hand by Shalika's result (\cite{Shalika} Theorem 5.9) every Hecke eigenfunction on $\Hom^{\text{inj}}_n(\bF_q)$ is in the image of $\Phi$ and there is a unique such function in the Whittaker space. Therefore the function $f_\sfE = \Phi(W_\sfE)$ is the pull back of a function on $\Bun_{n,S}(\bF_q)$.

Assume for the moment that $n\leq 3$. In this case we know that restricted to the maximal embeddings $\Hom^{max}_n\subset \Hom^{\text{inj}}_n$ the function $\tr_{\sfF_{\sfE,!}^n}= \Phi(W_\sfE)$. In particular, this function is not identically zero on $\Hom^{max}_n$ and descends to $\Bun_{n,S}(\bF_q)$.
  
Thus we can apply a variant of the argument given in \cite{FGV3}: 
Since $\sfF_\sfE^n$ is an irreducible perverse sheaf, there is a constructible subset 
$V\subset \Hom^{max}(\Omega^{\bullet,n-1},\cE^\bullet)$ such that $\sfF_\sfE^n|_V$ is an irreducible local system
and $\sfF_\sfE^n=j_{V,!*}(\sfF_\sfE^n|_V)$. 

Assume that $\sfE$ is pure (this is allowed by Lafforgue's theorem).  Then $\sfF_\sfE^n$
is pure as well. Therefore the restriction of $\sfF_\sfE^n$ to $V$ is constant on the fibres over $Bun_{n,S}$,
because the trace of $\sfF_\sfE^n$ is constant on the fibres (for any extension $\bF_{q^n}$ of the base field). 

The fibres of $V$ are open subsets of projective spaces, therefore the two pull backs of $\sfF_\sfE$ to $\Hom^{max}\times_{\Bun_{n,d}}\Hom^{max}$ are irreducible and isomorphic, because the corresponding trace functions are the same. Since the two systems are irreducible, there is only one isomorphism of these sheaves which induces the identity on the points of the diagonal $V\subset V\times_{Bun_n} V$.  Hence $\sfF_{\sfE}^n|_V$ 
descends to a perverse sheaf $A_{\sfE,V}$ on $pr_{\Bun_{n,S}}(V)$. 

Further, since $\sfF_\sfE^n=j_{V,!*}(\sfF_\sfE^n|_V)$, we also know that 
$\sfF_\sfE^n=pr_{\Bun_{n,S}}^* j_{pr(V),!*} A_{\sfE,V}$, i.e. $\sfF_\sfE^n$ descends to a sheaf $\sfA_\sfE^{\text{good}}$ on 
$\Bun_{n,S}^{\text{good}}$.

Note that in particular we have shown that $\tr_{\sfF_{\sfE,!}}= \Phi(W_\sfE)$ on the whole of $\Hom^{\text{inj}}_n$. Therefore we may apply $\Phi^{-1}$ to see that the trace function of the sheaf $\bR \forget_{Tor,!}(\quot^*(\cL_\sfE^{d_0})\tensor ext^*\sfL_\psi)$ on $\OmegaP$ is equal to $W_\sfE$. This allows us to drop the temporary assumption that $n\leq 3$, because we can apply the argument of Lemma \ref{spurFE} to show that the trace of $\sfF_\sfE^n$ is equal to $\Phi(W_\sfE)$ on the space of maximal embeddings for $n\leq 4$, and this gives an inductive argument for all $n$.

To finish the proof of the theorem we only need to extend the resulting sheaf $\sfA_\sfE^{\text{good}}$
 to the whole of $\Bun_{n,S}$. Again this works as in \cite{FGV3}(Section 7.8): 
For  $q\in C-S$ (we might allow $q \in S$) the maps $\tensor \cO(-rq): \Bun_{n,S}^{d+r,\text{good}} \to \Bun_{n,S}^{d}$ 
are a covering of $\Bun_{n,S}$.
We define $\sfA_\sfE := \varinjlim_{r} (\tensor\cO(-rq))_* \sfA_{\sfE}^{\text{good}} \tensor (det(\sfE)|_q)^{-\tensor r}.$
The Hecke property of $\sfA_\sfE^{\text{good}}$ (together with the $S_2$-equivariance of the isomorphism
$H^1\circ H^1 \sfA_\sfE^{\text{good}} \cong \sfA_\sfE^{\text{good}} \boxtimes \sfE \boxtimes \sfE$) gives that 
this is a well-defined Hecke eigensheaf on $\Bun_n^d$.  

\hfill$\square$


\section{The analogue of the vanishing theorem for $n\leq 3$}

The aim of the last two sections of this article is to prove that our assumption $\sfF_\sfE^k=\sfF^k_{\sfE,!}$
holds for $k\leq 3$ (Proposition \ref{sauber}). To do so we need an analogue of the vanishing theorem in \cite{FGV3} which
is given below (Proposition \ref{verschwinden}):

For any $i\in\bZ_{>0}$ consider the {\em total Hecke-} or {\em averaging functor } $H_{\sfE,tot}^{-i}$ defined as follows ($\underline{i}:=(i,\dots,i)$):
$$\xymatrix{
& \langle \cE^{\prime\bullet}\subset \cE^{\bullet} \; | \; \cE^{\prime\bullet}\in \Bun_{n,S}^{\underline{d}-\underline{i}} \; , \; \cE^{\bullet}\in \Bun_{n,S}^{\underline{d}}  \rangle \ar[dl]_{\pi_{\text{small}}}\ar[dr]^{\pi_{\text{big}}}\ar[d]_{\quot}& \\
\Bun_{n,S}^{\underline{d}-\underline{i}} & \Coh_{0,S}^{\underline{i}} & \Bun_{n,S}^{\underline{d}} 
}$$

We set
\begin{eqnarray*}
H^{-i}_{\sfE,tot}: D^b(\Bun_{n,S}^{\underline{d}-\underline{i}}) & \to & D^b(\Bun_{n,S}^{\underline{d}}) \\
\sfK & \mapsto & H^{-i}_{\sfE,tot}\sfK := \bR\pi_{\text{big},!} (\pi_{\text{small}}^*\sfK \tensor quot^* \cL^{i}_\sfE).
\end{eqnarray*}

{\bf Remark:} This definition is used for any $\underline{d}=(d_{i,p})_{0\leq i<n,p\in S}$, 
therefore it includes the case of bundles with not necessarily full parabolic structure. 
In particular for $\underline{d}=(d)_{0\leq i<n \atop p\in S}$ the stack $\Bun_{k,S}^{\underline{d}}\cong \Bun_{k}^d$
is the stack of vector bundles without extra structure.

\begin{satz}\label{verschwinden}
Let $\sfE$ be a pure irreducible rank $n$ local system with indecomposable unipotent
ramification at $S$. Then for any $k<\min(3,n)$ and any (mixed) complex 
$\sfK \in D^b(\Bun_{k,S})$ we have
$$ H^{-i}_{\sfE,tot} \sfK = 0 \text{ for all } i>(2g-2)nk + |S| $$
\end{satz}
\noindent{\bf Proof:} (almost the same as in \cite{FGV3})
We use that, by induction we already know the proposition for all $k^\prime < k$.

{\em Reductions:} Without loss of generality, we may assume that $\sfK$ is a pure complex,
 because any mixed complex has a filtration with pure filtration quotients.

For a pure complex $\sfK$ the complex $H^{-i}_{\sfE,tot} \sfK$ is pure as well, 
because $H^{-i}_{\sfE,tot}\sfK = \bR\pi_{\text{big},!} (\pi_{small}^*\sfK \tensor quot^* \cL^{i}_\sfE)$
and $\pi_{small}$ is smooth (Lemma \ref{hecke_abbildungs_eigenschaft}), therefore 
$\pi_{small}^*$ preserves purity (smooth morphisms are locally acyclic, i.e. 
$\pi_{small}^*=\pi_{small}^![d](d)$). The same is true for $quot^*$ and finally $\pi_{\text{big}}$ is
proper (Lemma \ref{hecke_abbildungs_eigenschaft}), therefore Deligne's theorem (\cite{Deligne_Weil_II}, 6.2.6) implies that $\bR \pi_{\text{big},*}=\bR \pi_{\text{big},!}$ also preserves purity.

Furthermore, a pure complex $H^{-i}_{\sfE,tot} \sfK$ is zero if and only if the 
associated function $\tr_{H^{-i}_{\sfE,tot}\sfK}$ on $\bF_{q^l}$-points is zero for all $l$.
Hence it is enough to prove that $h:=\tr_{H^{-i}_{\sfE,tot}\sfK}$ is the zero-function. 

Finally, to show that a function $h$ on $\Bun_{n,S}(\bF_q)$ is zero it is sufficient to show that
 (1) $h$ is cuspidal and (2) $<h,f> =0$ for all cuspidal functions $f$. In the proof of these statements we will reduce back to a statement for sheaves.

{\em 1$^{st}$ step: $H^{-i}_{\sfE,tot}\sfK$ is a cuspidal complex, therefore $\tr_{H^{-i}_{\sfE,tot} \sfK}$ is a cuspidal function,} i.e. for all $k_1+ k_2=k$ and all $\underline{d_1}+\underline{d_2}=\underline{d}$
let $C_{k_1,k_2}^{\underline{d_1},\underline{d_2}}$ be the functor defined as follows
$$\xymatrix{
& \langle \cE_{k_1}^{\bullet}\hookrightarrow \cE^{\bullet} \to \cE_{k_2}^{\bullet} \rangle \ar[dl]_{\forget} \ar[dr]_{\gr}^{(\cE_{n_i}^\bullet,\cE^\bullet) \mapsto (\cE_{k_1}^\bullet,\cE_{k_2}^\bullet)} & \\
\Bun_{k,S}^{\underline{d}} & & \Bun_{k_1,S}^{\underline{d_1}} \times \Bun_{k_2,S}^{\underline{d_2}}
}$$
\begin{eqnarray*}
C_{k_1,k_2}^{\underline{d_1},\underline{d_2}}: D^b(\Bun_{k,S}^{\underline{d}}) & \to & D^b(\Bun_{k_1,S}^{\underline{d_1}} \times \Bun_{k_2,D}^{\underline{d_2}}) \\
\sfK & \mapsto & C_{k_1,k_2}^{\underline{d_1},\underline{d_2}}\sfK := \bR \gr_!  \forget^* \sfK. 
\end{eqnarray*} 

\begin{definition} A complex $\sfK \in D^b(\Bun_{k,S})$ is called {\em cuspidal} if for all $\un{d_1},\un{d_2}$ and any non trivial partition $k_1+k_2=k$
we have $C_{k_1,k_2}^{\underline{d_1},\underline{d_2}} \sfK = 0$.
\end{definition}

\begin{satz}\label{cuspidal} Let $\sfE$ be a irreducible local system of arbitrary rank $n$ on $C-S$ with indecomposable
unipotent ramification at $S$. Then for all $\un{d}_1,\un{d}_2$ and any non trivial partition $k_1+k_2=k$ the complex $C_{k_1,k_2}^{\un{d}_1,\un{d}_2} \circ H^{-i}_{\sfE,tot} \sfK$ has a filtration with subquotients isomorphic to 
 $(H^{-i_1}_{\sfE,tot}\times H^{-i_2}_{\sfE,tot}) \circ C_{k_1,k_2}^{\un{d}_1,\un{d}_2} \sfK$ for some $i_1+i_2=i$. 
\end{satz}

Note that by induction on $k$ we can assume that the vanishing theorem \ref{verschwinden} holds for all $k_i<k$. Therefore we know that the filtration subquotients are all zero, because $k_1,k_2<k$ and either $i_1$ or $i_2$ is sufficiently big. Therefore the proposition proves that $H_{\sfE,tot}^{-i} \sfK$ is cuspidal if $i>(2g-2)nk+|S|$.

\noindent{\bf Proof:}(of Proposition \ref{cuspidal}) We define a diagram
$$\xymatrix{
\Bun_{k,S}^{\underline{d}}\ar@{}[dr]|{\square} & \langle \cE_{k_1}^\bullet \to \cE_k^{\bullet} \to \cE_{k_2}^{\bullet}\rangle \ar[l]_-{\forget}\ar[r]^{\gr} & \Bun_{k_1,S}^{\underline{d_1}} \times \Bun_{k_2,S}^{\underline{d_2}} \\
\langle \cE_k^{\prime \bullet}\subset \cE_k^{\bullet}\rangle\ar[u]^{\pi_{\text{big}}}\ar[d]_{\pi_{\text{small}}} &  {\underbrace{\langle ( \cE_k^{\prime\bullet}\subset \cE_k^{\bullet},\cE_{k_1}^{\bullet} \to \cE_k^{\bullet} \to \cE_{k_2}^{\bullet})  \rangle}_{=: \text{Middle}}} \ar[l]_-{\forget^\prime} \ar[u]^{\pi_{\text{big}}^\prime}& \\
\Bun_{k,S}^{\underline{d}-\underline{i}} & & 
}$$
to compute \begin{eqnarray*}
C_{k_1,k_2}^{\un{d}_1,\un{d}_2} \circ H^{-\underline{d}}_{\sfE,tot} \sfK & \eqweil{\text{definition}} & \bR \gr_! \forget^* \bR \pi_{\text{big},!} (\pi_{\text{small}}^* \sfK \tensor quot^* \cL_\sfE^i )\\
 & \eqweil{\text{base-change}} & \bR (\gr\circ \pi_{\text{big}}^\prime)_!\underbrace{(\forget^{\prime *}\pi_{\text{small}}^* \sfK\tensor \forget^{\prime *}\quot^* \cL_\sfE^i)}_{\sfK_1}.
\end{eqnarray*}

The stack $\text{Middle}$ is stratified by substacks indexed by $0\leq \underline{i}_1 \leq \underline{d}$, given by the condition $deg(\cE^{\prime (j,p)}\cap \cE_{k_1}^{(j,p)})=d_1^{(j,p)}-i_1^{(j,p)}$:
$$\text{Middle}_{\underline{i}_1}:= \left\langle \begin{array}{l}\xymatrix@C=1ex@R=1ex{
{\cE}_{k_1}^{\prime\bullet} \ar[r]\ar[d]  & {\cE}_{k_{\smash{}}}^{\prime\bullet} \ar[r]\ar[d] & {\cE}_{k_2}^{\prime\bullet} \ar[d] \\
{\cE}_{k_1}^{\bullet} \ar[r]\ar[d]  & {\cE}_{k_{\smash{}}}^{\bullet} \ar[r]\ar[d] & {\cE}_{k_2}^{\bullet} \ar[d]\\
{\cT}_{1}^\bullet \ar[r]  & {\cT}_{\smash{}}^{\bullet} \ar[r]& {\cT}_{2}^{\bullet}
}\end{array} \;\left|\; \cE_{k_1}^{\prime\bullet}=\cE_{k_1}^\bullet \cap \cE_k^{\prime\bullet} \textrm{ and } deg(\cE_{k_1}^{\prime\bullet})=\underline{i}_1 \right.\right\rangle. $$ 

Now $\gr\circ \pi_{\text{big}}^\prime$ restricted to $\text{Middle}_{\underline{i_1}}$ is the map forgetting everything but $\cE_{k_1}^{\bullet}$ and $\cE_{k_2}^{\bullet}$. We factor this as follows:

First consider the map $\forget_{\cE_{n}^{\bullet}}$ forgetting $\cE_k^\bullet$: This is an affine fibration,
 the fibres being homogeneous spaces for $\Ext^1_{para}(\cT_2^{\bullet}, \cE_{k_1}^{\prime\bullet})$ (because of the exact square of $\Ext^1$ groups we get from the extensions of the $\cT^\bullet_i$ by the $\cE^{\prime\bullet}_i$). 

Furthermore, both the map $\pi_{\text{small}}\circ\forget^\prime$ and $quot\circ\forget^\prime$ factor through $\forget_{\cE_{n}^\bullet}$, i.e. $\sfK_1|_{Middle_{\underline{i_1}}}=\forget_{\cE_k^\bullet}^* \sfK_2$ and therefore $\bR forget_{\cE_k^\bullet,!}\sfK_1= \sfK_2 [2c](c)$ for some $c$.

Now we can compose the map $\forget_{\cE^\bullet_k}$ with the forgetful map $\forget_{\cT^{\bullet}}$.  This is just the pull back of the corresponding map in the Hecke correspondence of torsion sheaves, and still $\pi_{\text{small}}\circ\forget^\prime$ factors through this map. Therefore by the Hecke property of $\cL_\sfE$ we get that $\bR \forget_{\cT^\bullet,!} \sfK_2$ is zero if $\underline{i_1}$ is not constant.

But if $\underline{i}_1=(i_1)$ is constant, we get that 
$$\bR (gr\circ\pi_{\text{big}}^\prime)_! (\sfK|_{\text{Middle}_{\underline{i_1}}})= H^{-i_1}_{\sfE,tot}\times H^{-(i-i_1)}_{\sfE,tot} \circ C_{k_1,k_2}^{d_1,d_2} (\sfK).$$
Thus the stratification of the stack $\text{Middle}$  induces a filtration as claimed.

\hfill $\square_{\textrm{\tiny Proposition \ref{cuspidal}}}$

{\em 2$^{nd}$ step: For every cuspidal function $f$ we have $<t_{H_{\sfE,tot}^{-i}\sfK }, f>=0$.}

Define $H_{\sfE,tot}^{i}\sfK:=\bR \pi_{\text{small},!}(\pi_{\text{big}}^* \sfK \tensor quot^* \cL_\sfE^i)$, and denote the analogous operator for functions on $\Bun_{k,S}^{\underline{d}}$ by the same symbol. Then for any cuspidal function $f$ 
$$ <\tr_{H_{\sfE,tot}^{-i}\sfK }, f> = <\tr_\sfK,H^{i}_{\sfE,tot}f>.$$

We want to show that $H_{\sfE,tot}^{i}f=0$ for all cuspidal functions $f$. Using the Langlands
 correspondence for $k<n$, we know that the space of cuspidal functions on $\Bun_{k,S}$ is spanned by cuspidal 
 Hecke eigenfunctions $f_{\sfE^\prime}$ corresponding to local systems $\sfE^\prime$ of dimension
 $k$ with at most unipotent ramification at $S$ and their images under the action of the Iwahori-Hecke algebra (note that for unramified local systems $\sfE^\prime$ on $C$ these functions do not have an eigenfunction property for the Iwahori-Hecke algebra). Furthermore, since $k<n$, we know already that these $f_{\sfE^\prime}$ are the traces of irreducible perverse sheaves $A_{\sfE^\prime}$ on $\Bun_{k,S^\prime}$ for some $S^\prime\subset S$. For this argument we need that $\mathbf{n\leq 3}$ , because for $k\geq 3$ we have not given a construction for representations with reducible unipotent monodromy.

To prove the $2^{\text{nd}}$ step it is therefore sufficient to show:\begin{enumerate}
\item  For all irreducible local systems on $C-S^\prime$ with indecomposable unipotent ramification at $S^\prime \subset S$ we have
$$H_{\sfE,tot}^{i}pr_{\Bun_{k,S^\prime}}^*A_{\sfE^\prime}=0,$$
where $pr_{\Bun_{k,S^\prime}}: \Bun_{k,S} \to \Bun_{k,S^\prime}$ is the map forgetting the parabolic structure at $S-S^\prime$.
\item Any element of the Iwahori-Hecke algebra commutes with the operator $H_{\sfE,tot}^i$. 
\end{enumerate}

We need another Hecke-operator $H_{\sfE,C}^{i}$. As before set $\un{i}:=(i)$.
$$\xymatrix{
& \langle \cE^{\prime\bullet}\subset \cE^{\bullet} \; | \; \cE^{\prime\bullet}\in \Bun_{k,S}^{\underline{d}-\un{i}} \; , \; \cE^{\bullet}\in \Bun_{k,S}^{\underline{d}}  \rangle \ar[dr]^-{\pi_{\text{small}}\times supp}\ar[dl]_-{\pi_{\text{big}}}\ar[d]_-{quot}& \\
\Bun_{k,S}^{\underline{d}} & \Coh_{0,S}^{\un{i}} & \Bun_{k,S}^{\underline{d}-\un{i}} \times C^{(i)} 
}$$
Here $supp(\cE^{\prime\bullet}\subset\cE^\bullet) := supp (\cE^\bullet/\cE^{\prime\bullet})$. We set
\begin{eqnarray*}
H^{i}_{\sfE,C}: D^b(\Bun_{k,S}^{\underline{d}}) & \to & D^b(\Bun_{k,S}^{\underline{d}-\un{i}}\times C^{(i)}) \\
\sfK & \mapsto & H^{i}_{\sfE,C}\sfK := \bR(\pi_{\text{small}}\times supp)! (\pi_{\text{big}}^*\sfK \tensor quot^* \cL^{i}_\sfE).
\end{eqnarray*}
Note that in the above we may assume that we are concerned with $k-$step parabolic structures since the image of $\quot$ is contained in the image of ($k-$step parabolic sheaves) $\subset$ ($n-$step parabolic sheaves).
Thus to prove the first claim we have to show:

\begin{satz}\label{HiEigenschaft} Let $\sfE^\prime$ a local system of rank $k<n$, possibly with unipotent ramification at $S^\prime\subset S$, and let $A_{\sfE^\prime}$ be a Hecke eigensheaf for $\sfE^\prime$ on $Bun_{k,S^\prime}$. Then
$$ H_{\sfE,tot}^{i} \pr_{\Bun_{k,S^\prime}}^* A_{\sfE^\prime} = 0 \;\; for \; i>kn(2g-2)+|S|. $$ 
More precisely, $$ H_{\sfE,C}^{i} \pr_{\Bun_{k,S^\prime}}^* A_{\sfE^\prime} = (j_*(\sfE\tensor\sfE^\prime))^{(i)} \boxtimes \pr_{\Bun_{k,S^\prime}}^* A_{\sfE^\prime} \text{ for all } i.$$
\end{satz}

\noindent{\bf Proof:} The first statement follows from the second, as in the proof of Deligne's Lemma in \cite{Drinfeld_fundamental_group}: 

$H^0(C,j_*(\sfE\tensor\sfE^\prime))=0$, because $\sfE$ is irreducible and not isomorphic to any subquotient of $\sfE^\prime$. By Poincar\'e duality
therefore $H^2(C,j_*(\sfE\tensor\sfE^\prime))=0$ and thus $\dim(H^1(C,j_*(\sfE\tensor\sfE^\prime)))=-\chi(j_*(\sfE\tensor\sfE^\prime))=kn(2g-2)+|S|$ by the formula for the Euler characteristic of Grothendieck-Ogg-Shafarevich (\cite{SGA5} Exp. X, 7.1).

Furthermore we can apply the symmetric K\"unneth formula (\cite{SGA4III} Exp. XVII, 5.5.21) and --- because $h^0=h^2=0$ --- we get that 
$$H^*(C,(j_*(\sfE\tensor\sfE^\prime))^{(i)})=\wedge^i H^1(C,j_*(\sfE\tensor\sfE^\prime)) =0 \text{ for } i>kn(2g-2)+|S|.$$
We are left with proving the second statement. 

{\em Reduction to the case that $i=1$:} Consider the resolution
$$\xymatrix@C=1ex{
                          & \left\langle { \begin{array}{l} \cE^{\prime\bullet}\subset\cE^\bullet \\ \cT_1^\bullet \subset \dots\subset\cT_i^\bullet=\cE^\bullet/\cE^{\prime\bullet}\end{array}} \right\rangle \ar[dr]^{\pi_{\text{small}}\times \widetilde{quot}}\ar[d]^{flag^\prime}\ar@/^1pc/[drr]^{\pi_{\text{small}}\times\widetilde{\supp}} & & \\
                          & \langle \cE^{\prime\bullet}\subset\cE^\bullet \rangle \ar[dl]_{\pi_{\text{big}}}\ar[dr]^{\pi_{\text{small}}\times \quot} & \Bun_{k,S}^{\underline{d}-\underline{i}} \times \widetilde{\Coh}_{0,S}^{\underline{i}}\ar[d]^{\text{flag}}\ar[r] &  \Bun_{k,S}^{\underline{d}-\underline{i}} \times C^i\ar[d]^{\text{sym}} \\   
\Bun_{k,S}^{\underline{d}} & & \Bun_{k,S}^{\underline{d}-\underline{i}} \times \Coh_{0,S}^{\underline{i}} \ar[r] & \Bun_{k,S}^{\underline{d}-\underline{i}} \times C^{(i)}.
}$$
Note that $(H_{\sfE,C}^{1})^{\circ i} \sfK = \bR \pi\times\widetilde{\supp}_! ((\text{flag}^\prime\circ \pi_{\text{big}})^*\sfK \tensor \widetilde{\quot}^* gr^* (\cL_{\sfE}^1)^{\boxtimes i})$. Further, by Lemma \ref{symmetrische_gruppe_operiert} the sheaf 
$\bR \text{flag}^\prime_!(\widetilde{\quot}^*gr^* (\cL_{\sfE}^1)^{\boxtimes i})$ carries an $S_i$-action 
and $$(\bR \text{flag}^\prime_!(\widetilde{quot}^*gr^* (\cL_{\sfE}^1)^{\boxtimes i}))^{S_i}=quot^* \cL_{\sfE}^i.$$ 
Therefore the projection formula implies that the complex $$\bR \text{flag}^\prime_!((\pi_{\text{big}}\circ \text{flag}^\prime)^* A_{\sfE} \tensor \widetilde{\quot}^*gr^* (\cL_{\sfE}^1)^{\boxtimes i})=(H_{\sfE,C}^{1})^{\circ i} A_{\sfE} $$ carries an 
$S_i$ action as well and that $((H_{\sfE,C}^{1})^{\circ i} A_\sfE)^{S_i}=H^{i}_{\sfE,C} A_\sfE$. 

We are reduced to prove that $H^{1}_{\sfE,C} A_{\sfE^\prime} = j_*(\sfE\tensor\sfE^\prime) \boxtimes A_{\sfE^\prime}$.
 
\begin{lemma}\label{H1Eigenschaft} With the notation of Proposition \ref{HiEigenschaft} we have
$$ H^{1}_{\sfE,C}(pr_{\Bun_{k,S^\prime}^{\un{d}}}^* A_{\sfE^\prime}) = pr_{\Bun_{k,S^\prime}^{\un{d}-\un{1}}}^* A_{\sfE^\prime} \boxtimes j_*(\sfE \tensor \sfE^\prime). $$
\end{lemma}

\noindent{\bf Proof:} In the proof we will denote sheaves with parabolic structure at $S$ by $\cE^{\bullet_S}$, and sheaves with 
parabolic structure at $S^\prime$ will be denoted $\cE^{\bullet_{S^\prime}}$ to distinguish the two 
sets of data.
We have
$$\xymatrix{
& \langle \cE^{\prime\bullet_S} \subset \cE^{\bullet_S} \rangle \ar[dl]_{\pi_{big}} \ar[d]^{\forget_{\Hecke}} \ar[dr]^{\pi_{small}\times pr_C} & \\
\Bun_{k,S}^{\un{d}} \ar[d]_{pr_{\Bun_{k,S^\prime}^{\un{d}}}} & \langle \cE^{\prime\bullet_{S^\prime}} \subset \cE^{\bullet_{S^\prime}} \rangle \ar[dl]_{\pi^\prime_{big}}\ar[dr] & \Bun_{k,S}^{\un{d}-\un{1}}\times C \ar[d] \\
\Bun_{k,S^\prime}^{\un{d}} & & \Bun_{k,S^\prime}^{\un{d}-\un{1}}\times C
.} $$
This induces a map
$$ pr_{\text{Fib}}: \langle \cE^{\prime \bullet_S} \subset \cE^{\bullet_S} \rangle \to 
\langle \cE^{\prime\bullet_{S^\prime}} \subset \cE^{\bullet_{S^\prime}} \rangle \times_{\Bun_{k,S^\prime}^{\un{d}-\un{1}}\times C}(\Bun_{k,S}^{\un{d}-\un{1}}\times C)=: \text{Fib}.$$

Denote by $pr_1,pr_2$ the projections from this fibre product to its factors, and let $quot,quot^\prime$
be the quotient maps from the Hecke correspondence to $\Coh_{0,S}^{\un{1}}$ and $\Coh_{0,S^\prime}^{\un{1}}$ respectively.

We can apply the projection formula to rewrite
$$ H^1_{\sfE,C}(pr_{\Bun_{k,S^\prime}^{\un{d}}}^* A_{\sfE^\prime}) = \bR pr_{2,!} ((\bR pr_{Fib,!} quot^* \cL_{\sfE})\tensor (\pi_{big}^\prime\circ pr_1)^* A_{\sfE^\prime} ).$$
The calculation of $\bR pr_{Fib,!} quot^* \cL_\sfE$ can be reduced to a 
calculation for torsion sheaves as follows. We have a map:
\begin{eqnarray*}
 Fib & \map{q} & \Coh_{0,S}^{\un{e}_n} \\
 (\cE^{\prime \bullet_S},\cE^{\bullet_{S^\prime}}) & \mapsto & 
                \cT^{(i,p)}:= \left\{ \begin{array}{ll}
                                        \cE^{(i,p)}/\cE^{\prime(i,p)} & \textrm{if } p\in S \textrm{ or } i=0\\
                                        \cE^{(0,p)}(p)/\cE^{\prime(i,p)} & \textrm{if } p\in S-S^\prime \textrm{ and } i\neq 0,
                              \end{array}\right. 
\end{eqnarray*}
where $\un{e}_n^{*,p}=\left\{\begin{array}{ll} (1,\dots,1) & \text{if } p \in S^\prime \\
                                             (1,n-1,\dots,2) & \text{if } p\in S-S^\prime \end{array}\right.$.
This gives rise to the cartesian diagram
$$\xymatrix{
\langle \cE^{\prime\bullet}\subset\cE^\bullet \rangle \ar[d]^{pr_{Fib}}\ar[r]^-{\tilde{q}}\ar@{}[dr]|{\square} & \langle (\cT^{\prime\bullet}\subset\cT^\bullet) | {\cT^{\prime\bullet} \in \Coh_{0,S}^{\un{1}}\hfill \atop \cT^\bullet \in \Coh_{0,S}^{\un{e}_n}} \rangle\ar[d]^{\forget_{\cT^{\prime\bullet}}}\ar[r]^-{pr_{\cT^{\prime\bullet}}} & \Coh_{0,S}^{\un{1}}\\
\text{Fib} \ar[r]^-{q}      & \Coh_{0,S}^{\un{e´}_n}\ar[r]^-{\forget_{S-S^\prime}} & \Coh_{0,S^\prime}^{\un{1}}, 
}$$
where $\tilde{q}(\cE^{\prime\bullet_S} \subset \cE^{\bullet_S})=\big(\cE^{\bullet_S}/\cE^{\prime\bullet_S}\subset q(pr_{Fib}(\cE^{\prime\bullet_S},\cE^{\bullet_S}))\big)$.
By the base change formula it will be sufficient to calculate:
\begin{lemma}\label{verflixtes}
$(\bR \forget_{\cT^\prime,!} pr_{\cT^\prime}^* \cL_\sfE^1)|_{Im(q)} \cong (\forget_{S-S^\prime}^* \cL_\sfE^1)|_{Im(q)}$, where by abuse of notation
we denoted by $\cL_\sfE^1$ the middle extensions of $\sfE$ on $C-S$ to $\Coh_{0,S}^{\un{1}}$ and $\Coh_{0,S^\prime}^{\un{1}}$.
\end{lemma}
\noindent{\bf Proof:} First note that the image of the map $q$ is the open substack of $\Coh_{0,S}^{\un{e}_n}$ defined by the condition that the maps $\phi^{i,p}$ are surjective for $1<i\leq n$ and $p\in S-S^\prime$. This follows because the image of $q$ is contained in this substack and the structure of torsion sheaves of degree $\un{e}_n$ (Lemma \ref{struktur}) shows that $q$ exhausts this substack.
 
Further, note that the map $\pr_{\cT^{\prime\bullet}}$ is smooth, since it can be factored into a generalized vector bundle over $\Coh_{0,S}^{\un{e}_n-\un{1}}\times \Coh_{0,S}^{\un{1}}$ and the projection onto the second factor.
The map $\forget_{\cT^\prime}$ is projective because the fibres are closed in a product of projective spaces and therefore $\bR\forget_{\cT^{\prime\bullet},*}=\bR\forget_{\cT^{\prime\bullet},!}$.

Denote by $\cT_0^{j,p} := \left\{\begin{array}{cl} 0 & \text{if } j=0 \text{ or } p\in S^\prime \\ k_p^{\oplus n-j} & \text{else} \end{array}\right.$ with $\phi^{i,p}:k_p^{\oplus n-i+1} \tto k_p^{n-i}$ some projection, and denote by $j_{S-S^\prime}:\Coh_{0,C-(S-S^\prime),S^\prime}^{\un{1}} \to \Coh_{0,S}^{\un{e_n}}$ the morphism $\cT^{\prime\bullet} \mapsto \cT^{\prime\bullet} \oplus \cT_0^\bullet$. 

Combining the two remarks above we get a canonical morphism 
$$F: j_{S-S^\prime,*} \cL^1_\sfE \to \bR\forget_*{\cT^{\prime\bullet}}\pr_{\cT^{\prime\bullet}}^* \cL^1_\sfE = \bR\forget_!{\cT^{\prime\bullet}}\pr_{\cT^{\prime\bullet}}^* \cL^1_\sfE,$$ 
and $j_{S-S^\prime,*}\cL^1_\sfE \cong \quot^{\prime,*}\cL_\sfE^1$.
We have to prove that $F$ is an isomorphism over the image of $q$. First note that $\forget_{\cT^{\prime\bullet}}$ is an isomorphism over the open substack where $supp(\cT^{(0,p)})\in C-(S-S^\prime)$, so the above sheaves are isomorphic on this substack. 
 
We are left to check that $F$ is an isomorphism on the fibres over points $\cT^\bullet$ with $\cT^{0,p}=k_p$ and $p\in S-S^\prime$. Since this problem is local on $\Coh_{0,S}^{\un{e}_n}$ we may assume that $(C,S,S^\prime)=(\bA^1,\{ 0\},\emptyset)$ and $\sfE=\sfE_n$ (see Section \ref{A1}). 

We know that $k_p\cong \cT^{\prime(i,p)}\subset\cT^{(i,p)}$, and we may factorize $\forget_{\cT^{\prime\bullet}}$, into the maps forgetting the choice of the subspaces $\cT^{\prime,(i,p)}$ for $i>k$. Consider for example the map forgetting the choice of $\cT^{\prime(n-1,p)}$. Its fibre is either a single point, if $\phi^{n-1,p}(\cT^{\prime(n-2,p)})\neq 0$, or it is isomorphic to the projective space $\bP(H^0(C,\cT^{(n-1,p)}))$, and the kernel of $\phi^{(n,p)}$ defines a linear subspace of codimension $1$. Thus we can apply the calculation of $\cL_{\sfE_n}^1$ (Lemma \ref{restriction_to_D}) to conclude that the cohomology of this fibre is isomorphic to the fibre of $\cL_\sfE^1$ at $\cT^{\prime\bullet}$ for any choice of $\cT^{\prime\bullet}$ not contained in this subspace. By induction we therefore get the claimed isomorphism.    
\hfill $\square_{\text{Lemma \ref{verflixtes}}}$


Continuing the proof of Lemma \ref{H1Eigenschaft} we can factor $pr_2$ as
$$ Fib 
 \map{\tilde{pr}_2} \Bun_{2,S}^{\un{d}-\un{1}} \times \Coh_{0,S^\prime}^{\un{1}} \map{id\times pr_C} \Bun_{2,S}^{\un{d}-\un{1}} \times C$$
and apply the projection formula again:
\begin{eqnarray*}
H^1_{\sfE,C} (pr_{\Bun_{k,S^\prime}^{\un{d}}}^* \sfA_{\sfE^\prime})
 & = & \bR (id\times pr_C)_! (quot^{\prime *} \cL_{\sfE}\tensor \bR \tilde{pr}_{2,!} \sfA_{\sfE^\prime})  \\
 & \eqweil{\text{base-change}} & \bR (id\times pr_C)_! (quot^{\prime *} \cL_{\sfE} \tensor (\sfA_{\sfE^\prime} \boxtimes \cL_{\sfE^\prime}))\\
 & \eqweil{\text{proj.fmla}} & \sfA_{\sfE^\prime} \boxtimes (\bR pr_{C,!} \cL_\sfE \tensor \cL_{\sfE^\prime}) \\  
 & \eqweil{\text{Corollary \ref{pr-zwei}}} & \sfA_{\sfE^\prime} \boxtimes (j_* (\sfE \tensor \sfE^\prime)).
\end{eqnarray*}\hfill $\square_{\text{Lemma }\ref{H1Eigenschaft} \text{ and Proposition } \ref{HiEigenschaft}}$

To finish the proof of the vanishing theorem \ref{verschwinden} we have to show that the operator $H^{i}_{\sfE,tot}$ commutes with all other Hecke operators at the level of functions. We may apply the reduction of Proposition \ref{HiEigenschaft} to reduce ourselves to prove this for the operator $H^1_{\sfE,C}$.

Fix a parabolic torsion sheaf $\cT^\bullet$ and define the Hecke operator 
$H_{\overline{\cT}^\bullet}$ as the sum over all Hecke operators corresponding to torsion sheaves 
contained in the closure of $\cT^\bullet$.

Let $\overline{\langle \cT^\bullet\rangle }\subset \Coh_{0,S}^{deg(\cT^\bullet)}$ be the closure
of the substack classifying parabolic torsion sheaves which are locally isomorphic to $\cT^\bullet$.
And define the stack
$$ \Hecke_{\overline{\cT}^\bullet}:=\langle \cE^{\prime\bullet}\subset\cE^\bullet | \cE^\bullet/\cE^{\prime\bullet} \in \overline{\langle\cT^\bullet\rangle}\subset \Coh_{0,S}^{\deg(\cT^\bullet)}\rangle, $$
which as before provides a Hecke operator
$$ H_{\overline{\cT^\bullet}} : D^b(\Bun_{n,S}^{\un{d}}) \to D^b(\Bun_{n,S}^{\un{d}-deg(\cT^\bullet)}).$$

By induction on the codimension of $\langle \cT^\bullet \rangle \subset \Coh_{0,S}^{\deg(\cT^\bullet)}$ it
is sufficient to prove that $H^1_{\sfE,C}$ commutes with $H_{\overline{\cT^\bullet}}$ for all $\cT^\bullet$.

\begin{lemma}\label{kommutiert} For any $\sfK\in D^b(\Bun_{n,S}^{\un{d}})$ we have
$$H^{1}_{\sfE,C} \circ H_{\overline{\cT^\bullet}} \sfK \cong H_{\overline{\cT^\bullet}} \circ H^{1}_{\sfE,C} \sfK$$
in $D^b(\Bun_{n,S}^{\un{d}-\un{1}-deg(\cT^\bullet)} \times C)$.
\end{lemma} 

\noindent{\bf Proof:} We may assume that $supp(\cT^\bullet)=p$ for a single point $p\in S$, since every
torsion sheaf is the direct sum of sheaves supported at a single point and for $p\not\in S$ 
the lemma is easy.

The claim follows easily from the corresponding lemma for parabolic torsion sheaves (the reduction
using the projection formula once again works as in the previous lemma, we skip it):
Denote 
$$\text{Flag}_{1,\cT^\bullet}:= \left\langle (0 \to \cT^{\prime\bullet} \to \cQ^\bullet \to \cT^{\pprime\bullet}\to 0)
\left| \begin{array}{l} \cT^{\prime\bullet} \in \Coh_{0,S}^{\un{1}}\\
                        \cQ^\bullet \in \Coh_{0,S}^{\un{1}+\deg(\cT^\bullet)} \\
                        \cT^{\pprime\bullet} \in \overline{\langle\cT^\bullet\rangle}
       \end{array}\right.\right\rangle.$$
Denote by $pr_{\cT^{\prime\bullet}}$,$pr_{\cQ^\bullet}$ the projections and by $pr_C$ the projection to
the curve $C$ defined by the support of $\cT^{\prime\bullet}$.     

Let $\text{Flag}_{\cT^\bullet,1}$ be the stack defined as above with the roles of
$\cT^{\prime\bullet}$ and $\cT^{\pprime\bullet}$ interchanged, i.e. $\cT^{\pprime\bullet}\in \Coh_{0,S}^{\un{1}}$ and $\cT^{\prime\bullet} \in \overline{\langle\cT^\bullet\rangle}$,
and denote its projections by $rp_{\cQ^\bullet}$ etc.

\begin{lemma} We have a canonical isomorphism of sheaves
$$\bR (pr_{\cQ^\bullet}\times pr_C)_! pr_{\cT^{\prime\bullet}}^* \cL_\sfE^1 \cong  \bR (rp_{\cQ^\bullet}\times rp_C)_! rp_{\cT^{\pprime\bullet}}^* \cL_\sfE^1$$
on $\Coh_{0,S}^{\un{1}+deg(\cT)} \times C$.
\end{lemma} 
\noindent{\bf Proof:} This is similar to the proof of Lemma \ref{verflixtes}:
Over the open substack of $\Coh_{0,S}^{\un{1}+ deg(\cT^\bullet)}$ where the support
of the torsion sheaf is not equal to $supp(\cT^\bullet)=p$ the stacks $Flag_{1,\cT^\bullet}$ and
$Flag_{\cT^\bullet,1}$ are isomorphic, because there are no extensions between sheaves supported
at different points. Therefore the claimed isomorphism exists over this subset.
To extend it, we again reduce to the case $(C,S)=(\bA^1,\{0\})$ and note that
the maps $pr_{\cQ^\bullet},rp_{\cQ^\bullet}$ are projective and the map $pr_{\cT^{\prime\bullet}}$
(resp. $rp_{\cT^{\pprime\bullet}}$) can be factored as
$$Flag_{1,\cT^\bullet} \to \Coh_{0,S}^{\un{1}} \times \overline{\langle \cT^\bullet\rangle} \to \Coh_{0,S}^{\un{1}}.$$
The first map is a generalized vector bundle, and the second one is the projection of a product, therefore both maps are locally acyclic. 
Hence we can use the exact triangle
$$ \to \cL_\sfE^1 \to \bR j_* \sfE_\infty \to j_! \sfE_{\infty}(-n) \map{[1]}$$
of Proposition \ref{e_unendlich}  once more. If we replace $\cL_\sfE^1$ by one of the 
sheaves $\bR j_*\sfE_\infty$ or $j_! \sfE_\infty(-n)$ of Proposition \ref{e_unendlich} the lemma
follows from the Leray spectral sequence. 
Therefore the lemma follows for $\cL_\sfE^1$ as well. \hfill $\square$


\section{The vanishing theorem implies that $j_{\Hom,!} \sfF^k_\sfE = j_{\Hom,!*}\sfF^k_\sfE = \bR j_{\Hom,*} \sfF^k_\sfE$}\label{section_sauberkeit}

With the notations of the fundamental diagram (\ref{Teil1}) of Section 2 we have

\begin{satz}\label{sauberkeit}
Assume that the vanishing theorem \ref{verschwinden} holds for $k<n$. 
Then for $k<n$ we have $j_{\Hom,!} \sfF^k_\sfE = j_{\Hom,!*} \sfF^k_\sfE$ and thus for $k\leq n$ we have $\sfF^k_{\sfE,!}=\sfF^k_\sfE$.
\end{satz}
Since we have shown the vanishing theorem for local systems of rank $\leq 3$, we get in particular:
\begin{kor}\label{sauber}
For $k\leq n\leq 3$ the sheaves $\sfF_\sfE^k\cong \sfF_{\sfE,!}^k$ are isomorphic. 

\hfill $\square_{\text{\em Corollary}}$
\end{kor}

\noindent{\bf Proof of Proposition \ref{sauberkeit}:} The Hecke-property of $\cL_\sfE^d$ allow us to copy the proof in \cite{FGV3} with some minor changes. We use induction, and assume that the proposition is true for all $k^\prime < k$.

{\em 1. Step:} The claim is true over the substack of parabolic vector bundles.

Here every nontrivial homomorphism is injective, that is 
$$\Hom^{\text{inj}}_k =\Hom_k-(\text{zero-section}) \text{ over } \Bun_{k,S}^{\text{good}}.$$ 
Furthermore $\sfF^k_\sfE$ is $\bG_m$-invariant, since the Fourier transform preserves this 
property by \cite{Laumon_transformation}, Proposition 1.2.3.4. Therefore we can apply Lemma \ref{brylinski-lemma} and get $$j_{\Hom,!} \sfF^k_\sfE = \bR j_{\Hom,*}\sfF^k_\sfE=j_{\Hom,!*}\sfF_\sfE^k  \Leftrightarrow \bR \pi_! \sfF^k_\sfE =0,$$ where $\pi : \Hom^{\text{inj}}_k \to \Bun_{k,S}^{d,\text{good}}$ is the projection.
 
As in Section 5 we can calculate $\sfF^k_\sfE=\sfF^{k}_{\sfE,!}$ in a different way:
{\small
$$\xymatrix@C=2ex{
 & {\left\langle
        \begin{array}{l}
           \cE_1^\bullet\subset \dots\subset \cE_k^\bullet\subset\cE^\bullet \\
           \cE_i^\bullet/\cE_{i-1}^\bullet \map{\cong} (\Omega^{\bullet,k-i})
        \end{array}\right\rangle}  \ar[d]\ar[dr]\ar@/_1pc/[ddl]_{\tilde{\pi}} \ar@/^1pc/[drr]^{\widetilde{ext}} & & \\
 & \langle \cE_k^\bullet \hookrightarrow \cE^\bullet \tto \cT^\bullet \rangle \ar[dl] \ar[dr] & 
  {\left\langle
     \begin{array}{l} 
       \cE_1^\bullet\subset\cE_2^\bullet\subset \dots\subset \cE_k^\bullet\\
       \cE_i^\bullet/\cE_{i-1}^\bullet \map{\cong} (\Omega^{\bullet,k-i})
      \end{array}
   \right\rangle} \times \Coh_{0,S}^{\underline{d}}\ar[d]^{\pi^\prime}\ar[r]_-{ext} & \bA^1\\
\Bun_{k,S}^d & & \Bun_{k,S}^0 \times \Coh_{0,S}^{\underline{d}}. &
}$$
}
Therefore $\bR \pi_! \sfF^k_\sfE=\bR \tilde{\pi}_! (\widetilde{ext}^*\sfL_\psi \tensor quot^*\cL_\sfE) = H^{-d}_{\sfE,tot} (\bR\pi^\prime_! ext^* \sfL_\psi)$, and the vanishing theorem \ref{verschwinden} implies that $H^{-d}_{\sfE,tot} (\bR\pi^\prime_!ext^*\sfL_\psi)  = 0.$ 

{\em 2. Step:} Induction on the length of the torsion of $\cF^\bullet$:

For any $\underline{r}=(r_{i,p})$ denote by
$$ \Coh_{k,S}^{\underline{d},\leq \underline{r}} := \langle \cF^{\bullet} \in \Coh_{k,S}^\bullet | \length(\torsion(\cF^\bullet))\leq \underline{r}\rangle $$
the stack of parabolic sheaves such that the length of the torsion of the coherent sheaves $\cF^{(i,p)}$ is bounded by $r_{i,p}$.

It is sufficient to prove the proposition after a smooth base change. To get a map to torsion free parabolic sheaves (we want to apply the vanishing theorem again) we use the stack
$$ \widetilde{\Coh}_{k,S}^{\underline{d},\leq \underline{r}}:= \langle \cE^\bullet\subset \cF^\bullet \to \cT^\bullet |\cE^\bullet \in \Bun_{k,S}^{\underline{d}-\underline{r}},\cF^\bullet\in \Coh_{k,S}^{\underline{d},\text{good}} \cT^{\bullet} \in \Coh_{0,S}^{\underline{r}} \rangle. $$
 
From Lemma \ref{hecke_abbildungs_eigenschaft} we know that the forgetful map $\widetilde{\Coh}_{k,S}^{\underline{d},\leq \underline{r}} \to \Coh_{k,S}^{\underline{d},\leq \underline{r}}$ is smooth. And the map 
\begin{eqnarray*}
\gr: \widetilde{\Coh}_{k,S}^{\underline{d},\leq \underline{r}} &\to& \Bun_{k,S}^{\underline{d}-\underline{r}}\times \Coh_{0,S}^{\underline{r}}\\
\cE^\bullet\subset\cF^\bullet & \mapsto & (\cE^\bullet,\cF^\bullet/\cE^\bullet)
\end{eqnarray*}
is a vector bundle, since $\dim(\Ext_{\text{para}}(\cT^\bullet,\cE^\bullet))$ is constant by Lemma (\ref{dim_Hom_und_Ext}).

Furthermore, over any point of $\widetilde{\Coh}^{\un{d}\leq \un{r}}_{k,S}$ we have $\Ext^1(\Omega^{\bullet,k-1},\cE^\bullet)=0$ (by assumption $\cF^\bullet\in \Coh_{k,S}^{\underline{d},\text{good}}$),
 therefore the dimension of $\Hom(\Omega^{\bullet,k-1},\cE^\bullet)$ is constant, so $\Hom(\Omega^{\bullet,k-1},\cE^\bullet)$ is a vector bundle over this stack.

Consider the base change $\widetilde{\Hom}_k$  of $\Hom_k$ to $\widetilde{\Coh}_{k,S}^{\underline{d},\leq \underline{r}}$, and 
analogously define 
$$\widetilde{\Hom}^{\text{inj}}_k:=\Hom^{\text{inj}}_k\times_{{\Coh}_{k,S}^{\underline{d}}} \widetilde{\Coh}_{k,S}^{\underline{d},\leq \underline{r}}.$$

By the above, the map
\begin{eqnarray*}
 gr_{\widetilde{\Hom}}:\widetilde{\Hom}_k& \to & \Bun_{0,S}^{\un{d}-\un{r}} \times \Hom(\Omega^{\bullet,k-1},\cT^\bullet) \\
(\Omega^{\bullet,k-1} \map{s} \cF^\bullet,\cE^\bullet\subset\cF^\bullet\map{p}\cT^\bullet) & \mapsto & (\cE^\bullet, \Omega^{\bullet,k-1} \map{p\circ s} \cT^\bullet)
\end{eqnarray*}
is also vector bundle, because it is the composition of the map induced by composing $p\circ s$, which has fibres $\Hom(\Omega^{\bullet,k},\cE^\bullet)$,  and $\quot$. The zero section 
$$(\cE^\bullet,\Omega^{\bullet,k-1}\map{s} \cT) \mapsto (\cT\map{(0,s)} \cE\oplus\cT, \cE\subset \cE\oplus\cT \to \cT)$$
of this bundle is the substack \footnote{Note that, if there is a splitting of  $\cF^\bullet\tto\cT^\bullet$, then there is a unique one, since $\torsion(\cF^\bullet)$ is a subsheaf of $\cF^\bullet$.}   
$$ \left\langle\left( {\begin{array}{l}\Omega^{\bullet,k-1} \map{s} \cF^\bullet \\ \cE^\bullet\subset\cF^\bullet\map{p}\cT^\bullet \end{array}} \right) \; |\; \length(\torsion(\cF^\bullet))=\un{r} \textrm{ and } \Omega^{\bullet,k-1}\to\cT^\bullet\to\cF^\bullet \right \rangle$$
and this is by induction hypothesis the substack to which we have to extend $\sfF^k_\sfE$.

Thus, denote $gr_{\widetilde{\Hom}^{\text{inj}}}:=gr_{\widetilde{\Hom}}|_{\Hom^{\text{inj}}}$ and again we have to show that $$\bR gr_{\widetilde{\Hom}^{\text{inj}},!} \sfF^k_\sfE = 0.$$ 
Since this can be checked fibrewise, we fix a point $(\cE^\bullet,\Omega^{\bullet,k-1}\to \cT^\bullet)\in \Bun_{k,S}^{\underline{d}-\underline{r}}$ and denote by $\text{Fibre}_{\cE^\bullet,\Omega^{\bullet,k-1}\to\cT^\bullet}$ the fibre of $\widetilde{\Hom}^{\text{inj}}_k$ over this point.

{\em Step 2.1} Reduction to the case that $\Omega^{\bullet,k-1} \tto \cT^{\bullet}$ is surjective.

Factor $s: \Omega^{\bullet,k-1} \to Im(s)=: \cT^{\prime\bullet} \subset \cT^\bullet$, denote $\cT^\bullet /\cT^{\prime\bullet} =:\cT^{\pprime\bullet}$ and set $\un{r}^\pprime:=deg(\cT^{\pprime\bullet})$. Then for any $(\Omega^{\bullet,k-1}\to\cF^\bullet\tto\cT^\bullet)\in \Fibre_{\cE^\bullet,\Omega^{\bullet,k-1}\map{s}\cT^\bullet}$ we get an extension $\cF^{\prime\bullet} \subset \cF^\bullet \tto \cT^{\pprime\bullet}$. Consider the Hecke operator for
$$\xymatrix@C=0.5ex{
 & \Hecke_{\Hom^{\text{inj}}}:=\left\langle{\begin{array}{l} \Omega^{\bullet,k-1} \hookrightarrow \cF^{\prime\bullet} \\ \cF^{\prime\bullet} \hookrightarrow \cF^\bullet \tto \cT^{\pprime\bullet}\end{array}} \right\rangle \ar[dl]\ar[dr] & \\
\Hom^{\text{inj}}_k & & \Hom^{\text{inj}}_k \times \Coh_{0,S}^{\underline{r}^\pprime}.
}$$
We know by Proposition \ref{hecke_eigenschaft} that $\sfF_{\sfE!}^k$ is a Hecke eigensheaf
and that $$\Fibre_{\cE^\bullet,\Omega^{\bullet,k-1}\map{s}\cT^\bullet}= \Fibre_{\cE^\bullet,\Omega^{\bullet,k-1}\to\cT^{\prime\bullet}}\times_{\Hom^{\text{inj}}_k\times \Coh_{0,S}^{\underline{r}^\pprime}}\Hecke_{\Hom^{\text{inj}}}.$$ 

Thus, in case that $\un{r}^\pprime$ is not constant,  we can establish our claim that $$H_c^*(\Fibre_{\cE^\bullet,\Omega^{\bullet,k-1}\map{s}\cT},\sfF_\sfE^k)=0,$$ since the above Hecke operator is zero by \ref{hecke_und_konstruktion}. 

If on the other hand $\un{r}^\pprime$ is constant, we know that it is sufficient to prove the claim for $\Fibre_{\cE^\bullet,\Omega^{\bullet,k-1}\to\cT^\prime}$. This has already been done in the case that $\cT^{\prime\bullet}\neq \cT^\bullet$. Therefore we may assume that $Im(\Omega^{\bullet,k-1})=\cT^{\prime\bullet}=\cT^\bullet$.

{\em Step 2.2.} Assume that $\Omega^{\bullet,k-1} \tto \cT^\bullet$ is surjective, i.e. $\cT^\bullet\cong \Omega^{\bullet,k-1}/\Omega^{\bullet,k-1}(-D)$ for some effective parabolic divisor $D$. 

In this case, giving an element $(\Omega^{\bullet,k-1}\hookrightarrow\cF^\bullet)\in \Fibre_{\cE^\bullet,\Omega^{\bullet,k-1}\to\cT^\bullet}$ is the same as 
to give a map $\Omega^{{k-1}\bullet}(-D)\to \cE^\bullet$:
$$ (\Omega^{\bullet,k-1}(-D) \to \cE^\bullet,\cT^\bullet) \mapsto ( \cE^\bullet \subset (\cE^\bullet\oplus \Omega^{\bullet,k-1})/\Omega^{\bullet,k-1}(-D)) $$ 

And indeed for any square
$$\xymatrix{ 
  {\Omega^{\bullet,k-1}(-D)} \ar@{^(->}[r] \ar@{^(->}[d] &  \Omega^{\bullet,k-1}\ar@{->>}[r] \ar@{^(->}[d] & {\cT}^\bullet \ar[d]^{\cong} \\
{\cE^\bullet} \ar@{^(->}[r] & {\cF^\bullet} \ar@{->>}[r] & {\cT^\bullet} }$$
we automatically get that $\cF^\bullet= \cE^\bullet \oplus_{\Omega^{\bullet,{k-1}}(-D)} \Omega^{\bullet,{k-1}}$.

Thus we get an isomorphism $\Fibre_{\cE^\bullet,\Omega^{\bullet,k-1}\to\cT^\bullet}\map{\cong}\Hom(\Omega^{\bullet,k-1}(-D),\cE^\bullet)$. Furthermore under this isomorphism $\sfF^k_\sfE|_{\Fibre}$ becomes the sheaf on $\Bun_{k,S}^{\underline{d}-\underline{r}}$ constructed in the same way as $\sfF^k_\sfE$, by replacing $\Omega^{\bullet,k-1}$ by $\Omega^{\bullet,k-1}(-D)$. More precisely, since $$\cE^\bullet /\Omega^{\bullet,k-1}(-D)\cong \cF^\bullet /\Omega^{\bullet,k-1},$$ we have again:
{\scriptsize
$$\xymatrix@C=0.5ex{
 & {\left\langle
        \begin{array}{l}
          \cE_1^\bullet\subset \dots\subset \cE_n^\bullet\subset\cE^\bullet \\           \cE_1^\bullet \map{\cong} \Omega^{\bullet,k-1}(-D) \\
          \cE_{i+1}^\bullet/\cE_{i}^\bullet \map{\cong} \Omega^{\bullet,k-i+1}
        \end{array}\right\rangle}  \ar[d]\ar[dr]\ar[dl]^{\textstyle\forget}\ar@/_5pc/[ddl]_{\textstyle \tilde{\pi}}\ar@/^1pc/[rr]^{\textstyle\widetilde{\ext}} & & \bA^1 \\
 \Hom^{\text{inj}}(\Omega^\bullet(-D),\cE^\bullet) \ar[d]^{\textstyle \pi} & \langle \cE_k^\bullet \hookrightarrow \cE^\bullet \tto \cT^\bullet \rangle \ar[dl] \ar[dr] & 
  {\left\langle
        \begin{array}{l}
           \cE_1^\bullet\subset \dots\subset \cE_n^\bullet\\
           \cE_1^\bullet \map{\cong} \Omega^{\bullet,k-1}(-D) \\
           \cE_{i+1}^\bullet/\cE_{i}^\bullet \map{\cong} \Omega^{\bullet,k-i+1}
        \end{array}\right\rangle} \times \Coh_{0,S}^{\underline{d}}\ar[d]^{\textstyle\pi^\prime}\ar[ur]_-{\textstyle\ext} & \\
\Bun_{k,S}^{\un{d}-\un{r}} & & \Bun_{k,S}^{-\un{r}} \times \Coh_{0,S}^{\underline{d}} &
}$$
}
Here $\ext$ is the composition

{\small$\xymatrix{
\langle \cE_i^\bullet,gr(\cE_k^\bullet)\cong \Omega^{\bullet,k-1}(-D)\oplus\bigoplus_{j=0}^{k-2} \Omega^{\bullet,j} \rangle \ar[r]\ar[ddr]_-{\ext}& H^1(C,\Omega^{1}(-D)) \oplus \bigoplus_{j=0}^{k-2} H^1(C,\Omega^{1}) \ar[d]\\
& H^1(C,\Omega) \oplus_{j=0}^{k-2} H^1(C,\Omega) \ar[d]^{sum} \\
& H^1(C,\Omega) \cong \bA^1}$}
 
and therefore 
$$\sfF^k_{\sfE}|_{\Fibre_{\cE^\bullet,\Omega^{\bullet,k-1}\tto\cT}}\cong(\bR \forget_! (quot^* \cL_\sfE\tensor \ext^* \sfL_\psi))|_{\textrm{\scriptsize Fibre of }\pi\textrm{\scriptsize\ over }\cE^\bullet}.$$
But here we can apply the vanishing theorem again, because 
$$\bR (\pi\circ \forget)_! ( quot^* \cL_\sfE\tensor \ext^* \sfL_\psi)= H_{\sfE,tot}^{-d} (\bR \pi^\prime_! \ext^*\sfL_\psi) =0.$$
\hfill $\square_{\textrm{\tiny Proposition \ref{sauberkeit}}}$



\end{document}